\documentclass[11pt,reqno,oneside]{amsart}
\usepackage{amscd}
\usepackage{amsmath}
\usepackage{latexsym}
\usepackage{amsfonts}
\usepackage{amssymb}
\usepackage{amsthm}
\usepackage{graphicx}
\usepackage{verbatim}
\usepackage{mathrsfs}
\usepackage{enumerate}
\usepackage{hyperref}
\usepackage{fullpage}
\usepackage{xcolor}
\usepackage[latin1]{inputenc}

\theoremstyle{plain}\newtheorem{definition}{Definition}[section]
\theoremstyle{plain}\newtheorem{theorem}{Theorem}[section]
\theoremstyle{plain}\newtheorem{lemma}[theorem]{Lemma}
\theoremstyle{plain}\newtheorem{coro}[theorem]{Corollary}
\theoremstyle{plain}\newtheorem{proposition}[theorem]{Proposition}
\theoremstyle{plain}\newtheorem{remark}{Remark}[section]

\numberwithin{equation}{section}

\newcommand{\norm}[1]{\left\|#1\right\|}

\newcommand{\R}{\mathbb{R}}
\newcommand{\be}{\begin{equation}}
\newcommand{\ee}{\end{equation}}
 \newcommand{\ba}{\begin{aligned}}
 \newcommand{\ea}{\end{aligned}}

  \newcommand{\ben}{\begin{enumerate}}
   \newcommand{\een}{\end{enumerate}}

\newcommand{\Rmnum}[1]{\expandafter\@slowromancap\romannumeral #1@}

\allowdisplaybreaks

\begin{document}
%%%%%%%%%%%%%%%%%%%%%%%%%%%%%%%%%%%%%%%%%%%%%%%%%%%%%%%%%%%%%%%%%%%%%%%%%%%%%%%%%%%%%%%%%%%%%%%%%%%%
\title{On a multi-dimensional transport equation with nonlocal velocity and fractional dissipation}
\author{Wanwan Zhang$^{1}$}
%\address{$^1$ School of Mathematical and Computational Science, Xiangtan University, Hunan, 411105, P.R.China}
%\email{thu3141@126.com}

%\address{$^2$ Department of Mathematics, Beijing Technology and Business University, Beijing, 100048, P.R.China}
%\email{liuyn@th.btbu.edu.cn}

\address{$^1$ School of Mathematics and Statistics, Jiangxi Normal University, Nanchang, 330022, Jiangxi,
P. R. China}
\email{zhangww@jxnu.edu.cn}

\thanks{$^*$Corresponding author}
\begin{abstract}
This paper aims to investigate a multi-dimensional transport equation with nonlocal velocity and fractional dissipation. The balance between the nonlinearity and dissipation gives rise to three different cases, namely the subcritical, critical and supercritical ranges. We study those three cases and  obtain a set of results containing local well-posedness, global smoothness, eventual regularity and finite-time blowup of smooth solutions.
\end{abstract}
\smallskip
\maketitle
%%%%%%%%%%
\section{Introduction and main results}
In this paper, we consider a multi-dimensional transport equation with nonlocal velocity and fractional dissipation
\begin{equation}\label{1.1}
\left\{\ba
&\partial_t\theta+u\cdot \nabla\theta +\kappa\Lambda^\gamma\theta=0,~(x,t)\in \mathbb{R}^n\times\mathbb{R}_+,\\
&u=\nabla\Lambda^{-2+2\alpha}\theta,\\
&\theta(x,0)=\theta_{0}(x),~x\in\mathbb{R}^n, \ea\ \right.
\end{equation}
 where $n\geq2$, $\alpha\in(0,\frac{1}{2}]$, $\gamma\in(0,2)$ and $\kappa>0$.  Here $\theta:\mathbb{R}^n\times\R_+\mapsto\R$ is the unknown function, and the fractional Laplacian $\Lambda^s=(-\Delta)^{\frac{s}{2}}$ with $s\in \R$ is defined by the Fourier transform:
 \begin{eqnarray*}
\widehat{\Lambda^sf}(\xi)=(2\pi|\xi|)^s\widehat{f}(\xi).
\end{eqnarray*}
The non-locality of velocity $u$ is due to the singular integral representation of $\eqref{1.1}_2$ (see, e.g., \cite{[Stein]}):
\begin{eqnarray}\label{1.2}
u(x,t)=C_{n,\alpha}P.V.\int_{\R^{n}}\frac{x-y}{|x-y|^{n+2\alpha}}\theta(y,t)dy,
\end{eqnarray}
where
\begin{eqnarray*}
C_{n,\alpha}
=-\frac{2^{2\alpha-1}\Gamma(\frac{n}{2}+\alpha)}{\pi^{\frac{n}{2}}\Gamma(1-\alpha)}.
\end{eqnarray*}
%with  $\Gamma$  the  Gamma function, and P.V. denotes the principle value integration.
\par
The one-dimensional analogous model of \eqref{1.1} has been extensively studied.
When $n=1$ and $\alpha=\frac{1}{2}$, \eqref{1.1} reduces to the well-known C\'{o}rdoba-C\'{o}rdoba-Fontelos model
\begin{eqnarray}\label{dissipative-CCF}
\partial_t\theta-H(\theta)\theta_x+\kappa\Lambda^\gamma\theta=0,
\end{eqnarray}
where $H\theta$ is the Hilbert transform of $\theta$. The inviscid case, i.e., \eqref{dissipative-CCF} with $\kappa=0$, was first proposed by C\'{o}rdoba, C\'{o}rdoba and Fontelos in \cite{[Cordoba-Cordoba-Fontelos05]} as a one-dimensional analogue of the two-dimensional surface quasi-geostrophic equation (SQG), which is an incompressible model (see e.g., \cite{[Constantin-Majda-Tabak],[Maida-Bertozzi]}).
We will not review here the known results for SQG and related model. One can refer to \cite{[Caffarelli-Vasseur],[Chae-Constantin-Cordoba-Gancedo],[Constantin-Vicol],[Constantin-Wu],[Cordoba],[Cordoba-Martinez],[CotiZelati-Vicol],[Jeong-Kim],[Jeong-Kim-Yao],[Kiselev-Nazarov-Volberg],[Kiselev-Ryzhik-Yao-Zlato],[Kiselev-Yao-Zlatos],[Zlatos]} and the references therein  for more details and the recent progress.
%There have been a number of mathematical studies on the well-posedness for SQG. We will not review here in detail the known results for the SQG and related equation. One can refer to \cite{[Caffarelli-Vasseur],[Chae-Constantin-Cordoba-Gancedo],[Constantin-Vicol],[Constantin-Wu],[Cordoba],[Cordoba-Martinez],[CotiZelati-Vicol],[Jeong-Kim],[Jeong-Kim-Yao],[Kiselev-Nazarov-Volberg],[Kiselev-Ryzhik-Yao-Zlato],[Kiselev-Yao-Zlatos],[Zlatos]} and the references therein  for more details and the recent progresses.
Here we briefly summarize the progress related to the CCF model.
Regarding the inviscid CCF model, C\'{o}rdoba et al.\cite{[Cordoba-Cordoba-Fontelos05]} first obtained an ingenious nonlinear weighed inequality
 for the Hilbert transform with the help of Meillin transform and complex analysis, and proved
 smooth solutions must blow up in finite time for a generic family of even initial data (see also \cite{[Cordoba-Cordoba-Fontelos06]} for another proof of the blowup for the initial data not necessarily even).
%By virtue of the Meillin transform and the complex analysis, they obtained some new bilinear
%estimates for the Hilbert transform and as a result proved the breakdown of local smooth solutions to \eqref{1.6}
%with $\kappa=0$ for a generic class of initial data.
%C\'{o}rdoba et al.\cite{[Cordoba-Cordoba-Fontelos05]} first proved the the local smooth solutions may blow up in finite time for a generic family of initial data (see also \cite{[Cordoba-Cordoba-Fontelos06],[Kiselev]}).
This blow-up phenomena was later proved by Silvestre and Vicol in \cite{[Silvestre-Vicol]} via four essentially different methods.
Based on completely real-variable arguments, Li and Rodrigo \cite{[Li-Rodrigo20]} recently provided a short proof of the nonlinear inequality first proved by C\'{o}rdoba et al. in \cite{[Cordoba-Cordoba-Fontelos05]}, and obtained several new weighted inequalities for the Hilbert transform and various nonlinear versions, which can be applied to show the finite-time blow-up of smooth solutions to \eqref{1.1} with $n=1$, $\alpha\in(0,1)$ and $\kappa=0$.
%and to \eqref{dissipative-CCF} for $0<\gamma<\frac12$.
%By virtue of the Meillin transform and the complex analysis, the authors obtained some new bilinear
%estimates for the Hilbert transform and as a result proved the breakdown of smooth solutions to \eqref{1.6}
%for a generic class of initial data.
%Similar result holds for the equation \eqref{M-CCF-T} with $n=1$ and $\alpha\in(0,\frac12)\cup(\frac12,1)$.
For the dissipative case, the authors in \cite{[Cordoba-Cordoba-Fontelos05]} also obtained the global well-posedness of \eqref{dissipative-CCF} for the positive $H^2$ initial data in the subcritical case $\gamma\in(1,2)$, as well as in the case $\gamma=1$ with small initial data. Later in \cite{[Dong2008]}, by adapting the method of continuity first used in \cite{[Kiselev-Nazarov-Volberg]}, Dong established the global regularity of solutions to \eqref{dissipative-CCF} in the critical case $\gamma=1$ for arbitrary initial data in appropriate critical Sobolev space. In \cite{[Li-Rodrigo08]}, Li and Rodrigo first proved the finite-time blowup of smooth solutions to \eqref{dissipative-CCF} in the supercritical case $\gamma\in(0,\frac12)$ (see also \cite{[Kiselev],[Li-Rodrigo20],[Silvestre-Vicol]} for different proofs for this blow-up result).
For the remaining range $\gamma\in[\frac12,1)$, whether smooth solutions to \eqref{dissipative-CCF} may blow up in finite time is currently still open. In \cite{[Silvestre-Vicol]}, the authors conjectured that solutions obtained as vanishing viscosity approximations could be bounded in $C^{\frac12}$ for all $t>0$, which would possibly yield H\"{o}lder regularization effect and then solve the problem of global regularity of solutions to \eqref{dissipative-CCF} with $\gamma\in[\frac12,1)$.
Recently, for each smooth nonnegative $\theta_0$, Ferreira and Moitinho \cite{[Ferreira-Moitinho20]} obtained the existence of global classical solutions to \eqref{dissipative-CCF} for $\gamma\in(\gamma_1,1)$ with $\gamma_1$ depending on the $H^{\frac32}$-norm of $\theta_0$.
Ferreira and Moitinho \cite{[Ferreira-Moitinho]} also studied the interpolation $\alpha$-CCF model, i.e., \eqref{1.1} with $n=1$, $\alpha\in(0,\frac12)$ and $\gamma\in(0,2)$, and obtained a set of results including local existence, global smoothness, eventual regularity and blow-up of solutions.
\par
We continue to review the existing results on the multi-dimensional \eqref{1.1}. The inviscid case of \eqref{1.1} was understood well.
The local well-posedness of \eqref{1.1} with $\kappa=0$ in the Sobolev space $H^s$ for appropriate $s>0$ was established by Chae in \cite{[Chae]}.
For $\kappa=0$ and $\alpha=\frac12$ in \eqref{1.1}, Balodis and C\'{o}rdoba \cite{[Balodis-Cordoba]} presented a weighted nonlinear inequality for the Riesz transform by using Meillin transform and spherical harmonic expansion, and  obtained the blow-up of smooth solution for any nonnegative, not-identically zero initial data.
When $n=2$, such finite-time singularity formation result was also proved for a similar model in \cite{[Dong-Li]} independently.
For the general dimension $n\geq2$ and the full range $\alpha\in (0,1)$, Dong later in \cite{[Dong]} was able to obtain the blow-up of smooth solutions to \eqref{1.1} with $\kappa=0$ for any smooth, radially symmetric and nonnegative $\theta_0$ with compact support and its positive maximum attained at the origin. Recently, motivated by \cite{[Silvestre-Vicol]}, Jiu and Zhang \cite{[Jiu-Zhang]} proved the finite-time singularity of solutions to \eqref{1.1} with $\kappa=0$, $\alpha\in(0,1)$ and $n\geq2$ for smooth initial data $\theta_0$ with $\displaystyle\sup_{x\in\R^n}\theta_0(x)>0$ via the De Giorgi iteration technique.
This iteration strategy in \cite{[Silvestre-Vicol]} was later adapted by Alonso-Or\'{a}n and Mart\'{i}nez \cite{[A-O-M]} to the proof of finite-time blow-up for non-local active scalar equations \eqref{1.1} with $\kappa=0$ on compact Riemannian manifolds.
For the two-dimensional and fractionally dissipative case , Li and Rodrigo \cite{[Li-Rodrigo09]} proved the finite-time blow-up of radial smooth solutions to \eqref{1.1} for $\alpha\in[\frac{1}{4},1)$ and $\gamma\in(0,\alpha)$.
Very recently, Li, Liu and Zhang \cite{[Li-Li-Zhang]} studied a related nonlocal transport model with the power type of damping, which reads as
\begin{equation}\label{M-CCF-T-damping}
\left\{\ba
&\partial_t\theta+u\cdot\nabla \theta +\kappa|\theta|^{\nu-1}\theta= 0, ~(x,t)\in \R^{n}\times\R_+,\\&u=\nabla\Lambda^{-2+2\alpha}\theta,\\
&\theta(x,0)=\theta_{0}(x),~x\in\R^n, \ea\ \right.
\end{equation}
where $n\geq1$, $\alpha\in(0,1)$, $\kappa\in \R$ and $\nu>0$. In \cite{[Li-Li-Zhang]}, by some change of time variable to implement the iteration technique originated from \cite{[Silvestre-Vicol]}, the authors showed that the damping term can not avoid the singularity formation in finite time, and that for particular $\theta_0$ depending on the size of its $L^1$ and $L^\infty$-norm, the solutions to \eqref{M-CCF-T-damping} must blow up in finite time, independent of the value of $\kappa$.
%Finally, It should be remarked that the two-dimensional dissipative SQG was extensively studied.
%We will not review here in detail the known results for the SQG and related equation. One can refer to \cite{[Caffarelli-Vasseur],[Constantin-Vicol],[Constantin-Wu],[CotiZelati-Vicol],[Kiselev-Nazarov-Volberg],[Zlatos]} and the references therein  for more details and the recent progresses.

In the present paper, we investigate the multi-dimensional nonlocal transport \eqref{1.1} with fractional
dissipation. Due to scaling property, if $\theta(x,t)$ solves \eqref{1.1} with the initial datum $\theta_0$, then
for any $\lambda>0$,
\begin{eqnarray*}
\theta_\lambda(x,t)\triangleq\lambda^{\gamma-2\alpha}\theta(\lambda x,\lambda^\gamma t)
\end{eqnarray*}
so does with the initial datum $\lambda^{\gamma-2\alpha}\theta_0(\lambda x)$. Therefore, the initial value problem \eqref{1.1} admits three basic cases: subcritical $\gamma\in(2\alpha,2)$, critical $\gamma=2\alpha$ and supercritical $\gamma\in(0,2\alpha)$.
For the critical and subcritical dissipation, our main result on \eqref{1.1} can be stated as
\begin{theorem}\label{the-1}
Let $n\geq2$, $\alpha\in(0,\frac12]$ and $\gamma\in[2\alpha,2)$.
For every nonnegative initial datum $\theta_0\in \mathcal{S}(\mathbb{R}^n)$, there exists a unique global-in-time smooth solutions to \eqref{1.1}.
\end{theorem}
In the case of supercritical dissipation with $\gamma$ close to $2\alpha$, we have the next type of global regularity result.
\begin{theorem}\label{the-2}
Let $n\geq2$, $\alpha\in(0,\frac12]$ and $\gamma\in(\alpha,2\alpha)$.
Let $R>0$ be arbitrary. For every nonnegative initial datum $\theta_0\in H^s(\mathbb{T}^n)$ with $s>\frac{n}{2}+1$ such that
\begin{eqnarray}\label{assumption-initial-data}
\|\theta_0\|^{\frac{n+4\alpha-2\gamma}{2s}}_{\dot{H}^s(\mathbb{T}^n)}
\|\theta_0\|^{1-\frac{n+4\alpha-2\gamma}{2s}}_{L^2(\mathbb{T}^n)}\leq R,
\end{eqnarray}
there exists $\gamma_1=\gamma_1(R)\in(\alpha,2\alpha)$ such that for each $\gamma\in[\gamma_1,2\alpha)$ the initial problem \eqref{1.1} admits a unique global $H^s$-solution.
\end{theorem}
For the case of supercritical dissipation with appropriately small $\gamma$, we obtain the finite time blowup result.
\begin{theorem}\label{the-3}
Let $n\geq2$, $\alpha\in(0,1)$ and $\gamma\in(0,\alpha)$.
For any radial initial datum $\theta_0\in\mathcal{S}(\mathbb{R}^n)$ such that for some constant $A>0$ depending only on $n$, $\alpha$ and $\gamma$,
\begin{eqnarray*}
\int_{\R^n}\frac{\theta_0(0)-\theta_0(x)}{|x|^n}e^{-|x|}dx\geq A(1+\|\theta_0\|_{L^\infty}),
\end{eqnarray*}
the corresponding smooth solution $\theta$ to \eqref{1.1} must blow up in finite time.
\end{theorem}
%\begin{remark}
%It should be remarked that Dong \cite{[Dong]} claimed without proofs that one can utilize the robust technique of the modulus of continuity originated from Kiselev, Nazarov and Volberg in \cite{[Kiselev-Nazarov-Volberg]} to show that the solution to \eqref{1.1} with $\alpha\in(0,\frac12]$ and $\gamma\in[2\alpha,2)$ is globally well-posed in time.
%Indeed, one can refer to \cite{[Kiselev-Nazarov-Volberg]} for the construction of the modulus of continuity for the case $\alpha=\frac12,\gamma=1$, and refer to \cite{[Ferreira-Moitinho]} for $\alpha\in(0,\frac12)$ and $\gamma\in[2\alpha,2)$.
%The general multi-dimensional case is completely similar.
%For the transparent way in which the dissipation is dominating the nonlinear term,
%we will instead utilize the nonlinear maximum principle presented by Constantin and Vicol in \cite{[Constantin-Vicol]} to show the global regularity of solutions to \eqref{1.1} in the subcritical and critical case.
%\end{remark}
\begin{remark}
We make some comments on Theorems \ref{the-1}-\ref{the-3}. Firstly, Dong claimed without proofs in \cite{[Dong]} that one can modify the arguments of nonlocal moduli of continuity firstly appeared in \cite{[Kiselev-Nazarov-Volberg]} to prove Theorem \ref{the-1}.
One can refer to \cite{[Kiselev-Nazarov-Volberg]} for the construction of the modulus of continuity for the case $\alpha=\frac12,\gamma=1$, and refer to \cite{[Ferreira-Moitinho]} for $\alpha\in(0,\frac12)$ and $\gamma\in[2\alpha,2)$. For the transparent way in which the dissipation dominates the nonlinear term,
we show Theorem \ref{the-1} by utilizing the nonlinear maximum principle firstly presented in \cite{[Constantin-Vicol]}. Secondly, the author also claimed without proofs in \cite{[Dong]}  that, for $\gamma\in(0,\alpha)$ and $\delta\in(0,2\alpha-2\gamma)$, there exists a constant $N$ depending only on $n,\alpha,\gamma$ and $\delta$, such that if
\begin{eqnarray*}
\int_{\mathbb{R}^n}\frac{\theta_0(0)-\theta_0(x)}{|x|^{n+\delta}}dx\geq N(1+\|\theta_0\|_{L^\infty}),
\end{eqnarray*}
then the solution to \eqref{1.1} with the radially symmetric $\theta_0$ blows up in finite time.
We revisit this problem of the finite time blowup of smooth solutions to \eqref{1.1} with $\gamma\in(0,\alpha)$ and obtain Theorem \ref{the-3} by some different arguments.
%We in this paper give an alternative proof of this blow-up result, i.e., Theorem \ref{the-3}.
Finally, the 1-d analogous conclusions of Theorem \ref{the-2} for $\alpha=\frac12$ and $\alpha\in(0,\frac12)$ were established in \cite{[Ferreira-Moitinho20]} and \cite{[Ferreira-Moitinho]}, respectively.
\end{remark}
%\begin{remark}
%Dong \cite{[Dong]} also claimed without proofs that by using the nonlinear inequality
%\begin{eqnarray*}
%\int_{\R^n}\frac{\Lambda^{-2+2\alpha}\nabla f(x)\cdot\nabla f(x)}{|x|^{n+\delta}}dx
%\geq C_{n,\alpha,\delta}\int_{\R^n}\frac{(f(0)-f(x))^2}{|x|^{n+2\alpha+\delta}}dx,
%\end{eqnarray*}
%and following the arguments in \cite{[Li-Rodrigo09]}, there exists a constant $N$ depending only on $n,\alpha,\gamma$ and $\delta$, such that if
%\begin{eqnarray*}
%\int_{\mathbb{R}^n}\frac{\theta_0(0)-\theta_0(x)}{|x|^{n+\delta}}dx\geq N(1+\|\theta_0\|_{L^\infty}),
%\end{eqnarray*}
%then the solution to \eqref{1.1} with the radially symmetric initial data $\theta_0$ blows up in finite time.
%Here we give an alternative proof of this blow-up result. More precisely,
%\end{remark}
Now we explain our approaches in some details. Considering all basic cases of $\gamma$, that is, for $\alpha\in(0,\frac12]$ and $\gamma\in(0,2)$, we firstly study the local well-posedness of solutions to \eqref{1.1} and after provide a Beale-Kato-Majda blow-up criterion and control on the existence time for initial datum $\theta_0$ in $H^s$ with $s>\frac{n}{2}+1$ (see Propositions \ref{local}-\ref{BKM} and Corollary \ref{lower-bound-time}). In the critical and subcritical cases we show the global existence of smooth solutions via an approach based on the nonlinear maximum principle inspired by \cite{[Constantin-Vicol]}. For the supercritical case, with the help of a conditional regularity result (see Lemma \ref{conditional-regularity}), we prove an eventual regularization result, i.e., Proposition \ref{eventual-regularity}, and then conclude the global existence of solutions when $\gamma$ is sufficiently close to $2\alpha$ for nonnegative initial datum.
%Finally, motivated by \cite{[Li-Rodrigo20]}, we establish the finite-time blow-up of solutions for a generic class of radial initial data in the part $\gamma\in(0,\alpha)$ of the supercritical range.
%Finally, for the supercritical case $\gamma\in(0,\alpha)$, we use a nonlinear weighted inequality from \cite{[Zhang]} (see Lemma \ref{nonlinear-inequality-exponential-weight}) and establish a new inequality involving the dissipation (see Lemma \ref{dissipation}) to show the finite time blowup result.
Finally, for the supercritical case $\gamma\in(0,\alpha)$, we establish a new weighted inequality involving the dissipation along with the use of a nonlinear weighted inequality from \cite{[Zhang]} to show the finite time blowup result (see Lemmas \ref{nonlinear-inequality-exponential-weight} and \ref{dissipation}).

The outline of this paper is organized as follows. In Section $2$, we recall some definitions and basic facts needed later.
Section 3 is devoted to the local well-posedness, the Beale-Kato-Majda criterion, and a suitable control on the existence time for the model \eqref{1.1}.
%In the subsequent sections, we will prove the main theorems.
In Section 4, we prove Theorem \ref{the-1}.
%In Section 4, based on the nonlinear maximum principle, we will establish the global regularity of the solutions to \eqref{1.1} with the subcritical and critical dissipation.
%Finally, in Section 6, we prove the finite-time blow-up
The proof of Theorems \ref{the-2} and \ref{the-3} is given in Sections $5$ and $6$, respectively.
\section{Preliminaries}
In this section, we will present some notations and basic facts needed later.
\subsection{Notations}
For $p\in[1,\infty]$, we denote $L^p(\Omega)$ the standard $L^p$-space and its  norm by  $\|\cdot\|_{L^p(\Omega)}$.
For $s\in\R$, we use $\dot{H}^s(\mathbb{T}^n)$ to denote the homogeneous  Sobolev space of $s$ order,
whose endowed norm is denoted by $\|\cdot\|_{\dot{H}^s(\mathbb{T}^n)}=\|\Lambda^s(\cdot)\|_{L^2(\mathbb{T}^n)}$. Here the fractional Laplacian $\Lambda^s=(-\Delta)^{\frac{s}{2}}$ with $s\in \R$ is defined by the Fourier transform (see e.g., \cite{[Bahouri-Chemin-Danchi]}).  For $s\in(0,2)$, another equivalent definition is given as follows: for $f\in C^\infty(\mathbb{T}^n)$,
\begin{eqnarray}\label{Fractional-Laplacian}
(\Lambda^s f)(x)
=C_{n,s}\sum_{k\in\mathbb{Z}^n}\int_{\mathbb{T}^n}
\frac{f(x)-f(x+y)}{|y-2\pi k|^{n+s}}dy=C_{n,s}P.V.\int_{\R^n}\frac{f(x)-f(x+y)}{|y|^{n+s}}dy,
\end{eqnarray}
where $C_{n,s}$ is a normalization constant (see, e.g., \cite{[Di Nezza-Palatucci-Valdinoci]}).
In the above identity, and throughout this paper, we abuse notation and still denote by $f$ the periodic extension of $f$ to the whole space.
Moreover, for $\beta\in(0,1)$, we recall the H\"{o}lder space $C^\beta(\Omega)$ whose norm is given by
$\|f\|_{C^\beta(\Omega)}=\|f\|_{L^\infty(\Omega)}+[f]_{C^\beta(\Omega)}$, where
$[\cdot]_{C^\beta(\Omega)}$ means the semi-norm
\begin{eqnarray*}
[f]_{C^\beta(\Omega)}
=\sup_{x,y\in\Omega,x\neq y}\frac{|f(x)-f(y)|}{|x-y|^\beta}.
\end{eqnarray*}
For convenience of the notation, the $L^{p}(\Omega)$-norm of a function $f$ is always abbreviated as $\norm{f}_{L^p}$ and the $\dot{H}^s(\mathbb{T}^n)$-norm as $\norm{f}_{\dot{H}^s}$. The positive part of $\ln a$ with $a>0$ is given by
\begin{eqnarray*}
\ln^+a
=
\begin{cases}
\ln a,
&  \mbox{{\rm for} $a\geq1,$ }\\
0,
& \mbox{{\rm for} $0<a<1.$ } \\
\end{cases}
\end{eqnarray*}
It is clear that, for any $a>0$,
\begin{eqnarray}\label{inequality-log}
1+\ln^+a\leq2\ln(a+e).
\end{eqnarray}
Throughout this paper, we will denote by $\chi$ a smooth radially nondecreasing cut-off function that vanishes for $|x|\leq\frac12$, is equal to 1 for $|x|\geq1$, and its gradient satisfies $\|\nabla\chi\|_{L^\infty}\leq2$. We will use $C$ to denote a positive constant, whose value may change from line to line, and write $C_{n,\alpha}$ or $C(n,\alpha)$ to emphasize the dependence of a constant on $n$ and $\alpha$.
\subsection{Some basic facts}
In the following, we present some basic facts needed later.
We begin with a classical commutator estimate for fractional Laplacian (see, e.g., \cite{[Ju]}).
%, whose proof can be found in \cite{[Ju]} and their references. The following lemma contains such a type of estimate
\begin{lemma}\label{commutator}
Let $s>0$. Let $p,p_1,p_4\in(1,\infty)$, $p_2,p_3\in(1,\infty]$ such that
\begin{eqnarray*}
\frac{1}{p}=\frac{1}{p_1}+\frac{1}{p_2}=\frac{1}{p_3}+\frac{1}{p_4}.
\end{eqnarray*}
For $f,g\in C^\infty(\mathbb{T}^n)$ or $\mathcal{S}(\mathbb{R}^n)$, there exists a constant $C=C_{n,s,p,p_1,p_3}>0$ such that
\begin{eqnarray*}
\|\Lambda^{s}(fg)-f\Lambda^{s}g\|_{L^p}\leq C (\|\Lambda^{s}f\|_{L^{p_1}}\|g\|_{L^{p_2}}+\|\nabla f\|_{L^{p_3}}\|\|\Lambda^{s-1}g\|_{L^{p_4}}).
\end{eqnarray*}
\end{lemma}
For the fractional operator $\Lambda^\gamma$, we next recall a pointwise identity, which follows from the arguments in
\cite{[Cordoba-Cordoba-CMP04]} (see also \cite{[Constantin]}).
%Next we recall a pointwise identity of the operator. One can refer to \cite{[Cordoba-Cordoba-CMP04]} for more details.
%The reader is referred to \cite{[Cordoba-Cordoba-CMP04]} for further details.
\begin{lemma}\label{pointwise-property}
Let $\gamma\in(0,2)$ and  $f\in C^\infty(\mathbb{T}^n)$. Then, for all $x\in\mathbb{T}^n$,
\begin{eqnarray*}
f(x)\Lambda^\gamma f(x)=\frac12\Lambda^\gamma(f^2)(x)+\frac12D_\gamma(f)(x),
\end{eqnarray*}
where the operator $D_\gamma$ is defined by
\begin{eqnarray}\label{definition-dissipation-remainder}
D_\gamma(f)(x)
=c_{n,\gamma}\sum_{k\in\mathbb{Z}^n}\int_{\mathbb{T}^n}
\frac{|f(x)-f(x+y)|^2}{|y-2\pi k|^{n+\gamma}}dy=c_{n,\gamma} P.V.\int_{\R^n}\frac{|f(x)-f(x+y)|^2}{|y|^{n+\gamma}}dy.
\end{eqnarray}
\end{lemma}
We proceed to present the $L^\infty$ and $C^\beta$ nonlinear lower bounds for the operator $D_\gamma$, whose proof can be found in \cite{[Constantin-Vicol]}.
\begin{lemma}\label{maximum-lower-bound}
Let $\gamma\in(0,2)$ and $f\in\mathcal{S}(\R^n)$. Then, for all $x\in\R^n$,
\begin{eqnarray*}
D_\gamma(\nabla f)(x)\geq\frac{|\nabla f(x)|^{2+\gamma}}{c_{n,\gamma}\|f\|^\gamma_{L^\infty}}.
\end{eqnarray*}
In addition, for any $\beta\in(0,1)$ and for all $x\in\R^n$,
\begin{eqnarray*}
D_\gamma(\nabla f)(x)\geq C\frac{|\nabla f(x)|^{2+\frac{\gamma}{1-\beta}}}{\|f\|^{\frac{\gamma}{1-\beta}}_{C^\beta}}.
\end{eqnarray*}
\end{lemma}
We end this subsection with an auxiliary inequality needed later.
\begin{lemma}\label{Inequality}
For any $\rho>0$, we have
\begin{eqnarray*}
\int^\infty_{\rho}\frac{e^{-r}}{r}dr
\leq2\ln\Big(e+\frac{1}{\rho}\Big).
\end{eqnarray*}
\end{lemma}
\textbf{Proof}.
For $\rho>1$,
\begin{eqnarray*}
\int^\infty_{\rho}\frac{e^{-r}}{r}dr\leq\int^\infty_{1}\frac{e^{-r}}{r}dr
<1.
\end{eqnarray*}
For the remaining value $\rho\in(0,1]$, by integration by parts, we derive that
\begin{eqnarray*}
\int^\infty_{\rho}\frac{e^{-r}}{r}dr
=\frac{\ln(\frac{1}{\rho})}{e^{\rho}}+\int^\infty_{\rho}\frac{\ln r}{e^r}dr
\leq\ln\Big(\frac{1}{\rho}\Big)+\int^\infty_1\frac{ \ln r}{e^r}dr
\leq\ln\Big(\frac{1}{\rho}\Big)+1.
\end{eqnarray*}
Therefore, it follows from \eqref{inequality-log} that, for any $\rho>0$,
\begin{eqnarray*}
\int^\infty_{\rho}\frac{e^{-r}}{r}dr
\leq\ln^+\Big(\frac{1}{\rho}\Big)+1
\leq2\ln\Big(e+\frac{1}{\rho}\Big),
\end{eqnarray*}
which concludes the proof of Lemma \ref{Inequality}.
\subsection{Some properties of solutions}
\begin{lemma}\label{maximum}
Let $n\geq2$, $\alpha\in(0,1)$, $\gamma\in(0,2)$ and $\Omega=\mathbb{T}^n$ or $\mathbb{R}^n$.
If $\theta$ is a smooth solution to \eqref{1.1} on $[0,T]$ for nonnegative $\theta_0\in H^s(\Omega)$ with $s>\frac{n}{2}+1$, then
\begin{eqnarray*}
0\leq\theta(x,t)\leq\|\theta_0\|_{L^\infty(\Omega)}, ~for ~all~ (x,t)\in\Omega\times[0,T],
\end{eqnarray*}
and
\begin{eqnarray}\label{L2-decreasing}
\|\theta(t)\|_{L^2(\Omega)}\leq\|\theta_0\|_{L^2(\Omega)}, ~for ~all~t\in [0,T].
\end{eqnarray}
\end{lemma}
\textbf{Proof}. Define
\begin{eqnarray*}
m(t)=\min_{x\in\Omega}\theta(x,t),~M(t)=\max_{x\in\Omega}\theta(x,t).
\end{eqnarray*}
Following completely analogous arguments as Lemma 3.2 in \cite{[Jiu-Zhang]}, we can show that $m(t)$ is non-decreasing, and that $M(t)$ is non-increasing on $[0,T]$. Thus, for all $(x,t)\in\mathbb{T}^n\times[0,T]$,
\begin{eqnarray}\label{non-negative}
0\leq\min_{x\in\Omega}\theta_0(x)\leq m(t)\leq\theta(x,t)\leq M(t)\leq \max_{x\in\Omega}\theta_0(x).
\end{eqnarray}
Multiplying \eqref{1.1} by $\theta$ and integrating by parts lead to
\begin{eqnarray}\label{L2}
\frac{1}{2}\frac{d}{dt}\|\theta(t)\|^2_{L^2}+\kappa\|\Lambda^{\frac{\gamma}{2}}\theta(t)\|^2_{L^2}
=\frac12\int_{\Omega}\theta^2 {\rm div} u dx.
\end{eqnarray}
By the parity, \eqref{Fractional-Laplacian} and \eqref{non-negative}, we note that
\begin{eqnarray*}
\begin{split}
\int_{\Omega}\theta^2 {\rm div} u dx
&=-\int_{\Omega}\theta^2\Lambda^{2\alpha}\theta dx
=-C_{n,\alpha}\iint\limits_{\R^n\times\R^n}\theta^2(x,t)\frac{\theta(x,t)-\theta(y,t)}{|x-y|^{n+2\alpha}}dxdy\\
&=-\frac{C_{n,\alpha}}{2}\iint\limits_{\R^n\times\R^n}\frac{(\theta(x,t)+\theta(y,t))(\theta(x,t)-\theta(y,t))^2}{|x-y|^{n+2\alpha}}dxdy
\leq0,
\end{split}
\end{eqnarray*}
which along with \eqref{L2} yields that
\begin{eqnarray*}
\frac{d}{dt}\|\theta(t)\|_{L^2(\Omega)}\leq0.
\end{eqnarray*}
By direct integration, it immediately gives \eqref{L2-decreasing}. We complete the proof of Lemma \ref{maximum}.

At the end of this section, we show that the radial symmetry of initial data is preserved by solutions to \eqref{1.1}.
\begin{lemma}\label{radial-property-preserved}
Let $n\geq2$, $\alpha\in(0,1)$ and $\gamma\in(0,2)$.
If $\theta$ is a smooth local solution to \eqref{1.1} with radially symmetric $\theta_0\in\mathcal{S}(\mathbb{R}^n)$, then $\theta(x,t)$ is radial for all time $t>0$.
\end{lemma}
\textbf{Proof}. Let $\mathbf{O}$ be any orthogonal matrix. By the uniqueness of solutions to \eqref{1.1} (c.f. Proposition \ref{local} below) and  the radial property of $\theta_0$, it suffices to show that $\theta_{\mathbf{O}}(x,t):=\theta(\mathbf{O}x,t)$ is also a solution to \eqref{1.1} with the initial data $\theta_0(\mathbf{O}x)$. Indeed, following some computations in \cite{[Zhang]}, we have
\begin{eqnarray*}
(\partial_t\theta_\mathbf{O})(x,t)
%=\partial_t(\theta(\mathbf{O}x,t))
=(\partial_t\theta)(\mathbf{O}x,t),~(\nabla_x\theta_\mathbf{O})(x,t)=\mathbf{O}^T(\nabla_x\theta)(\mathbf{O}x,t),
\end{eqnarray*}
%\begin{eqnarray*}
%(\nabla_x\theta_\mathbf{O})(x,t)
%=\nabla_x(\theta(\mathbf{O}x,t))
%=\mathbf{O}^T(\nabla_x\theta)(\mathbf{O}x,t),
%\end{eqnarray*}
and
\begin{eqnarray*}
u_\mathbf{O}(x,t)=\nabla\Lambda^{-2+2\alpha}\theta_\mathbf{O}=\mathbf{O}^{-1} u(\mathbf{O} x,t).
\end{eqnarray*}
%\begin{eqnarray*}
%u_\mathbf{O}(x,t)
%&=&C_{n,\alpha}P.V.\int_{\R^{n}}\frac{x-y}{|x-y|^{n+2\alpha}}\theta_{\mathbf{O}}(y,t)dy \nonumber \\
%&=&C_{n,\alpha}P.V.\int_{\R^{n}}\frac{x-y}{|x-y|^{n+2\alpha}}\theta(\mathbf{O} y, t)dy\nonumber\\
%&=& C_{n,\alpha}P.V.\int_{\R^{n}}\frac{x-\mathbf{O}^{-1}z}{|x-\mathbf{O}^{-1}z|^{n+2\alpha}}\theta(z, t)dz \nonumber \\
%&=& \mathbf{O}^{-1}C_{n,\alpha}P.V.\int_{\R^{n}}\frac{\mathbf{O} x-z}{|\mathbf{O} x-z|^{n+2\alpha}}\theta(z,t)dz \nonumber \\
%&=&\mathbf{O}^{-1} u(\mathbf{O} x,t),
%\end{eqnarray*}
By \eqref{Fractional-Laplacian} and a change of variables, we also have
\begin{eqnarray*}
&&(\Lambda^\gamma\theta_\mathbf{O})(x,t)=C_{n,\gamma}P.V.\int_{\R^n}\frac{\theta(\mathbf{O}x,t)-\theta(\mathbf{O}y,t)}{|x-y|^{n+\gamma}}dy\\
%&=&C_{n,\gamma}P.V.\int_{\R^n}\frac{\theta(\mathbf{O}x,t)-\theta(z,t)}{|x-O^{-1}z|^{n+\gamma}}dz\\
&&\quad
=C_{n,\gamma}P.V.\int_{\R^n}\frac{\theta(\mathbf{O}x,t)-\theta(z,t)}{|Ox-z|^{n+\gamma}}dz=\Lambda^\gamma\theta(\mathbf{O}x,t).
\end{eqnarray*}
Finally, we derive that
\begin{eqnarray*}
&&(\partial_t\theta_{\mathbf{O}}+u_{\mathbf{O}}\cdot\nabla\theta_{\mathbf{O}}+\Lambda^\gamma\theta_\mathbf{O})(x,t)\\
&=&(\partial_t\theta)(\mathbf{O}x,t)+\mathbf{O}^{-1} u(\lambda x,t)\cdot \mathbf{O}^T(\nabla_x\theta)(\mathbf{O}x,t)+\Lambda^\gamma\theta(\mathbf{O}x,t)\\
%&=&(\partial_t\theta)(\mathbf{O}x,t)+u(\lambda x,t)\cdot \mathbf{O}\mathbf{O}^T(\nabla_x\theta)(\mathbf{O}x,t)+\Lambda^\gamma\theta(\mathbf{O}x,t)\\
&=&(\partial_t\theta+u\cdot\nabla_x\theta+\Lambda^\gamma\theta)(\mathbf{O}x,t)=0.
\end{eqnarray*}
The proof of Lemma \ref{radial-property-preserved} is then finished.
\section{Local well-posedness and Beale-Kato-Majda criterion}
In this section, we study the local existence of solutions to \eqref{1.1}. Also, we obtain a regularity criterion and a control on the existence time.
\begin{proposition}\label{local}
Let $n\geq2$, $\alpha\in(0,\frac12]$, $\gamma\in(0,2)$ and $\Omega=\mathbb{T}^n$ or $\mathbb{R}^n$.
For each $\theta_0\in H^s(\Omega)$ with $s>\frac n2+1$, there exists $T=T(\|\theta_0\|_{H^s})>0$ and a unique solution to \eqref{1.1} such that
\begin{eqnarray}\label{solution-class}
\theta\in L^\infty([0,T);H^s(\Omega)\cap L^2([0,T);H^{s+\frac{\gamma}{2}}(\Omega)).
\end{eqnarray}
\end{proposition}
\textbf{Proof}. We derive the local energy estimate.
By \eqref{L2} and H\"{o}lder's inequality, we have
\begin{eqnarray}\label{3.2}
\frac{1}{2}\frac{d}{dt}\|\theta(t)\|^2_{L^2}+\kappa\|\Lambda^{\frac{\gamma}{2}}\theta(t)\|^2_{L^2}
\leq\frac12\|\theta(t)\|^2_{L^2}\|{\rm div} u\|_{L^\infty}.
\end{eqnarray}
Applying $\Lambda^s$ on $\eqref{1.1}_1$, multiplying by $\Lambda^s\theta$, integrating by parts and utilizing H\"{o}lder's inequality, we can arrive at
\begin{eqnarray}\label{3.3}
&&\frac{1}{2}\frac{d}{dt}\|\Lambda^s\theta(t)\|^2_{L^2}+\kappa\|\Lambda^{s+\frac{\gamma}{2}}\theta(t)\|^2_{L^2}\nonumber\\
&=&-\int_{\Omega}\Lambda^s\theta\Big[\Lambda^s(u\cdot\nabla\theta)-u\cdot\nabla\Lambda^s\theta\Big] dx+\frac12\int_{\Omega}{\rm div}u(\Lambda^s\theta)^2 dx\nonumber\\
&\leq&\|\Lambda^s\theta\|_{L^2}\|\Lambda^s(u\cdot\nabla\theta)-u\cdot\nabla\Lambda^s\theta\|_{L^2}
+\frac12\|{\rm div}u\|_{L^\infty}\|\Lambda^s\theta\|^2_{L^2}.
\end{eqnarray}
%By H\"{o}lder inequality, it follows that
%\begin{eqnarray}\label{3.3}
%\frac{1}{2}\frac{d}{dt}\|\Lambda^s\theta(t)\|^2_{L^2}+\kappa\|\Lambda^{s+\frac{\gamma}{2}}\theta(t)\|^2_{L^2}
%&\leq&\|\Lambda^s\theta\|_{L^2}\|\Lambda^s(u\cdot\nabla\theta)-u\cdot\nabla\Lambda^s\theta\|_{L^2}\nonumber\\
%&&\ \ \ \+\frac12\|{\rm div}u\|_{L^\infty}\|\Lambda^s\theta\|^2_{L^2}.
%\end{eqnarray}
By virtue of Lemma \ref{commutator}, the commutator in \eqref{3.3} can be bounded as
\begin{eqnarray}\label{3.4}
\|\Lambda^s(u\cdot\nabla\theta)-u\cdot\nabla\Lambda^s\theta\|_{L^2}
\leq C(\|\Lambda^su\|_{L^2}\|\nabla\theta\|_{L^\infty}
+\|\Lambda^{s-1}\nabla\theta\|_{L^2}\|\nabla u\|_{L^\infty}).
\end{eqnarray}
For $\alpha\in(0,\frac{1}{2}]$ and $s>\frac{n}{2}+1$, $L^2$-boundedness of the Riesz transform and the Gagliardo-Nirenberg inequality lead us to following estimates
\begin{eqnarray}\label{3.5}
\|\Lambda^su\|_{L^2}
\leq C\|\Lambda^{s-1+2\alpha}\theta\|_{L^2}
\leq C\|\theta\|^{\frac{1-2\alpha}{s}}_{L^2}\|\Lambda^{s}\theta\|^{1-\frac{1-2\alpha}{s}}_{L^2},
\end{eqnarray}
\begin{eqnarray}\label{3.6}
\|\nabla\theta\|_{L^\infty}
\leq C\|\theta\|^{1-\frac{n+2}{2s}}_{L^2}\|\Lambda^{s}\theta\|^{\frac{n+2}{2s}}_{L^2},
\end{eqnarray}
and
\begin{eqnarray}\label{3.7}
\begin{split}
\|\nabla u\|_{L^\infty}
&=\|\Lambda^{2\alpha}(\mathcal{R}\otimes\mathcal{R})\theta\|_{L^\infty}\\
&\leq C\|(\mathcal{R}\otimes\mathcal{R})\theta\|^{1-\frac{n+4\alpha}{2s}}_{L^2}\|\Lambda^s(\mathcal{R}\otimes\mathcal{R})\theta\|^{\frac{n+4\alpha}{2s}}_{L^2}\\
&\leq C\|\theta\|^{1-\frac{n+4\alpha}{2s}}_{L^2}\|\Lambda^s\theta\|^{\frac{n+4\alpha}{2s}}_{L^2}.
\end{split}
\end{eqnarray}
Substituting estimates \eqref{3.4}-\eqref{3.7} into \eqref{3.3}, we obtain that
\begin{eqnarray}\label{3.8}
\frac{1}{2}\frac{d}{dt}\|\Lambda^s\theta(t)\|^2_{L^2}+\kappa\|\Lambda^{s+\frac{\gamma}{2}}\theta(t)\|^2_{L^2}
\leq C\|\theta(t)\|^{1-\frac{n+4\alpha}{2s}}_{L^2}\|\Lambda^s\theta(t)\|^{2+\frac{n+4\alpha}{2s}}_{L^2}.
\end{eqnarray}
Finally, inserting \eqref{3.7} into \eqref{3.2}, and adding with \eqref{3.8}, we conclude that
\begin{eqnarray*}
\frac{1}{2}\frac{d}{dt}\|\theta(t)\|^2_{H^s}+\kappa\|\Lambda^{s+\frac{\gamma}{2}}\theta(t)\|^2_{L^2}
\leq C\|\theta(t)\|^{3}_{H^s}.
\end{eqnarray*}
With the help of this energy estimate, we then use the standard technique (see, e.g., \cite{[Temam]}) to obtain the local existence of solutions to \eqref{1.1}.

We proceed to prove the uniqueness of solutions. Suppose that $\theta_1$ and $\theta_2$ are solutions to \eqref{1.1} with the initial data $\theta_0\in H^s(\Omega)$ in the class \eqref{solution-class}.
%The subtraction $\vartheta=\theta_1-\theta_2$ then satisfies
%\begin{equation*}
%\left\{\ba
%&\partial_t\vartheta+(u_1-u_2)\cdot\nabla\theta_1+(u_1-u_2)\cdot\nabla\theta_2+\kappa\Lambda^\gamma\vartheta=0, ~x\in \R^{2}_+,\\
%&\Psi^n(x,t)=0,~x\in \partial\R^{2}_+. \ea\ \right.
%\end{equation*}
Thus $\theta_1-\theta_2$ satisfies
\begin{eqnarray}\label{3.9}
\partial_t(\theta_1-\theta_2)+(u_1-u_2)\cdot\nabla\theta_1+u_2\cdot\nabla(\theta_1-\theta_2)+\kappa\Lambda^\gamma(\theta_1-\theta_2)=0.
\end{eqnarray}
Multiplying \eqref{3.9} by $\theta_1-\theta_2$,  integrating over $\mathbb{T}^n$ and integrating by parts, we can arrive at
\begin{eqnarray*}
\begin{split}
&\frac{1}{2}\frac{d}{dt}\|\theta_1(t)-\theta_2(t)\|^2_{L^2}+\kappa\|\theta_1(t)-\theta_2(t)\|^2_{\dot{H}^{\frac\gamma2}}\\
%=&-\int_{\mathbb{R}^n}((u_1-u_2)\cdot\nabla\theta_1)(\theta_1-\theta_2) dx-\int_{\mathbb{R}^n}(u_2\cdot\nabla(\theta_1-\theta_2))(\theta_1-\theta_2)dx\\
=&\underbrace{-\int_{\Omega}((u_1-u_2)\cdot\nabla\theta_1)(\theta_1-\theta_2)
dx}_{I_1}+\frac{1}{2}\underbrace{\int_{\Omega}(\theta_1-\theta_2)^2\Lambda^{2\alpha}\theta_2dx}_{I_2}.
\end{split}
\end{eqnarray*}
By H\"{o}lder's inequality and the $L^2$-boundedness of Riesz's transform, we can derive that
\begin{eqnarray*}
\begin{split}
|I_1|
&\leq\|\Lambda^{-(1-2\alpha)}\mathcal{R}(\theta_1-\theta_2)\|_{L^{\frac{2n}{n-2+4\alpha}}}\|\nabla\theta_1\|_{L^{\frac{n}{1-2\alpha}}}\|\theta_1-\theta_2\|_{L^2}\\
&\leq C\|\theta_1\|_{H^s}\|\theta_1-\theta_2\|^2_{L^2},
\end{split}
\end{eqnarray*}
where in the last inequality we used the Hardy-Littlewood-Sobolev inequality
\begin{eqnarray*}
\|\Lambda^{-(1-2\alpha)}f\|_{L^{\frac{2n}{n-2+4\alpha}}}\leq C\|f\|_{L^2}~{\rm~for}~\alpha\in\Big(0,\frac12\Big)
\end{eqnarray*}
and the embedding inequality $\|\nabla\theta_1\|_{L^{\frac{n}{1-2\alpha}}}\leq C\|\theta_1\|_{H^s}$ for $\alpha\in(0,\frac12]$ and $s>\frac{n}{2}+1$.

Similarly, we also have
\begin{eqnarray*}
|I_2|
\leq\|\theta_1-\theta_2\|^2_{L^2}\|\Lambda^{2\alpha}\theta_2\|_{L^\infty}
\leq C\|\theta_2\|_{H^s}\|\theta_1-\theta_2\|^2_{L^2}
\end{eqnarray*}
Therefore,
\begin{eqnarray*}
\frac{d}{dt}\|\theta_1(t)-\theta_2(t)\|_{L^2}
\leq C[\|\theta_1\|_{L^\infty([0,T);H^s)}+\|\theta_2\|_{L^\infty([0,T);H^s)}]\|\theta_1(t)-\theta_2(t)\|_{L^2},
\end{eqnarray*}
which along with Gronwall's inequality yields the uniqueness of solutions.
The proof of Proposition \ref{local} is now completed.

Next, we prove the Beale-Kato-Majda blow-up criterion for \eqref{1.1}, which can be stated as
\begin{proposition}\label{BKM}
Let $n\geq2$, $\alpha\in(0,\frac{1}{2}]$ and $\gamma\in(0,2)$. Let $\theta$ be the solution to \eqref{1.1} constructed as in Proposition \ref{local} for nonnegative $\theta_0\in H^s$ with $s>\frac{n}{2}+1$.
Assume $T_\ast$ is the first time such that $\theta$ cannot be continued in the class \eqref{solution-class} to $T=T_\ast$. Then
\begin{eqnarray*}
\displaystyle\limsup_{t\rightarrow T_\ast}\|\theta(t)\|_{H^s}=\infty~{\rm if~and~only~if}~\int_0^{T_\ast}\|\nabla\theta(t)\|_{L^\infty}dt=\infty.
\end{eqnarray*}
\end{proposition}
To prove Proposition \ref{BKM}, we need the $L^\infty$ estimate of the gradient of the velocity, whose proof will be postponed into Appendix A.
\begin{lemma}\label{gradient-velocity}
Let $n\geq2$, $\alpha\in(0,\frac12]$, $\theta\in H^s$ with $s>\frac{n}{2}+1$, and $u$ be defined via $\theta$ by \eqref{1.2}. Then
\begin{eqnarray}\label{gradient-velocity-3.10}
\|\nabla u\|_{L^\infty}
\leq
\begin{cases}
C\Big[1+(1+\ln^+\|\theta\|_{H^s})\|\nabla\theta\|_{L^\infty}+\|\nabla\theta\|_{L^2}\Big],
&  \mbox{{\rm for} $\alpha=\frac12,$ }\\
C(1+\|\nabla\theta\|_{L^\infty}+\|\nabla\theta\|_{L^2}),
& \mbox{{\rm for} $\alpha\in(0,\frac12)$. } \\
\end{cases}
\end{eqnarray}
\end{lemma}
\textbf{Proof of Proposition \ref{BKM}}.
To show the necessity, we will assume $\int_0^{T_\ast}\|\nabla\theta(t)\|_{L^\infty}dt\triangleq M_0<\infty,$
%\begin{eqnarray*}
%\int_0^{T_\ast}\|\nabla\theta(t)\|_{L^\infty}dt\triangleq M_0<\infty,
%\end{eqnarray*}
and then show that there exists some $C_0>0$ such that
\begin{eqnarray}\label{contradiction}
\|\theta(t)\|_{H^k}\leq C_0,~ t<T_\ast,
\end{eqnarray}
which is a contradiction with our assumption.
Indeed, by \eqref{3.2}-\eqref{3.5}, we have
\begin{eqnarray}\label{3.10}
\frac{d}{dt}\|\theta(t)\|_{H^s}\leq C(\|\nabla\theta\|_{L^\infty}+\|\nabla u\|_{L^\infty})\|\theta(t)\|_{H^s}.
\end{eqnarray}
In view of Lemma \ref{gradient-velocity}, we may need the estimate of $\nabla\theta$ in $L^2$.
Applying $\nabla$ on $\eqref{1.1}_1$ and taking its dot product by $\nabla\theta$, we obtain that
\begin{eqnarray}\label{4.1-gradient}
\frac{1}{2}(\partial_t+u\cdot\nabla)|\nabla\theta|^2
+\kappa\nabla\theta\cdot\Lambda^\gamma\nabla \theta+\nabla u:\nabla\theta\cdot\nabla\theta=0,
\end{eqnarray}
%\begin{eqnarray*}
%\partial_t\Big(\frac{1}{2}|\nabla\theta|^2\Big)+u\cdot\nabla\Big(\frac{1}{2}|\nabla\theta|^2\Big)
%+\kappa\sum^n_{i=1}\partial_i\theta\Lambda^\gamma\partial_i\theta+\sum_{1\leq i,j\leq n} \partial_iu^j\partial_i\theta\partial_j\theta=0,
%\end{eqnarray*}
which follows from direct integration, integration by parts and \eqref{L2-decreasing} that,
\begin{eqnarray*}
&&\frac12\frac{d}{dt}\|\nabla\theta\|^2_{L^2}+\kappa\|\Lambda^{\frac{\gamma}{2}}\nabla\theta\|^2_{L^2}
%&=&-\int_{\R^n}u\cdot\nabla\Big(\frac{1}{2}|\nabla\theta|^2\Big)dx-\int_{\R^n}\sum_{1\leq i,j\leq2} \partial_iu^j\partial_i\theta\partial_j\theta dx\\
=\frac{1}{2}\int_{\R^n}{\rm div}u|\nabla\theta|^2dx-\int_{\R^n}\nabla u:\nabla\theta\cdot\nabla\theta dx\\
&&\quad\quad\quad
\leq%\|{\rm div}u\|_{L^2}\|\nabla\theta\|_{L^2}\|\nabla\theta\|_{L^\infty}
C\|\nabla u\|_{L^2}\|\nabla\theta\|_{L^2}\|\nabla\theta\|_{L^\infty}
\leq C\|\Lambda^{2\alpha}\theta\|_{L^2}\|\nabla\theta\|_{L^2}\|\nabla\theta\|_{L^\infty}\\
%&\leq&
%\begin{cases}
%C\|\nabla\theta\|_{L^\infty}\|\nabla\theta\|^2_{L^2},
%&  \mbox{for $\alpha=\frac12,$ }\\
%C\|\theta\|^{1-2\alpha}_{L^2}\|\Lambda\theta\|^{2\alpha}_{L^2}\|\nabla\theta\|_{L^\infty}\|\nabla\theta\|_{L^2},
%& \mbox{for $\alpha\in(0,\frac12),$ } \\
%\end{cases}\\
&&\quad\quad\quad
\leq
\begin{cases}
C\|\nabla\theta\|_{L^\infty}\|\nabla\theta\|^2_{L^2},
&  \mbox{for $\alpha=\frac12,$ }\\
C\|\theta_0\|^{1-2\alpha}_{L^2}\|\nabla\theta\|_{L^\infty}\|\nabla\theta\|^{1+2\alpha}_{L^2},
& \mbox{for $\alpha\in(0,\frac12).$ } \\
\end{cases}
\end{eqnarray*}
%\begin{eqnarray*}
%\frac12\frac{d}{dt}\|\nabla\theta\|^2_{L^2}+\kappa\|\Lambda^{\frac{\gamma}{2}}\nabla\theta\|^2_{L^2}
%&=&-\int_{\R^n}u\cdot\nabla\Big(\frac{1}{2}|\nabla\theta|^2\Big)dx-\int_{\R^n}\sum_{1\leq i,j\leq2} \partial_iu^j\partial_i\theta\partial_j\theta dx\\
%&=&\frac{1}{2}\int_{\R^n}{\rm div}u|\nabla\theta|^2dx-\int_{\R^n}\nabla u:\nabla\theta\cdot\nabla\theta dx\\
%&\leq&%\|{\rm div}u\|_{L^2}\|\nabla\theta\|_{L^2}\|\nabla\theta\|_{L^\infty}
%C\|\nabla u\|_{L^2}\|\nabla\theta\|_{L^2}\|\nabla\theta\|_{L^\infty}\\
%&\leq&C\|\Lambda^{2\alpha}\theta\|_{L^2}\|\nabla\theta\|_{L^2}\|\nabla\theta\|_{L^\infty}\\
%&\leq&
%\begin{cases}
%C\|\nabla\theta\|_{L^\infty}\|\nabla\theta\|^2_{L^2},
%&  \mbox{for $\alpha=\frac12,$ }\\
%C\|\theta\|^{1-2\alpha}_{L^2}\|\Lambda\theta\|^{2\alpha}_{L^2}\|\nabla\theta\|_{L^\infty}\|\nabla\theta\|_{L^2},
%& \mbox{for $\alpha\in(0,\frac12),$ } \\
%\end{cases}\\
%&\leq&
%\begin{cases}
%C\|\nabla\theta\|_{L^\infty}\|\nabla\theta\|^2_{L^2},
%&  \mbox{for $\alpha=\frac12,$ }\\
%C\|\theta_0\|^{1-2\alpha}_{L^2}\|\nabla\theta\|_{L^\infty}\|\nabla\theta\|^{1+2\alpha}_{L^2},
%& \mbox{for $\alpha\in(0,\frac12).$ } \\
%\end{cases}
%\end{eqnarray*}
It follows that
\begin{eqnarray*}
\frac{d}{dt}\|\nabla\theta\|_{L^2}
\leq
\begin{cases}
C\|\nabla\theta\|_{L^\infty}\|\nabla\theta\|_{L^2},
&  \mbox{for $\alpha=\frac12,$ }\\
C\|\theta_0\|^{1-2\alpha}_{L^2}\|\nabla\theta\|_{L^\infty}\|\nabla\theta\|^{2\alpha}_{L^2},
& \mbox{for $\alpha\in(0,\frac12).$ } \\
\end{cases}
\end{eqnarray*}
Solving this differential inequality directly yields that, for any $t\in(0,T_\ast)$,
\begin{eqnarray*}
\begin{split}
\|\nabla\theta(t)\|_{L^2}
&\leq
\begin{cases}
\|\nabla\theta_0\|_{L^2}e^{C\int_0^t\|\nabla\theta(\tau)\|_{L^\infty}d\tau},
&  \mbox{for $\alpha=\frac12,$ }\\
\Big[\|\nabla\theta_0\|^{1-2\alpha}_{L^2}+C_{\alpha,\|\theta_0\|_{L^2}}\int_0^t\|\nabla\theta(\tau)\|_{L^\infty}d\tau\Big]^{\frac{1}{1-2\alpha}},
& \mbox{for $\alpha\in(0,\frac12).$ } \\
\end{cases}\\
&\leq
\begin{cases}
\|\nabla\theta_0\|_{L^2}e^{CM_0},
&  \mbox{for $\alpha=\frac12,$ }\\
\Big(\|\nabla\theta_0\|^{1-2\alpha}_{L^2}+CM_0\Big)^{\frac{1}{1-2\alpha}},
& \mbox{for $\alpha\in(0,\frac12).$ } \\
\end{cases}
\end{split}
\end{eqnarray*}
By virtue of this bound for $\|\nabla\theta\|_{L^2}$, Lemma \ref{gradient-velocity} and \eqref{inequality-log}, we can get
\begin{eqnarray}\label{3.11}
\|\nabla u\|_{L^\infty}
\leq
\begin{cases}
C\Big[1+\|\nabla\theta\|_{L^\infty}\log(e+\|\theta\|_{H^s})\Big],
&  \mbox{for $\alpha=\frac12,$ }\\
C(1+\|\nabla\theta\|_{L^\infty}),
& \mbox{for $\alpha\in(0,\frac12).$ } \\
\end{cases}
\end{eqnarray}
Here $C$ denotes a constant depending on $M_0$ and $T_\ast$.
Substituting \eqref{3.11} into \eqref{3.10}, we can derive that
\begin{eqnarray*}
\frac{d}{dt}\|\theta(t)\|_{H^s}
%&\leq&C\Big(\|\nabla\theta\|_{L^\infty(\R^n)}+\|\nabla u\|_{L^\infty(\R^n)}\Big)\|\theta(t)\|_{H^s}\\
\leq
\begin{cases}
C\Big[1+\|\nabla\theta\|_{L^\infty}+\|\nabla\theta\|_{L^\infty}\log(e+\|\theta\|_{H^s})\Big]\|\theta(t)\|_{H^s},
&  \mbox{for $\alpha=\frac12,$ }\\
C(1+\|\nabla\theta\|_{L^\infty})\|\theta(t)\|_{H^s},
& \mbox{for $\alpha\in(0,\frac12)$, } \\
\end{cases}
\end{eqnarray*}
It follows from Gronwall's inequality that, for $t\in(0,T_\ast)$,
\begin{eqnarray*}
\begin{split}
\|\theta(t)\|_{H^s}
&\leq
\begin{cases}
e^{[\log(e+\|\theta_0\|_{H^s})+C\int_0^t(1+\|\nabla\theta(\tau)\|_{L^\infty})d\tau]e^{C\int_0^t\|\nabla\theta(\tau)\|_{L^\infty}d\tau}},
&  \mbox{for $\alpha=\frac12,$ }\\
\|\theta_0\|_{H^s}e^{C\int_0^t(1+\|\nabla\theta(\tau)\|_{L^\infty})d\tau},
& \mbox{for $\alpha\in(0,\frac12),$ } \\
\end{cases}\\
&\leq
\begin{cases}
e^{[\log(e+\|\theta_0\|_{H^s})+C(T_\ast+M_0)]e^{CM_0}},
&  \mbox{for $\alpha=\frac12,$ }\\
\|\theta_0\|_{H^s}e^{C(T_\ast+M_0)},
& \mbox{for $\alpha\in(0,\frac12),$ } \\
\end{cases}
\end{split}
\end{eqnarray*}
%which exactly shows that $\|\theta(t)\|_{H^s}$ is bounded by a constant depending on $M_0$, $T_\ast$, and $\|\theta_0\|_{H^s}$, and \eqref{contradiction} is established.
which gives \eqref{contradiction}.
The sufficient part follows easily by Sobolev's inequality: for $s>\frac n2+1$,
\begin{eqnarray*}
\int_{0}^{T_\ast}\|\nabla\theta(t)\|_{L^\infty}dt
\leq CT_\ast\sup_{t\in[0,T_\ast]}\|\theta(t)\|_{H^s}.
\end{eqnarray*}
The proof of Proposition \ref{BKM} is now finished.

We conclude this section by providing a lower bound for the existence time.
\begin{coro}\label{lower-bound-time}
Let $n\geq2$, $\alpha\in(0,\frac{1}{2}]$ and $\gamma\in(0,2)$. For every nonnegative $\theta_0\in H^s(\mathbb{T}^n)$ with $s>\frac n2+1$, consider the unique local-in-time solution
\begin{eqnarray*}
\theta\in L^\infty([0,T);H^s(\mathbb{T}^n))\cap L^2([0,T);H^{s+\frac{\gamma}{2}}(\mathbb{T}^n))
\end{eqnarray*}
to \eqref{1.1} originating from $\theta_0$ provided by Proposition \ref{local}. Then
\begin{eqnarray*}
T\geq \frac{1}{C\|\theta_0\|^{1-\frac{n+4\alpha}{2s}}_{L^2(\mathbb{T}^n)}\|\theta_0\|^{\frac{n+4\alpha}{2s}}_{\dot{H}^s(T^n)}},
\end{eqnarray*}
where $C$ is the constant appearing in \eqref{3.8}.
\end{coro}
\textbf{Proof.}
Without loss of generality, let $[0,T)$ be the maximal time interval.
By \eqref{3.8} and \eqref{L2-decreasing} in Lemma \ref{maximum}, we derive that
\begin{eqnarray*}
\frac{d}{dt}\|\theta(t)\|_{\dot{H}^s}
\leq C\|\theta_0\|^{1-\frac{n+4\alpha}{2s}}_{L^2}\|\theta(t)\|^{1+\frac{n+4\alpha}{2s}}_{\dot{H}^s}.
\end{eqnarray*}
Solving this differential inequality directly yields that
\begin{eqnarray*}
\|\theta(t)\|_{\dot{H}^s}
\leq \frac{\|\theta_0\|_{\dot{H}^s}}{\Big[1-C\Big(\frac{n}{2s}+\frac{2\alpha}{s}\Big)\|\theta_0\|^{1-\frac{n+4\alpha}{2s}}_{L^2}\|\theta_0\|^{\frac{n+4\alpha}{2s}}_{\dot{H}^s}t\Big]^{\frac{2s}{n+4\alpha}}},
\end{eqnarray*}
which along with \eqref{L2-decreasing} and $s>\frac{n}{2}+1\geq\frac{n}{2}+2\alpha$ implies that $\|\theta(t)\|_{H^s}$ cannot blow up until
\begin{eqnarray*}
T_0\triangleq\frac{1}{C\|\theta_0\|^{1-\frac{n+4\alpha}{2s}}_{L^2}\|\theta_0\|^{\frac{n+4\alpha}{2s}}_{\dot{H}^s}}
<\frac{1}{C\Big(\frac{n}{2s}+\frac{2\alpha}{s}\Big)\|\theta_0\|^{1-\frac{n+4\alpha}{2s}}_{L^2}\|\theta_0\|^{\frac{n+4\alpha}{2s}}_{\dot{H}^s}}.
\end{eqnarray*}
This fact shows that $T\geq T_0$. The proof of Corollary \ref{lower-bound-time} is then completed.
%Therefore, for any $t\in[0,T_0]$ with
%\begin{eqnarray*}
%T_0=\frac{1}{\frac{4sC}{n+4\alpha}\|\theta_0\|^{1-\frac{n+4\alpha}{2s}}_{L^2}\|\theta_0\|^{\frac{n+4\alpha}{2s}}_{\dot{H}^s}},
%\end{eqnarray*}
%we have
%\begin{eqnarray*}
%\|\theta(t)\|_{\dot{H}^s}
%\leq2^{\frac{2s}{n+4\alpha}}\|\theta_0\|_{\dot{H}^s},
%\end{eqnarray*}
%which implies $T\geq T_0$. The proof of Corollary \ref{lower-bound-time} is then completed.
\section{Global regularity for the subcritical and critical dissipation}
In this section, we will prove Theorem \ref{the-1}, which is on the global regularity of the solutions to \eqref{1.1} with critical and subcritical dissipation.
\par
\textbf{Proof of Theorem \ref{the-1}.}
By Proposition \ref{BKM},
%virtue of the Beale-Kato-Majda blow-up criterion, i.e., Proposition \ref{BKM},
it is sufficient to obtain a control of $\|\nabla \theta\|_{L^\infty}$ uniformly in time.
%By Beale-Kato-Majda criterion in Proposition \ref{BKM}, it is equivalent to obtain a control of $\|\nabla\theta(t)\|_{L^\infty}$ uniformly in time. To do this, we will utilize the nonlinear maximum principle presented by Constantin and Vicol in \cite{[Constantin-Vicol]}.
%Without loss of generality, we can assume that $\theta$ is smooth on $[0,T)$. By Lemma \ref{maximum}, $\|\theta(\cdot,t)\|_{L^\infty}\leq\|\theta_0\|_{L^\infty}$, and hence a suitable a priori estimate on $\|\nabla\theta(\cdot,t)\|_{L^\infty}$
%implies that $\theta$ is in fact a smooth solution.
For this purpose, by \eqref{4.1-gradient}, Lemmas \ref{pointwise-property} and \ref{maximum-lower-bound}, we can derive that
%Applying $\nabla$ to \eqref{1.1} and multiplying by $\nabla\theta$ yield that
%\begin{eqnarray}\label{4.1-gradient}
%\frac{1}{2}(\partial_t+u\cdot\nabla)|\nabla\theta|^2
%+\nabla \theta\cdot\Lambda^\gamma\nabla \theta+\sum_{1\leq i,j\leq n} \partial_iu^j\partial_i\theta\partial_j\theta=0,
%\end{eqnarray}
%\begin{eqnarray}\label{4.1-gradient}
%\frac{1}{2}(\partial_t+u\cdot\nabla)|\nabla\theta|^2
%+\nabla \theta\cdot\Lambda^\gamma\nabla \theta+\nabla u:\nabla\theta\cdot\nabla\theta=0,
%\end{eqnarray}
%which follows from Lemmas \ref{pointwise-property} and \ref{maximum-lower-bound}, the $L^\infty$ maximum principle for $\theta$ that
\begin{eqnarray}\label{4.1}
\frac{1}{2}(\partial_t+u\cdot\nabla+\Lambda^\gamma)|\nabla\theta|^2
+c_1\frac{|\nabla \theta|^{2+\gamma}}{\|\theta_0\|^\gamma_{L^\infty}}+\frac{1}{4}D_\gamma(\nabla \theta)
\leq|\nabla u||\nabla\theta|^2.
\end{eqnarray}
Here, for simplicity, we consider $\kappa=1$.
%Let $\chi$ be a smooth radially nondecreasing cut-off function that vanishes for $|x|\leq\frac12$ and is equal to 1 for $|x|\geq1$, and such that its gradient satisfies $\|\nabla\chi\|_{L^\infty}\leq2$.
We proceed to estimate the term $\nabla u$ in \eqref{4.1}, which will be divided into two different cases, i.e., \eqref{1.1} with sub-critical and critical dissipation.

\textbf{Subcritical dissipation }: $2\alpha<\gamma<2$.
Similar to \eqref{1.2}, we split $\nabla u$ into two pieces:
\begin{eqnarray*}
\begin{split}
\nabla u(x,t)&=\nabla\Lambda^{-2+2\alpha}(\nabla\theta)(x,t)=\int_{\R^{n}}\frac{x-y}{|x-y|^{n+2\alpha}}(\nabla\theta(y,t)-\nabla\theta(x,t))dy\\
&=\underbrace{\int_{\R^{n}}\Big[1-\chi\Big(\frac{x-y}{\rho}\Big)\Big]F(x,y,t)dy}_{J_1}
+\underbrace{\int_{\R^{n}}\chi\Big(\frac{x-y}{\rho}\Big)F(x,y,t)dy}_{J_2},
\end{split}
\end{eqnarray*}
where $\rho=\rho(x,t)>0$ will be chosen later and
\begin{eqnarray}\label{F-definition-4.3}
F(x,y,t)=\frac{x-y}{|x-y|^{n+2\alpha}}(\nabla\theta(y,t)-\nabla\theta(x,t))
\end{eqnarray}
Also, the constant depending on $n$ and $\alpha$ before the integral has been omitted.

By H\"{o}lder's inequality, $\alpha\in(0,\frac12]$ and \eqref{definition-dissipation-remainder}, $J_1$ can be bounded by
\begin{eqnarray}\label{4.2}
\begin{split}
|J_1|
&\leq\int\limits_{|x-y|\leq\rho}\frac{|\nabla\theta(y,t)-\nabla\theta(x,t)|}{|x-y|^{n-1+2\alpha}}dy\\
&\leq\Big[\int\limits_{|y|\leq\rho}\frac{dy}{|y|^{n-2+4\alpha-\gamma}}\Big]^{\frac12}\Big[\int\limits_{|x-y|\leq\rho}\frac{|\nabla\theta(y,t)-\nabla\theta(x,t)|^2}{|x-y|^{n+\gamma}}\Big]^{\frac12}\\
&\leq c_2(n,\alpha,\gamma)\rho^{1-2\alpha+\frac{\gamma}{2}}(D_\gamma(\nabla \theta)(x))^{\frac12}.
\end{split}
\end{eqnarray}
By integration by parts, H\"{o}lder's inequality, the choice of cut-off function $\chi$ and Lemma \ref{maximum}, we can derive that
\begin{eqnarray}\label{4.3}
\begin{split}
|J_2|
&\lesssim \int_{\R^{n}}\Big[\frac{1}{\rho}\frac{|\nabla\chi(\frac{x-y}{\rho})|}{|x-y|^{n-1+2\alpha}}+\frac{|\chi(\frac{x-y}{\rho})|}{|x-y|^{n+2\alpha}}\Big]|\theta(y,t)-\theta(x,t)|dy\\
%&\leq&C_{n,\alpha}\Big[L^{-1}\int_{\R^{n}}\frac{|\nabla\chi(\frac{x-y}{L})|}{|x-y|^{n-1+2\alpha}}dy+\int_{\R^{n}}\frac{|\chi(\frac{x-y}{L})|}{|x-y|^{n+2\alpha}}dy \Big]\\
%&\leq&2\|\theta(t)\|_{L^\infty}\Big[\|\nabla\chi\|_{L^\infty}\int_{\frac{\rho}{2}<|z|\leq \rho}\frac{dz}{\rho|z|^{n-1+2\alpha}}
%+\|\chi\|_{L^\infty}\int\limits_{|z|>\frac{\rho}{2}}\frac{dz}{|z|^{n+2\alpha}} \Big]\nonumber\\
&\leq2\|\theta(t)\|_{L^\infty}\Big[\int\limits_{\frac{\rho}{2}<|z|\leq \rho}\frac{2dz}{\rho|z|^{n-1+2\alpha}}
+\int\limits_{|z|>\frac{\rho}{2}}\frac{dz}{|z|^{n+2\alpha}} \Big]
\leq c_3(n,\alpha)\frac{\|\theta_0\|_{L^\infty}}{\rho^{2\alpha}}.
\end{split}
\end{eqnarray}
Gathering up \eqref{4.2} and \eqref{4.3} gives that
\begin{eqnarray*}
|\nabla u(x,t)|
\leq c_2\rho^{1-2\alpha+\frac{\gamma}{2}}(D_\gamma(\nabla \theta)(x))^{\frac12}
+c_3\frac{\|\theta_0\|_{L^\infty}}{\rho^{2\alpha}},
\end{eqnarray*}
which follows form Young's inequality and \eqref{4.1} that
\begin{eqnarray*}
\frac{1}{2}(\partial_t+u\cdot\nabla+\Lambda^\gamma)|\nabla\theta|^2
+c_1\frac{|\nabla \theta|^{2+\gamma}}{\|\theta_0\|^\gamma_{L^\infty}}
\leq c^2_2\rho^{2-4\alpha+\gamma}|\nabla\theta|^4
+c_3\frac{\|\theta_0\|_{L^\infty}}{\rho^{2\alpha}}|\nabla\theta|^2.
\end{eqnarray*}
In this inequality, we explicitly choose $\rho>0$ such that
\begin{eqnarray*}
\frac{c_1}{2}\frac{|\nabla \theta|^{2+\gamma}}{\|\theta_0\|^\gamma_{L^\infty}}
=c^2_2\rho^{2-4\alpha+\gamma}|\nabla\theta|^4,~{\rm i.e.,}~\rho=\Big[\frac{c_1}{2c^2_2}\Big]^{\frac{1}{2-4\alpha+\gamma}}\frac{|\nabla\theta|^{\frac{\gamma-2}{2-4\alpha+\gamma}}}{\|\theta_0\|^{\frac{\gamma}{2-4\alpha+\gamma}}_{L^\infty}}
\end{eqnarray*}
to yield that
\begin{eqnarray}\label{4.4}
&&\frac{1}{2}(\partial_t+u\cdot\nabla+\Lambda^\gamma)|\nabla\theta|^2
\leq c_3\Big[\frac{2c^2_2}{c_1}\Big]^{\frac{2\alpha}{2-4\alpha+\gamma}}\|\theta_0\|^{1+\frac{2\alpha\gamma}{2-4\alpha+\gamma}}_{L^\infty}|\nabla\theta|^{\frac{2(1-\alpha)(2+\gamma)}{2-4\alpha+\gamma}}
-\frac{c_1}{2}\frac{|\nabla\theta|^{2+\gamma}}{\|\theta_0\|^\gamma_{L^\infty}},\nonumber\\
&&\quad
=\frac{c_1|\nabla\theta|^{\frac{2(1-\alpha)(2+\gamma)}{2-4\alpha+\gamma}}}{2\|\theta_0\|^\gamma_{L^\infty}}
\Big[\frac{2c_3}{c_1}\Big(\frac{2c^2_2}{c_1}\Big)^{\frac{2\alpha}{2-4\alpha+\gamma}}\|\theta_0\|^{\frac{(\gamma+2)(\gamma+1-2\alpha)}{2-4\alpha+\gamma}}_{L^\infty}
-|\nabla\theta|^{\frac{(2+\gamma)(\gamma-2\alpha)}{2-4\alpha+\gamma}}\Big].
\end{eqnarray}
For further arguments, we set
\begin{eqnarray*}
c_4=\Big[\frac{2c_3}{c_1}\Big]^{\frac{2-4\alpha+\gamma}{(2+\gamma)(\gamma-2\alpha)}}
\Big[\frac{2c^2_2}{c_1}\Big]^{\frac{2\alpha}{(2+\gamma)(\gamma-2\alpha)}},
\end{eqnarray*}
and introduce a nondecreasing $C^2$ convex function $\phi:[0,\infty)\rightarrow[0,\infty)$ satisfying
\begin{eqnarray*}
\phi(r)
=
\begin{cases}
0,
&  \mbox{{\rm for} $0\leq r\leq\max\Big\{\|\nabla\theta_0\|^2_{L^\infty},c^2_4
\|\theta_0\|^{2+\frac{2}{\gamma-2\alpha}}_{L^\infty}\Big\},$ }\\
{\rm positive},
& \mbox{{\rm otherwise}. } \\
\end{cases}
\end{eqnarray*}
By virtue of the convexity of $\phi$, we see that (see e.g., \cite{[Constantin]})
\begin{eqnarray*}
\phi'(|\nabla\theta|^2)\Lambda^\gamma(|\nabla\theta|^2)\geq\Lambda^\gamma(\phi(|\nabla\theta|^2)).
\end{eqnarray*}
Thus, we may multiply \eqref{4.4} by $\phi'(|\nabla\theta|^2)$ and obtain that
\begin{eqnarray*}
\begin{split}
&(\partial_t+u\cdot\nabla+\Lambda^\gamma)\phi(|\nabla\theta|^2)\nonumber\\
\leq&
\frac{c_1\phi'(|\nabla\theta|^2)|\nabla\theta|^{\frac{2(1-\alpha)(2+\gamma)}{2-4\alpha+\gamma}}}{\|\theta_0\|^\gamma_{L^\infty}}
\Big[\frac{2c_4}{c_1}\Big[\frac{2c_3}{c_1}\Big]^{\frac{2\alpha}{2-4\alpha+\gamma}}\|\theta_0\|^{\frac{(\gamma+2)(\gamma+1-2\alpha)}{2-4\alpha+\gamma}}_{L^\infty}
-|\nabla\theta|^{\frac{(2+\gamma)(\gamma-2\alpha)}{2-4\alpha+\gamma}}\Big]\leq0,
\end{split}
\end{eqnarray*}
%since $\phi'(|\nabla\theta|^2)=0$ for $|\nabla\theta|<c_5
%\|\theta_0\|^{1+\frac{1}{\gamma-2\alpha}}_{L^\infty}$.
which follows from $\alpha\in(0,\frac12]$ and $\gamma>2\alpha$ that
$\phi(|\nabla\theta|^2)=0$ a.e.. Indeed, we can show that
\begin{eqnarray*}
M(t)=:\displaystyle\max_{x\in\mathbb{T}^n} \phi(|\nabla\theta|^2)(x,t)=\phi(|\nabla\theta|^2)(x_t,t)
\end{eqnarray*}
is almost differentiable in time. Also, we can derive that
\begin{eqnarray*}
M'(t)
\leq-C_{n,\gamma}P.V.\int_{\R^n}\frac{\phi(|\nabla\theta|^2)(x_t,t)-\phi(|\nabla\theta|^2)(y,t)}{|x_t-y|^{n+\gamma}}dy
\leq0,
\end{eqnarray*}
which implies that
\begin{eqnarray*}
0\leq M(t)\leq M(0)=\max \phi(|\nabla\theta_0|^2)(x)=0.
\end{eqnarray*}
%where
%\begin{eqnarray*}
%c_5=\Big(\frac{2c_{4}}{c_1}\Big)^{\frac{2-4\alpha+\gamma}{(2+\gamma)(\gamma-2\alpha)}}
%\Big(\frac{2c_{3}}{c_1}\Big)^{\frac{2\alpha}{(2+\gamma)(\gamma-2\alpha)}}.
%\end{eqnarray*}
This along with the choice of $\phi$ yields that, for a.e. $x$ and all $t\in[0,T]$,
\begin{eqnarray*}
|\nabla\theta(x,t)|\leq\max\Big\{\|\nabla\theta_0\|_{L^\infty},c_4
\|\theta_0\|^{1+\frac{1}{\gamma-2\alpha}}_{L^\infty}\Big\},
\end{eqnarray*}
which is the desired estimate. Next we turn to the proof of critical part of Theorem \ref{the-1}.

\textbf{Critical dissipation}: $\gamma=2\alpha$.
Following \cite{[Constantin-Vicol]}, the proof is split in two steps. The first step shows that a solution to \eqref{1.1} having ``only small shocks" is regular, while
the second step shows that if the initial data has only small shocks, then so does the solution for all later time.
To be precise, we begin with the definition of Only Small Shocks.
\begin{definition}\label{definition-4.1}
%{\rm (Only Small Shocks)}
Let $\delta>0$ and $t>0$. We say $\theta(x,t)$ has the $OSS_\delta$ property, if there exists an $L>0$ such that
\begin{eqnarray*}
\sup_{\{(x,y):|x-y|<L\}}|\theta(x,t)-\theta(y,t)|\leq\delta.
\end{eqnarray*}
Moreover, for $T>0$, we say $\theta(x,t)$ has the uniform $OSS_\delta$ property on $[0,T]$, if
there exists an $L>0$ such that
\begin{eqnarray}\label{uniform}
\sup_{\{(x,y,t):|x-y|<L,0\leq t\leq T\}}|\theta(x,t)-\theta(y,t)|\leq\delta.
\end{eqnarray}
\end{definition}
The proof will be divided into a sequence of lemmas. The first states that the uniform $OSS_\delta$ property implies regularity of the solution:
\begin{lemma}\label{From-Only-Small-Shocks-to-regularity}
Let $n\geq2$, $\alpha\in(0,\frac12]$ and $\gamma=2\alpha$. For nonnegative initial data $\theta_0\in H^s(\R^n)$ with $s>\frac{n}{2}+1$, there exists a $\delta=\delta(\|\theta_0\|_{L^\infty})>0$ such that if $\theta$ is a solution to \eqref{1.1} with the uniform $OSS_{\delta}$ property on $[0,T]$,
then it is a smooth solution on $[0,T]$. Moreover
\begin{eqnarray*}
\sup_{t\in[0,T]}\|\nabla\theta(t)\|_{L^\infty}\leq C(L,\|\theta_0\|_{L^\infty},\|\nabla\theta_0\|_{L^\infty}).
\end{eqnarray*}
\end{lemma}
\textbf{Proof}. By hypothesis, there exists an $L>0$ such that \eqref{uniform} holds, where $\delta>0$ will be chosen later. By \eqref{1.2} and \eqref{F-definition-4.3}, we split the gradient of the velocity into three parts as
\begin{eqnarray*}
\begin{split}
\nabla u(x)
&=\underbrace{\int\limits_{\R^{n}}\Big[1-\chi\Big(\frac{x-y}{\rho}\Big)\Big]F(x,y)dy}_{K_1}
+\underbrace{\int\limits_{\R^{n}}\chi\Big(\frac{x-y}{L}\Big)F(x,y)dy}_{K_3}\\
&\ \ \ \ \ \ \ \ \ \ \
+\underbrace{\int\limits_{\R^{n}}\Big[\chi\Big(\frac{x-y}{\rho}\Big)-\chi\Big(\frac{x-y}{L}\Big)\Big]
F(x,y)dy}_{K_2},
\end{split}
\end{eqnarray*}
where $\rho=\rho(x,t)\leq L$ will be determined later.
Similar to \eqref{4.2}, $K_1$ can be bounded by
\begin{eqnarray}\label{4.7}
|K_1|
\leq c_2\rho^{1-2\alpha+\frac{\gamma}{2}}(D_\gamma(\nabla \theta)(x))^{\frac12}.
\end{eqnarray}
Integrating by parts on $K_2$ and utilizing H\"{o}lder's inequality, $\rho\leq L$, \eqref{uniform} and Lemma \ref{maximum} to yield that
\begin{eqnarray}\label{J2-4.7}
\begin{split}
|K_2|
&\lesssim \int_{\R^{n}}\Big|\frac{1}{\rho}\nabla\chi\Big(\frac{x-y}{\rho}\Big)-\frac{1}{L}\nabla\chi\Big(\frac{x-y}{L}\Big)\Big|\frac{|\theta(y,t)-\theta(x,t)|}{|x-y|^{n-1+2\alpha}}dy\\
&\ \ \
+\int_{\R^{n}}\Big|\chi\Big(\frac{x-y}{\rho}\Big)-\chi\Big(\frac{x-y}{L}\Big)\Big|\frac{|\theta(y,t)-\theta(x,t)|}{|x-y|^{n+2\alpha}}dy\\
%&=&
%P.V.\int_{\R^{n}}\Big|\rho^{-1}\nabla\chi(\frac{x-y}{\rho})-L^{-1}\nabla\chi(\frac{x-y}{L})\Big|\frac{|\theta(y,t)-\theta(x,t)|}{|x-y|^{n-1+2\alpha}}dy\\
%&&\ \ \
%+C_{n,\alpha}P.V.\int_{\min\{\frac\rho2,\frac L2\}<|x-y|<\max\{\rho,L\}}\Big|\chi(\frac{x-y}{\rho})-\chi(\frac{x-y}{L})\Big|\frac{|\theta(y,t)-\theta(x,t)|}{|x-y|^{n+2\alpha}}dy\\
&\lesssim
\int\limits_{\frac{\rho}{2}<|x-y|\leq \rho}\frac{|\theta(y,t)-\theta(x,t)|}{\rho|x-y|^{n-1+2\alpha}}dy
+\int\limits_{\frac{L}{2}<|x-y|\leq L}\frac{|\theta(y,t)-\theta(x,t)|}{L|x-y|^{n-1+2\alpha}}dy\\
&\ \ \ \ \ \ \
+\int\limits_{\frac\rho2<|x-y|<L}\frac{|\theta(y,t)-\theta(x,t)|}{|x-y|^{n+2\alpha}}dy\\
&\leq\int\limits_{\frac{\rho}{2}<|z|\leq \rho}\frac{\delta dz}{\rho|z|^{n-1+2\alpha}}
+\int\limits_{\frac{L}{2}<|z|\leq L}\frac{2\|\theta_0\|_{L^\infty}dz}{L|z|^{n-1+2\alpha}}
+\int\limits_{\frac\rho2<|z|<L}\frac{\delta dz}{|z|^{n+2\alpha}}\\
&\leq C_{n,\alpha}\Big[\frac{\delta}{\rho^{2\alpha}}+\frac{\|\theta_0\|_{L^\infty}}{L^{2\alpha}}\Big].
\end{split}
\end{eqnarray}
Similar to \eqref{4.3}, we can deduce that
\begin{eqnarray*}
|K_3|\leq c_3\frac{\|\theta_0\|_{L^\infty}}{L^{2\alpha}},
\end{eqnarray*}
which along with \eqref{4.7} and \eqref{J2-4.7} implies that
\begin{eqnarray*}
|\nabla u(x,t)|
\leq c_2\rho^{1-2\alpha+\frac{\gamma}{2}}(D_\gamma(\nabla \theta)(x))^{\frac12}
+c_4\Big[\frac{\delta}{\rho^{2\alpha}}+\frac{\|\theta_0\|_{L^\infty}}{L^{2\alpha}}\Big].
\end{eqnarray*}
Inserting this estimate into \eqref{4.1} and utilizing Young's inequality yields that
\begin{eqnarray*}
\frac{1}{2}(\partial_t+u\cdot\nabla+\Lambda^{2\alpha})|\nabla\theta|^2
+c_1\frac{|\nabla \theta|^{2+\gamma}}{\|\theta_0\|^\gamma_{L^\infty}}
\leq
c^2_2\rho^{2-4\alpha+\gamma}|\nabla\theta|^4
+c_4\Big[\frac{\delta}{\rho^{2\alpha}}+\frac{\|\theta_0\|_{L^\infty}}{L^{2\alpha}}\Big]|\nabla\theta|^2.
\end{eqnarray*}
In this differential inequality, we choose $\rho\in(0,L]$ such that $c^2_2\rho^{2-4\alpha+\gamma}|\nabla\theta|^4\leq\frac{c_1}{4}\frac{|\nabla \theta|^{2+\gamma}}{\|\theta_0\|^\gamma_{L^\infty}}$, that is,
\begin{eqnarray*}
\begin{split}
\rho
&=
\min\Big\{\Big(\frac{c_1}{4c^2_2}\frac{|\nabla \theta|^{\gamma-2}}{\|\theta_0\|^\gamma_{L^\infty}}\Big)^{\frac{1}{2-4\alpha+\gamma}},L\Big\}\\
&=\Big(\frac{c_1}{4c^2_2}\frac{|\nabla \theta|^{\gamma-2}}{\|\theta_0\|^\gamma_{L^\infty}}\Big)^{\frac{1}{2-4\alpha+\gamma}},~{\rm whenever}~|\nabla\theta|\geq\frac{1}{L}\Big[\frac{c_1}{4c^2_2\|\theta_0\|^{\gamma}_{L^\infty}}\Big]^{\frac{1}{2-2\alpha}},
\end{split}
\end{eqnarray*}
and then choose $\delta>0$ such that $c_4\frac{\delta}{\rho^{2\alpha}}|\nabla\theta|^2=\frac{c_1}{4}\frac{|\nabla \theta|^{2+\gamma}}{\|\theta_0\|^\gamma_{L^\infty}}$, that is,
\begin{eqnarray}\label{4.11}
\delta
=\frac{c_1}{4c_4}\Big[\frac{c_1}{4c^2_2}\Big]^{\frac{2\alpha}{2-4\alpha+\gamma}}
\frac{|\nabla\theta|^{\frac{(\gamma+2)(\gamma-2\alpha)}{2-4\alpha+\gamma}
}}{\|\theta_0\|^{\frac{\gamma(\gamma-2\alpha+2)}{2-4\alpha+\gamma}}_{L^\infty}}
=\frac{\frac{c_1}{4c_4}\Big[\frac{c_1}{4c^2_2}\Big]^{\frac{\alpha}{1-\alpha}}}{\|\theta_0\|^{\frac{2\alpha}{1-\alpha}}_{L^\infty}}.
\end{eqnarray}
Then, when $|\nabla\theta|\geq\frac{1}{L}\Big[\frac{c_1}{4c^2_2\|\theta_0\|^{\gamma}_{L^\infty}}\Big]^{\frac{1}{2-2\alpha}}$, we can obtain that
\begin{eqnarray*}
\begin{split}
\frac{1}{2}(\partial_t+u\cdot\nabla+\Lambda^\gamma)|\nabla\theta|^2
&\leq c_4\frac{\|\theta_0\|_{L^\infty}}{L^{2\alpha}}|\nabla\theta|^2-\frac{c_1}{2}\frac{|\nabla \theta|^{2+\gamma}}{\|\theta_0\|^\gamma_{L^\infty}}\\
&=\frac{c_1|\nabla\theta|^2}{2\|\theta_0\|^\gamma_{L^\infty}}
\Big[\frac{2c_4\|\theta_0\|^{1+\gamma}_{L^\infty}}{c_1L^{2\alpha}}-|\nabla\theta|^\gamma\Big],
\end{split}
\end{eqnarray*}
which immediately implies that
\begin{eqnarray*}
(\partial_t+u\cdot\nabla+\Lambda^\gamma)|\nabla\theta(x,t)|^2\leq0,~{\rm whenever}~|\nabla\theta(x,t)|\geq C_\ast
\end{eqnarray*}
where
\begin{eqnarray*}
C_\ast=\max\Big\{\frac{1}{L}\Big[\frac{2c_4}{c_1}\Big]^{\frac{1}{\gamma}}\|\theta_0\|^{1+\frac1\gamma}_{L^\infty},\frac{1}{L}\Big[\frac{c_1}{4c^2_2\|\theta_0\|^{2\alpha}_{L^\infty}}\Big]^{\frac{1}{2-2\alpha}}\Big\}.
\end{eqnarray*}
Conducting the similar arguments as in the subcritical dissipation, we can derive that
\begin{eqnarray*}
\sup_{t\in[0,T]}\|\nabla\theta(\cdot,t)\|_{L^\infty}\leq C(L,\|\theta_0\|_{L^\infty},\|\nabla\theta_0\|_{L^\infty}),
\end{eqnarray*}
which concludes the proof of Lemma \ref{From-Only-Small-Shocks-to-regularity}.

We next show that the stability of only small shocks for \eqref{1.1} with the critical dissipation:
\begin{lemma}\label{Stability-of-Only-Small-Shocks}
Let $\alpha\in(0,\frac12]$ and $\gamma=2\alpha$. Let $\sigma>0$ and $T>0$ be arbitrary. If the initial data $\theta_0$ has $OSS_{\frac{\sigma}{4}}$ property, then a solution $\theta$ to \eqref{1.1} has uniform $OSS_{\sigma}$ property on $[0,T]$.
\end{lemma}
In order to prove Lemma \ref{Stability-of-Only-Small-Shocks}, we consider the finite difference of the solution $\theta$ to \eqref{1.1}:
for $x,h\in\mathbb{T}^n$, we define $\delta_h\theta(x):=\theta(x+h)-\theta(x)$.
Then we have
%\begin{eqnarray*}
%\delta_h\theta(x)=\theta(x+h)-\theta(x).
%\end{eqnarray*}
\begin{lemma}\label{difference-maximum-lower-bound}
Let $n\geq2$, $\gamma\in(0,2)$ and nonnegative $\theta_0\in H^s(\mathbb{T}^n)$ with $s>\frac{n}{2}+1$. Then, for all $x,h\in\mathbb{T}^n$,
\begin{eqnarray*}
D_\gamma(\delta_h\theta)(x)\geq C_{n,\gamma}\frac{|\delta_h\theta(x)|^{\gamma+2}}{|h|^\gamma\|\theta_0\|^\gamma_{L^\infty}}.
\end{eqnarray*}
%\begin{eqnarray*}
%D_\gamma(\delta_h\theta)(x)\geq\frac{c^{\gamma+1}_5}{2^{\gamma+1}c^\gamma_{7}}\cdot\frac{|\delta_h\theta(x)|^{\gamma+2}}{|h|^\gamma\|\theta_0\|^\gamma_{L^\infty}}.
%\end{eqnarray*}
\end{lemma}
\textbf{Proof}. By \eqref{definition-dissipation-remainder} and Lemma \ref{maximum}, we derive that
\begin{eqnarray*}
&&D_\gamma(\delta_h\theta)(x)
\geq c_{n,\gamma} P.V.\int_{\R^n}\frac{(\delta_h\theta(x)-\delta_h\theta(x+y))^2}{|y|^{n+\gamma}}\chi\Big(\frac{y}{\rho}\Big)dy\\
&\geq&c_{n,\gamma}|\delta_h\theta(x)|^2\int_{\R^n}\frac{\chi(\frac{y}{\rho})}{|y|^{n+\gamma}}dy
-2c_{n,\gamma}|\delta_h\theta(x)|\Big|\int_{\R^n}\frac{\delta_h\theta(x+y)}{|y|^{n+\gamma}}\chi\Big(\frac{y}{\rho}\Big)dy\Big|\\
&\geq&c_{n,\gamma}|\delta_h\theta(x)|^2\int_{|y|\geq\rho}\frac{dy}{|y|^{n+\gamma}}
-2c_{n,\gamma}|\delta_h\theta(x)|\Big|\int_{\R^n}\theta(x+y)\delta_{-h}\Big(\frac{\chi(\frac{y}{\rho})}{|y|^{n+\gamma}}\Big)dy\Big|\\
&\geq&
c_{n,\gamma}\frac{|\delta_h\theta(x)|^2}{\rho^\gamma}
-2c_{n,\gamma}\|\theta_0\|_{L^\infty}|\delta_h\theta(x)|\int_{\R^n}\Big|\delta_{-h}\Big(\frac{\chi(\frac{y}{\rho})}{|y|^{n+\gamma}}\Big)\Big|dy.
\end{eqnarray*}
For $\rho\geq 8|h|$, by the mean value theorem, we obtain that
\begin{eqnarray*}
\begin{split}
\int_{\R^n}\Big|\delta_{-h}\Big(\frac{\chi(\frac{y}{\rho})}{|y|^{n+\gamma}}\Big)\Big|dy
&=\int_{|y|>\frac{\rho}{4}}\Big|\delta_{-h}\Big(\frac{\chi(\frac{y}{\rho})}{|y|^{n+\gamma}}\Big)\Big|dy\\
&\lesssim \int_{|y|>\frac{\rho}{4}}\frac{|h|dy}{|y|^{n+\gamma+1}}\leq C_{n,\gamma}\frac{|h|}{\rho^{\gamma+1}}.
\end{split}
\end{eqnarray*}
Therefore, for $\rho\geq 8|h|$, we have
\begin{eqnarray*}
D_\gamma(\delta_h\theta)(x)
\geq
c_5\frac{|\delta_h\theta(x)|^2}{\rho^\gamma}
-c_{6}\frac{\|\theta_0\|_{L^\infty}|h|}{\rho^{\gamma+1}}|\delta_h\theta(x)|.
%&\geq&c_5\frac{|\delta_h\theta(x)|^2}{\rho^\gamma}
%-c_{7}\frac{\|\theta_0\|_{L^\infty}|h|}{\rho^{\gamma+1}}|\delta_h\theta(x)|,
\end{eqnarray*}
Finally, denoting $c_7=\max\{8c_5,c_6\}$ and choosing
\begin{eqnarray*}
\rho:=\frac{2c_7|h|\|\theta_0\|_{L^\infty}}{c_5|\delta_h\theta(x,t)|}\geq8|h|,
\end{eqnarray*}
immediately yield the desired lower bound. We finish the proof of Lemma \ref{difference-maximum-lower-bound}.

We also need an estimate on the finite difference of the velocity stated in following lemma.
\begin{lemma}\label{maximum-difference-velocity}
Let $n\geq2$, $\alpha\in(0,\frac12]$ and $\gamma=2\alpha$. Suppose that $\theta$ is a smooth solution to \eqref{1.1} for nonnegative initial data $\theta_0\in H^s(\mathbb{T}^n)$ with $s>\frac{n}{2}+1$. Then, for all $x,h\in\mathbb{T}^n$,
\begin{eqnarray*}
|\delta_hu(x)|
\leq
\begin{cases}
C_{n,\gamma}+C_n\|\theta_0\|_{L^\infty}\ln^+ D_\gamma(\delta_h\theta)(x)+C'_n\|\theta_0\|_{L^2},
&  \mbox{{\rm for} $\alpha=\frac12,$ }\\
C_{n,\alpha}(\|\theta_0\|_{L^\infty}+\|\theta_0\|_{L^2}),
& \mbox{{\rm for} $\alpha\in(0,\frac12),$ } \\
\end{cases}
\end{eqnarray*}
\end{lemma}
\textbf{Proof}. By \eqref{1.2}, we have the expression
\begin{eqnarray}\label{velocity-difference}
\delta_hu(x)=C_{n,\alpha}P.V.\int_{\R^{n}}\frac{x-y}{|x-y|^{n+2\alpha}}\delta_h\theta(y)dy.
\end{eqnarray}
\textbf{Case 1}: $0<\alpha<\frac12$.
By H\"{o}lder's inequality and Lemma \ref{maximum}, we deduce that
\begin{eqnarray*}
\begin{split}
|\delta_hu(x)|
%&\leq& C_{n,\alpha}P.V.\int_{|x-y|<1}\frac{|\delta_h\theta(y)|dy}{|x-y|^{n-1+2\alpha}}
%+C_{n,\alpha}\int_{|x-y|\geq1}\frac{|\delta_h\theta(y)|dy}{|x-y|^{n-1+2\alpha}}\\
&\lesssim\|\delta_h\theta\|_{L^\infty}\int_{|z|<1}\frac{dz}{|z|^{n-1+2\alpha}}+\|\delta_h\theta\|_{L^2}\Big[\int_{|z|\geq1}\frac{dz}{|z|^{2(n-1+2\alpha)}}\Big]^{\frac12}\\
&\leq C_{n,\alpha}(\|\theta_0\|_{L^\infty}+\|\theta_0\|_{L^2}).
\end{split}
\end{eqnarray*}
\textbf{Case 2}: $\alpha=\frac12$.
By \eqref{velocity-difference}, we divide $\delta_hu(x)$ into three parts and it can be bounded by
\begin{eqnarray*}
|\delta_hu(x)|
%&=&C_{n,\alpha}\Big|P.V.\int_{\R^{n}}\frac{(x-y)\delta_h\theta(y)}{|x-y|^{n+1}}dy\Big|\\
%&\leq&C_{n,\alpha}\Big|P.V.\int_{|x-y|\leq\varrho}\frac{(x-y)(\delta_h\theta(y)-\delta_h\theta(x))}{|x-y|^{n+1}}dy\Big|\\
%&&\ \ \ \ \ \ \ \
%+C_{n,\alpha}\Big|\Big(\int_{\varrho<|x-y|\leq R}+\int_{|x-y|>R}\Big)\frac{(x-y)\delta_h\theta(y)}{|x-y|^{n+1}}dy\Big|\\
\lesssim \underbrace{\int_{|x-y|\leq\varrho}\frac{|\delta_h\theta(y)-\delta_h\theta(x)|}{|x-y|^{n}}dy}_{\mathcal{I}_1}
+\Big[\underbrace{\int_{\varrho<|x-y|\leq R}}_{\mathcal{I}_2}+\underbrace{\int_{|x-y|>R}}_{\mathcal{I}_3}\Big]\frac{|\delta_h\theta(y)|}{|x-y|^{n}}dy,
\end{eqnarray*}
where $0<\varrho\leq R$ will be chosen later.

By H\"{o}lder's inequality and \eqref{definition-dissipation-remainder}, $\mathcal{I}_1$ can be bounded as
\begin{eqnarray*}
\mathcal{I}_1
\leq\Big[\int\limits_{|z|\leq\varrho}\frac{dz}{|z|^{n-\gamma}}\Big]^{\frac12}\Big[\int\limits_{|x-y|\leq\varrho}\frac{|\delta_h\theta(y)-\delta_h\theta(x)|^2}{|x-y|^{n+\gamma}}dy\Big]^{\frac12}
\leq C_{n,\gamma}\varrho^{\frac\gamma2}(D_\gamma(\delta_h\theta)(x))^{\frac12}.
\end{eqnarray*}
Furthermore, by H\"{o}lder's inequality and Lemma \ref{maximum}, we have
\begin{eqnarray*}
\begin{split}
\mathcal{I}_2+\mathcal{I}_3
&\lesssim\|\theta(t)\|_{L^\infty}\int_{\varrho<|z|\leq R}\frac{dz}{|z|^{n}}+\|\delta_h\theta(t)\|_{L^2}\Big[\int_{|z|>R}\frac{dy}{|z|^{2n}}\Big]^{\frac12}\\
&\leq C_{n}\|\theta_0\|_{L^\infty}\ln\frac{R}{\varrho}+C'_n\|\theta_0\|_{L^2}R^{-\frac n2}.
\end{split}
\end{eqnarray*}
%We finally estimate $A_3$ as follows
%\begin{eqnarray*}
%A_3
%&\leq&\|\delta_h\theta(t)\|_{L^2}\Big(\int_{|z|>R}\frac{dy}{|z|^{2n}}\Big)^{\frac12}\\
%&\leq&C_n\|\theta_0\|_{L^2}R^{-\frac n2}.
%\end{eqnarray*}
Therefore, we obtain that
\begin{eqnarray*}
|\delta_hu(x)|
\leq C_{n,\gamma}\varrho^{\frac12}(D_\gamma(\delta_h\theta)(x))^{\frac12}+C_{n}\|\theta_0\|_{L^\infty}\ln\frac{R}{\varrho}+C'_n\|\theta_0\|_{L^2}R^{-\frac n2}.
\end{eqnarray*}
%We distinguish between the cases $D_\gamma(\delta_h\theta)(x)\leq1$ and $D_\gamma(\delta_h\theta)(x)>1$.
If $D_\gamma(\delta_h\theta)(x)\leq1$, we then take $\varrho=R=1$ to obtain that
\begin{eqnarray*}
|\delta_hu(x)|\leq C_{n,\gamma}+C'_n\|\theta_0\|_{L^2}.
\end{eqnarray*}
If $D_\gamma(\delta_h\theta)(x)>1$, we then choose $\rho=\frac{1}{D_\gamma(\delta_h\theta)(x)}$ and $R=1$ to get that
\begin{eqnarray*}
|\delta_hu(x)|\leq C_{n,\gamma}+C_n\|\theta_0\|_{L^\infty}\ln D_\gamma(\delta_h\theta)(x)+C'_n\|\theta_0\|_{L^2}.
\end{eqnarray*}
Summarizing these estimates, we can derive that
\begin{eqnarray*}
|\delta_hu(x)|\leq C_{n,\gamma}+C_n\|\theta_0\|_{L^\infty}\ln^+ D_\gamma(\delta_h\theta)(x)+C'_n\|\theta_0\|_{L^2}.
\end{eqnarray*}
%where we set $C_1=c_1+c_2\|\theta_0\|_{L^2}$ and $C_\infty=c_{\infty}\|\theta_0\|_{L^\infty}$.
We thus complete the proof of Lemma \ref{maximum-difference-velocity}.

With the help of Lemmas  \ref{difference-maximum-lower-bound} and \ref{maximum-difference-velocity}, we are ready to prove Lemma \ref{Stability-of-Only-Small-Shocks}.

\textbf{Proof of Lemma \ref{Stability-of-Only-Small-Shocks}.}
From $\eqref{1.1}_1$, the finite difference $\delta_h\theta$ satisfies
\begin{eqnarray}\label{solution-difference-equation}
\mathcal{L}_{\alpha,\gamma}(\delta_h\theta)=0,
\end{eqnarray}
where
\begin{eqnarray*}
\mathcal{L}_{\alpha,\gamma}=\partial_t+u\cdot\nabla_x+\delta_hu\cdot\nabla_h+\Lambda^\gamma.
\end{eqnarray*}
Following \cite{[Constantin-Vicol]}, we introduce $\Phi(h)=e^{-\Psi(|h|)}$,
multiply \eqref{solution-difference-equation} by $\Phi(h)\delta_h\theta$ and utilize Lemma \ref{pointwise-property} to obtain that
\begin{eqnarray}\label{4.13}
\mathcal{L}_{\alpha,\gamma}((\delta_h\theta)^2\Phi(h))
+D_{\gamma}(\delta_h\theta)\Phi(h)
=(\delta_h\theta)^2\delta_hu\cdot\nabla_h\Phi(h).
\end{eqnarray}
where $\Phi(h)$ and $\Psi(|h|)$ will be constructed at the end of the proof.

We define $\omega\triangleq(\delta_h\theta)^2\Phi(h)$, note that  $|\nabla_h\Phi|\leq\Phi(h)\Psi'(|h|)$,
%\begin{eqnarray*}
%|\nabla_h\Phi|\leq\Phi(h)\Psi'(|h|)
%\end{eqnarray*}
and so, by \eqref{4.13}, we have
\begin{eqnarray}\label{4.14}
\mathcal{L}_{\alpha,\gamma}\omega+\Phi(h)D_{\gamma}(\delta_h\theta)\leq\Psi'(|h|)|\delta_hu|\omega.
\end{eqnarray}

\textbf{Case 1}: $\alpha=\frac12$. In this case, $\gamma=1$.
By Lemma \ref{maximum-difference-velocity} and the inequality
\begin{eqnarray*}
bc\ln^+ a\leq\frac{a}{2}+bc\ln^+(2bc)~{\rm~for}~a,b,c>0,
\end{eqnarray*}
and Lemma \ref{maximum}, we can derive that
\begin{eqnarray*}
\begin{split}
\Psi'(|h|)\omega|\delta_hu|
%&\leq& \Psi'(|h|)\omega[c_8+c_9\log^+ D_\gamma(\delta_h\theta)(x)]\\
&\leq c_8\Psi'(|h|)\omega+\Phi(h)\Psi'(|h|)c_9(\delta_h\theta)^2\log^+ D_\gamma(\delta_h\theta)(x)\\
%&\leq&c_8\Psi'(|h|)\omega+\Phi(h)\Big[\frac{1}{2}D_\gamma(\delta_h\theta)(x)+\Psi'(|h|)c_9(\delta_h\theta)^2\log^+(2\Psi'(|h|)c_9(\delta_h\theta)^2)\Big]\\
&\leq\frac{\Phi(h)}{2}D_\gamma(\delta_h\theta)(x)+[c_8+c_9\log^+(2\Psi'(|h|)c_9(\delta_h\theta)^2)]\Psi'(|h|)\omega\\
&\leq\frac{\Phi(h)}{2}D_\gamma(\delta_h\theta)(x)+[c_8+c_9\ln^+(8\Psi'(|h|)C_n\|\theta_0\|^3_{L^\infty})]\Psi'(|h|)\omega,
\end{split}
\end{eqnarray*}
where $c_8=C_{n,\gamma}+C'_n\|\theta_0\|_{L^2}$ and $c_9=C_n\|\theta_0\|_{L^\infty}$.
It follows from \eqref{4.14} and Lemma \ref{difference-maximum-lower-bound} that
\begin{eqnarray*}
\mathcal{L}_{\alpha,\gamma}\omega+\frac{c_{10}\omega^{\frac32}(\Phi(h))^{-\frac12}}{|h|\|\theta_0\|_{L^\infty}}
%&\leq& (C_1+C_\infty\log^+(8\Psi'(|h|)c_\infty\|\theta_0\|^3_{L^\infty}))\Psi'\omega\\
\leq c_{11}\Psi'(|h|)[1+\ln(1+\Psi'(|h|))]\omega.
\end{eqnarray*}
Let us denote $q\triangleq\frac{c_{10}\sigma}{4c_{11}\|\theta_0\|_{L^\infty}}$,
pick $l>0$ and take a smooth function $\phi_l(y)=0$ for
$y\in[0,\frac l2]$, $\phi_l(y)=1$ for $y\geq l$, satisfying $0\leq\phi_l(y)\leq1$ for all $y\geq0$ and define
\begin{eqnarray*}
\Psi(y)=\int^y_0\frac{q\phi_l(x)}{x(1+\ln(1+\frac{q}{x}))}dx.
\end{eqnarray*}
It is clear that, for all $y>0$,
\begin{eqnarray*}
\Psi'(y)\leq \frac{q}{y(1+\ln(1+\frac{q}{y}))}.
\end{eqnarray*}
%It in turn follows that
%\begin{eqnarray*}
%\Psi'(y)(1+\ln(1+\Psi'(y)))\leq \frac{q}{y}
%\end{eqnarray*}
%for all $y>0$.
Then we deduce that
\begin{eqnarray*}
\begin{split}
\mathcal{L}_{\alpha,\gamma}\omega
&\leq c_{11}\frac{q\omega}{|h|}-\frac{c_{10}\omega^{\frac32}}{|h|\|\theta_0\|_{L^\infty}}(\Phi(h))^{-\frac12}\\
&\leq\frac{c_{10}\omega}{|h|\|\theta_0\|_{L^\infty}}\Big(\frac{\sigma }{4}-\omega^{\frac12}\Big)\leq0
\end{split}
\end{eqnarray*}
holds provided $\omega\geq(\frac\sigma4)^{2}$.
Because $\mathcal{L}_{\alpha,\gamma}$ has a weak maximum principle we have
\begin{eqnarray}\label{weak-maximum-principle}
\|\omega(t)\|_{L^\infty_{x,h}}\leq\|\omega_0\|_{L^\infty},~{\rm for~all}~t\in[0,T],~{\rm whenever}~\|\omega_0\|_{L^\infty}<\Big(\frac\sigma4\Big)^{2}.
\end{eqnarray}
%\begin{eqnarray*}
%\|\omega(\cdot,t;\cdot)\|_{L^\infty}\leq\|\omega_0\|_{L^\infty}
%\end{eqnarray*}
%for all $0\leq t\leq T$ if $\|\omega_0\|_{L^\infty}<(\frac\sigma4)^{2}$.
Now, if $\theta_0$ has the $OSS_{\frac{\sigma}{4}}$ property, by choosing $l>0$ small enough we can assure that $\|\omega_0\|_{L^\infty}<\frac{\sigma^2}{16}$,
%\begin{eqnarray*}
%\|\omega_0\|_{L^\infty}<\frac{\sigma^2}{16}
%\end{eqnarray*}
and so, for all $x,h$ and $t\leq T$,
\begin{eqnarray}\label{bound}
|\omega(x,t;h)|<\frac{\sigma^2}{16},~{\rm i.e.},~|\delta_h\theta(x,t)|\leq\frac{\sigma}{4}e^{\frac12\Psi(|h|)}.
\end{eqnarray}
%\begin{eqnarray*}
%|\omega(x,t;h)|<\frac{\sigma^2}{16},
%\end{eqnarray*}
%which means that
%\begin{eqnarray*}
%|\delta_h\theta(x,t)|\leq\frac{\sigma}{4}e^{\frac12\Psi(|h|)}.
%\end{eqnarray*}
Therefore, $|\delta_h\theta(x,t)|\leq\sigma$ for all $x,t\leq T$, and $|h|\leq \Psi^{-1}(4\ln2)$.

\textbf{Case 2}: $\alpha\in(0,\frac12)$. In this case, $\gamma=2\alpha\in(0,1)$. By Lemma \ref{maximum-difference-velocity} and \eqref{4.14},
we can derive that
\begin{eqnarray*}
\mathcal{L}_{\alpha,\gamma}\omega+\Phi(h)D_{\gamma}(\delta_h\theta)\leq c_{12}\Psi'(|h|)\omega.
\end{eqnarray*}
Let
\begin{eqnarray*}
\Psi(y)=\frac{c_{10}\sigma^{\gamma}y^{1-\gamma}}{(1-\gamma)c_{12}\|\theta_0\|^\gamma_{L^\infty}},
\end{eqnarray*}
and deduce from Lemma \ref{difference-maximum-lower-bound} that,
\begin{eqnarray*}
\begin{split}
\mathcal{L}_{\alpha,\gamma}\omega
&\leq c_{12}\Psi'(|h|)\omega-\frac{c_{10}\omega^{\frac\gamma2+1}}{|h|^\gamma\|\theta_0\|^\gamma_{L^\infty}}(\Phi(h))^{-\frac\gamma2}\\
&\leq\frac{c_{10}\omega}{|h|^\gamma\|\theta_0\|^\gamma_{L^\infty}}(\sigma^{\gamma}-\omega^{\frac\gamma2}),
\end{split}
\end{eqnarray*}
provided $\omega\geq\sigma^{2}$.
Similarly, by a weak maximum principle of $\mathcal{L}_{\alpha,\gamma}$, \eqref{weak-maximum-principle} also holds.
%Because $\mathcal{L}_{\alpha,\gamma}$ has a weak maximum principle we have
%\begin{eqnarray*}
%\|\omega(\cdot,t;\cdot)\|_{L^\infty}\leq\|\omega_0\|_{L^\infty}
%\end{eqnarray*}
%for all $0\leq t\leq T$ if $\|\omega_0\|_{L^\infty}<\sigma^{2}$.
If $\theta_0$ has the $OSS_{\frac{\sigma}{4}}$ property with the length $L$, then $\|\omega_0\|_{L^\infty(\mathbb{R}^n\times B_{L}(0))}<\frac{\sigma^2}{16}$,
%\begin{eqnarray*}
%\|\omega_0\|_{L^\infty(\mathbb{R}^n\times B_{L}(0))}<\frac{\sigma^2}{16}
%\end{eqnarray*}
and so, for all $x\in\R^n$, $|h|<L$ and $t\leq T$, \eqref{bound} holds.
%\begin{eqnarray*}
%|\omega(x,t;h)|<\frac{\sigma^2}{16},
%\end{eqnarray*}
%which means that
%\begin{eqnarray*}
%|\delta_h\theta(x,t)|\leq\frac{\sigma}{4}e^{\frac12\Psi(|h|)}.
%\end{eqnarray*}
Therefore, $|\delta_h\theta(x,t)|\leq\sigma$ for all $x,t\leq T$ and $|h|\leq \min\{\Psi^{-1}(4\ln2),L\}$.
We then complete the proof of Lemma \ref{Stability-of-Only-Small-Shocks}.

Combining Lemmas \ref{From-Only-Small-Shocks-to-regularity} and \ref{Stability-of-Only-Small-Shocks} concludes the proof of the critical part of Theorem \ref{the-1}.

\section{Smoothness for the supercritical case with $\gamma$ close to $2\alpha$}
In this section, we show Theorem \ref{the-2}.
We begin with a conditional regularity result of \eqref{1.1} with $\gamma\in(0,2\alpha)$, which shows that the control of the H\"{o}lder  seminorm of the solution is sufficient to obtain smoothness.
It is obtained by following the steps of Constantin and Vicol \cite{[Constantin-Vicol]} for the two-dimensional quasi-geostrophic equation without the use of Besov-space technique, which is utilized in \cite{[Constantin-Wu]} for this conditional regularity result.
%The main difference between \eqref{1.1} and 2-d QG is the incompressibility.
%Since we have nonnegative initial data, solutions of \eqref{1.1} are bounded in $L^2$ without the need for incompressibility, a core condition for showing the analogous result for 2-d QG (see also \cite{[Do]}).

%which will be utilized to obtain an eventual regularity result.
%Firstly, we obtain an eventual regularity result for the model \eqref{1.1} with the supercritical dissipation. This result is the main ingredient to prove Theorem \ref{the-2} (1).
\begin{lemma}\label{conditional-regularity}
Let $n\geq2$, $\alpha\in(0,\frac12]$, $\gamma\in(0,2\alpha)$ and $\Omega=\mathbb{T}^n$ or $\mathbb{R}^n$. Suppose that $\theta$ is a solution of \eqref{1.1} in $C([0,T),H^s(\Omega))$ for nonnegative initial data $\theta_0\in H^s(\Omega)$ with $s>\frac{n}{2}+1$.
%If $0<t_0<t_1\leq T$ and $\theta\in L^\infty((t_0,t_1);C^\beta(\mathbb{T}^n))$ with $2\alpha-\gamma<\beta<2\alpha$,
If $\theta\in L^\infty((t_0,t_1);C^\beta(\Omega))$ with $0<t_0<t_1\leq T$ and $\beta\in(2\alpha-\gamma,2\alpha)$,
then $\theta\in C^\infty(\Omega\times(t_0,t_1])$.
\end{lemma}
\textbf{Proof}.
We only consider the case of whole space for simplicity and clarity.
It is sufficient to prove that $\theta\in L^\infty((t_0,t_1];W^{1,\infty}(\mathbb{R}^n))$. This is sub-critical information which may be used to bootstrap and show that $\theta\in C^\infty(\mathbb{R}^n\times(t_0,t_1])$.

To this end, by \eqref{4.1-gradient}, Lemmas \ref{pointwise-property} and \ref{maximum-lower-bound}, we can derive that
\begin{eqnarray}\label{inequality-5.1}
\frac{1}{2}(\partial_t+u\cdot\nabla+\Lambda^\gamma)|\nabla\theta|^2
+\mathfrak{c}_1\frac{|\nabla \theta|^{2+\frac{\gamma}{1-\beta}}}{\|\theta\|^{\frac{\gamma}{1-\beta}}_{C^\beta}}+\frac{1}{4}D_\gamma(\nabla\theta)
\leq |\nabla u||\nabla\theta|^2.
\end{eqnarray}
%which follows from Lemma \ref{pointwise-property} that
%\begin{eqnarray*}
%\frac{1}{2}(\partial_t+u\cdot\nabla+\Lambda^\gamma)|\nabla\theta|^2
%+\frac12D_\gamma(\nabla\theta)(x)
%\leq |\nabla u(x,t)||\nabla\theta(x,t)|^2.
%\end{eqnarray*}
%It follows from Lemma \ref{Holder-lower-bound} that
%\begin{eqnarray*}
%\frac{1}{2}(\partial_t+u\cdot\nabla+\Lambda^\gamma)|\nabla\theta|^2
%+c_1\frac{|\nabla \theta(x)|^{2+\frac{\gamma}{1-\beta}}}{\|\theta\|^{\frac{\gamma}{1-\beta}}_{C^\beta}}+\frac14D_\gamma(\nabla\theta)(x)
%\leq |\nabla u(x,t)||\nabla\theta(x,t)|^2.
%\end{eqnarray*}
To bound $|\nabla u(x)|$, by \eqref{1.2} and \eqref{F-definition-4.3}, we split it into two pieces, according to an undetermined $\varrho>0$.
\begin{eqnarray*}
\nabla u(x,t)
%&=&C_{n,\alpha}P.V.\int_{\R^{n}}\frac{x-y}{|x-y|^{n+2\alpha}}\nabla\theta(y,t)dy\\
=\underbrace{\int_{\R^{n}}\Big[1-\chi\Big(\frac{x-y}{\varrho}\Big)\Big]F(x,y)dy}_{\mathcal{J}_1}
+\underbrace{\int_{\R^{n}}\chi\Big(\frac{x-y}{\varrho}\Big)F(x,y)dy}_{\mathcal{J}_2}.
\end{eqnarray*}
Similar to \eqref{4.2}, the inner piece $\mathcal{J}_1$ can be bounded as
\begin{eqnarray*}
|\mathcal{J}_1|
%&\leq& C_{n,\alpha}\int_{|x-y|\leq\varrho}\frac{|\nabla\theta(y,t)-\nabla\theta(x,t)|}{|x-y|^{n-1+2\alpha}}dy\\
\leq\mathfrak{c}_2\varrho^{\frac{\gamma}{2}+1-2\alpha}(D_\gamma(\nabla\theta)(x))^{\frac12}.
\end{eqnarray*}
For the outer piece, integrating by parts, we bound $\mathcal{J}_2$ as
\begin{eqnarray*}
\begin{split}
|\mathcal{J}_2|
&=\Big|\int_{\R^{n}}\nabla_y\Big[\chi\Big(\frac{x-y}{\rho}\Big)\frac{x-y}{|x-y|^{n+2\alpha}}\Big](\theta(y,t)-\theta(x,t))dy\Big|\\
&\leq\|\theta\|_{C^\beta}\int_{\R^{n}}\Big|\nabla_y\Big[\chi\Big(\frac{x-y}{\rho}\Big)\frac{x-y}{|x-y|^{n+2\alpha}}\Big]\Big||x-y|^\beta dy\\
&\lesssim\|\theta\|_{C^\beta}\int_{|z|\geq\frac{\varrho}{2}}\frac{dz}{|z|^{n+2\alpha-\beta}}
\leq \mathfrak{c}_3\frac{\|\theta\|_{C^\beta}}{\varrho^{2\alpha-\beta}}.
\end{split}
\end{eqnarray*}
Therefore, we obtain that
\begin{eqnarray*}
|\nabla u(x)|\leq \mathfrak{c}_2\varrho^{\frac{\gamma}{2}+1-2\alpha}(D_\gamma(\nabla\theta)(x))^{\frac12}
+\mathfrak{c}_3\frac{\|\theta\|_{C^\beta}}{\varrho^{2\alpha-\beta}}.
\end{eqnarray*}
%which follows from Young's inequality that
%\begin{eqnarray*}
%|\nabla u||\nabla\theta|^2
%\leq\frac{D_\gamma(\nabla\theta)(x)}{4}+c^2_1\varrho^{\gamma+2-4\alpha}|\nabla\theta(x,t)|^4+\frac{c_2\|\theta\|_{C^\beta}}{\varrho^{2\alpha-\beta}}|\nabla\theta|^2.
%\end{eqnarray*}
Substituting this estimate into \eqref{inequality-5.1} and utilizing Young's inequality yield that
\begin{eqnarray*}
\frac{1}{2}(\partial_t+u\cdot\nabla+\Lambda^\gamma)|\nabla\theta|^2
+\mathfrak{c}_1\frac{|\nabla \theta(x)|^{2+\frac{\gamma}{1-\beta}}}{\|\theta\|^{\frac{\gamma}{1-\beta}}_{C^\beta}}
\leq \mathfrak{c}^2_2\varrho^{\gamma+2-4\alpha}|\nabla\theta|^4+\mathfrak{c}_3\frac{\|\theta\|_{C^\beta}}{\varrho^{2\alpha-\beta}}|\nabla\theta|^2.
\end{eqnarray*}
In this equality, we take $\varrho>0$ such that
\begin{eqnarray*}
\mathfrak{c}^2_2\varrho^{\gamma+2-4\alpha}|\nabla\theta|^4=\mathfrak{c}_3\frac{\|\theta\|_{C^\beta}}{\varrho^{2\alpha-\beta}}|\nabla\theta|^2,
~{\rm i.e.,}~\varrho=\Big(\frac{\mathfrak{c}_3\|\theta\|_{C^\beta}}{\mathfrak{c}^2_2|\nabla\theta|^2}\Big)^{\frac{1}{\gamma+2-2\alpha-\beta}}
\end{eqnarray*}
to imply that
\begin{eqnarray*}
&&\frac{1}{2}(\partial_t+u\cdot\nabla+\Lambda^\gamma)|\nabla\theta|^2
\leq \mathfrak{c}_4\|\theta\|^{\frac{\gamma+2-4\alpha}{\gamma+2-2\alpha-\beta}}_{C^\beta}
|\nabla \theta|^{\frac{2(\gamma+2-2\beta)}{\gamma+2-2\alpha-\beta}}-\mathfrak{c}_1\frac{|\nabla \theta|^{2+\frac{\gamma}{1-\beta}}}{\|\theta\|^{\frac{\gamma}{1-\beta}}_{C^\beta}}\\
%&\leq&c_4M^{\frac{\gamma+2-4\alpha}{\gamma+2-2\alpha-\beta}}
%|\nabla \theta|^{\frac{2(\gamma+2-2\beta)}{\gamma+2-2\alpha-\beta}}
%-c_1\frac{|\nabla \theta|^{2+\frac{\gamma}{1-\beta}}}{M^{\frac{\gamma}{1-\beta}}}\\
&&\quad
\leq
\frac{\mathfrak{c}_1|\nabla \theta|^{\frac{2(\gamma+2-2\beta)}{\gamma+2-2\alpha-\beta}}}{\mathfrak{M}^{\frac{\gamma}{1-\beta}}}
\Big[\frac{\mathfrak{c}_4}{\mathfrak{c}_1}\mathfrak{M}^{\frac{\gamma}{1-\beta}+\frac{\gamma+2-4\alpha}{\gamma+2-2\alpha-\beta}}
-|\nabla\theta|^{\frac{\gamma}{1-\beta}+\frac{2(\beta-2\alpha)}{\gamma+2-2\alpha-\beta}}\Big],
\end{eqnarray*}
where $\mathfrak{M}:=\|\theta\|_{L^\infty_tC^\beta_x}$.
Note that, for $\beta>2\alpha-\gamma$,
\begin{eqnarray*}
\frac{\gamma}{1-\beta}+\frac{2(\beta-2\alpha)}{\gamma+2-2\alpha-\beta}
>\frac{\gamma}{1-(2\alpha-\gamma)}+\frac{2(2\alpha-\gamma-2\alpha)}{\gamma+2-2\alpha-(2\alpha-\gamma)}=0.
\end{eqnarray*}
Thus, we derive that
\begin{eqnarray*}
(\partial_t+u\cdot\nabla+\Lambda^\gamma)|\nabla\theta(x,t)|^2\leq0,
\end{eqnarray*}
provided
\begin{eqnarray*}
|\nabla\theta(x,t)|\geq\Big[\frac{\mathfrak{c}_4}{\mathfrak{c}_1}\mathfrak{M}^{\frac{\gamma}{1-\beta}+\frac{\gamma+2-4\alpha}{\gamma+2-2\alpha-\beta}}\Big]^{\frac{1}{\frac{\gamma}{1-\beta}+\frac{2(\beta-2\alpha)}{\gamma+2-2\alpha-\beta}}}.
\end{eqnarray*}
Then the maximum of $|\nabla\theta|$ can not exceed a certain constant which depends on $\mathfrak{M}$, showing that
$\nabla\theta\in L^\infty((t_0,t_1);L^\infty)$. We then complete the proof of Lemma \ref{conditional-regularity}.

Based on the conditional regularity result, we establish an eventual regularity result for \eqref{1.1} with the supercritical dissipation.
\begin{proposition}\label{eventual-regularity}
Let $n\geq2$, $\alpha\in(0,\frac12]$, $\gamma\in(0,2\alpha)$ and nonnegative $\theta_0\in H^s(\mathbb{T}^n)$ with $s>\frac{n}{2}+1$.
Consider $\beta\in(2\alpha-\gamma,2\alpha)$ and define
\begin{eqnarray}\label{5.1}
T^\ast_\alpha:=C\beta^{\frac{2\alpha}{2\alpha-\gamma}}\|\theta_0\|^{\frac{\gamma}{2\alpha-\gamma}}_{L^\infty},
\end{eqnarray}
where $C>0$ is a constant independent of $\beta$ and $\theta_0$.
If $\theta\in C([0,T),H^s(\mathbb{T}^n))$ is a solution to \eqref{1.1} and $T^\ast_\alpha<T<+\infty$, then $\theta\in C^\infty(\mathbb{T}^n\times(T^\ast_\alpha,T])$.
\end{proposition}
In view of Lemma \ref{conditional-regularity}, it is reduced to show that $\theta\in L^\infty((T^\ast_\alpha,T);C^\beta(\mathbb{T}^n))$.
%\begin{eqnarray*}
%\theta\in L^\infty((T^\ast_\alpha,T);C^\beta(\mathbb{T}^n)).
%\end{eqnarray*}
%with $\beta>2\alpha-\gamma$ and $T^\ast_\alpha$ is as in \eqref{5.1}.
For this purpose, we consider
%introduce $\delta_h\theta(x,t)=\theta(x+h,t)-\theta(x,t)$ and
\begin{eqnarray}\label{5.2}
\vartheta(x,t;h)\triangleq\frac{\delta_h\theta(x,t)}{(\eta^2(t)+|h|^2)^{\frac\beta2}},
\end{eqnarray}
%\begin{eqnarray}\label{5.2}
%\omega(x,t;h)=\frac{\delta_h\theta(x,t)}{(\eta^2(t)+|h|^2)^{\frac\beta2}},
%\end{eqnarray}
where $\eta:[0,\infty)\rightarrow[0,\infty)$ is a bounded decreasing differentiable function determined later. In order to control the $C^\beta$-seminorm of $\theta$, it is sufficient to estimate $\|\vartheta(t)\|_{L^\infty_{x,h}}$ when $\eta(t)=0$.
%Note that it is sufficient to estimate $\|\vartheta(t)\|_{L^\infty_{x,h}}$ when $\eta(t)=0$ in order to control the $C^\beta$-seminorm of $\theta$.
%We denote
%\begin{eqnarray*}
%\mathcal{L_{\alpha,\gamma}}=\partial_t+u\cdot\nabla_x+\delta_hu\cdot\nabla_h+\Lambda^\gamma.
%\end{eqnarray*}
%We first compute $\mathcal{L_{\alpha,\gamma}}\vartheta^2$ where $\mathcal{L_{\alpha,\gamma}}$ is the operator of the corresponding equation satisfied by $\delta_h\theta$.
%In order to find $\mathcal{L_{\alpha,\gamma}}$, taking the differences in \eqref{} evaluated in $x+h$ and $x$, it follows that
%\begin{eqnarray}
%\partial_t\delta_h\theta+u\cdot\nabla_x\delta_h\theta+\delta_hu\cdot\nabla_h\delta_h\theta+\Lambda^\gamma\delta_h\theta=0,
%\end{eqnarray}
%which gives that
%\begin{eqnarray*}
%\mathcal{L}_{\alpha,\gamma}=\partial_t+u\cdot\nabla_x+\delta_hu\cdot\nabla_h+\Lambda^\gamma
%\end{eqnarray*}
%with $u=\nabla\Lambda^{-2+2\alpha}\theta$.
By Lemma \ref{pointwise-property} and \eqref{solution-difference-equation}, a standard computation yields that
\begin{eqnarray}\label{5.3}
\mathcal{L}_{\alpha,\gamma}\vartheta^2+\frac{D_\gamma(\delta_h\theta)}{(\eta^2(t)+|h|^2)^{\beta}}
=-\frac{2\beta\eta'(t)\eta(t)}{\eta^2(t)+|h|^2}\vartheta^2-\frac{2\beta\delta_hu\cdot h}{\eta^2(t)+|h|^2}\vartheta^2.
\end{eqnarray}
Next we bound below the second term in the left-hand side of \eqref{5.3}.
\begin{lemma}\label{nonlinear-lower-bound-difference-of-solution}
Let $n\geq2$, $\alpha\in(0,\frac12]$, $0<\gamma_0\leq\gamma<2\alpha$ and $\beta\in(2\alpha-\gamma,2\alpha)$.
Then, for $x,h\in\mathbb{T}^n$,
\begin{eqnarray*}
\frac{D_\gamma(\delta_h\theta)(x)}{(\eta^2(t)+|h|^2)^{\beta}}
\geq
\frac{1}{\mathfrak{c}_0|h|^\gamma}\Big[\frac{|\vartheta(x,t;h)|}{\|\vartheta(t)\|_{L^\infty_{x,h}}}\Big]^{\frac{\gamma}{1-\beta}}|\vartheta(x,t;h)|^2.
\end{eqnarray*}
\end{lemma}
\textbf{Proof}. From the proof of Lemma \ref{difference-maximum-lower-bound}, we can see that, for $\rho\geq8|h|$,
\begin{eqnarray}\label{5.5}
&&D_\gamma(\delta_h\theta)(x)
-c_{n,\gamma}\frac{|\delta_h\theta(x)|^2}{\rho^\gamma}
\geq
-2c_{n,\gamma}|\delta_h\theta(x)|\Big|\int_{|y|>\frac{\rho}{4}}\delta_y\theta(x)\delta_{-h}\Big(\frac{\chi(\frac{y}{\rho})}{|y|^{n+\gamma}}\Big)dy\Big|\nonumber\\
&&\qquad\quad
\geq-2c_{n,\gamma}|\delta_h\theta(x)|\cdot C_{n,\gamma}|h|\int_{|y|>\frac{\rho}{4}}\frac{|\delta_y\theta(x)|}{|y|^{n+\gamma+1}}dy\nonumber\\
&&\qquad\quad
\geq-2c_{n,\gamma}|\delta_h\theta(x)|\cdot C_{n,\gamma}|h|\|\vartheta(t)\|_{L^\infty_{x,h}}\int_{|y|>\frac{\rho}{4}}\frac{(\eta^2(t)+|y|^2)^{\frac\beta2}}{|y|^{n+\gamma+1}}dy\nonumber\\
&&\qquad\quad
\geq-2c_{n,\gamma}|\delta_h\theta(x)|\cdot \underbrace{C_{n,\gamma}\frac{|h|\|\vartheta(t)\|_{L^\infty_{x,h}}}{\rho^\gamma}\Big[\frac{\eta^\beta(t)}{\rho}+\frac{1}{\rho^{1-\beta}}\Big]}_{\mathcal{K}}.
\end{eqnarray}
Now choose $\rho>0$ given by
\begin{eqnarray*}
\rho=\Big[\frac{8C_{n,\gamma}\|\vartheta(t)\|_{L^\infty_{x,h}}}{|\vartheta(x,t;h)|}\Big]^{\frac{1}{1-\beta}}|h|.
\end{eqnarray*}
Since we can assume $C_{n,\gamma}\geq1$ and $|\vartheta(x,t;h)|\leq\|\vartheta(t)\|_{L^\infty_{x,h}}$, $\rho\geq8^{\frac{1}{1-\beta}}|h|\geq8|h|$.
By the choice of $\rho$ and the inequality $a^{\frac{1}{1-\beta}}\leq a$ for all $a\leq1$,
we can bound $\mathcal{K}$ as
\begin{eqnarray*}
\begin{split}
|\mathcal{K}|
%&\leq&C_{n,\gamma}\frac{|h|\|\vartheta(x,t;\cdot)\|_{L^\infty}}{\rho^\gamma}\Big(\frac{\eta^\beta(t)}{\rho}+\frac{1}{\rho^{1-\beta}}\Big)\\
&\leq\frac{1}{\rho^\gamma}\Big[\frac{C_{n,\gamma}\eta^\beta(t)\|\vartheta(t)\|_{L^\infty_{x,h}}}{(8C_{n,\gamma})^{\frac{1}{1-\beta}}}\Big(\frac{|\vartheta(x,t;h)|}{\|\vartheta(t)\|_{L^\infty_{x,h}}}\Big)^{\frac{1}{1-\beta}}+\frac{|\vartheta(x,t;h)||h|^\beta}{8}\Big]\\
&\leq
\frac{\eta^\beta(t)+|h|^\beta}{8\rho^\gamma}|\vartheta(x,t;h)|
%=\frac{\eta^\beta(t)+|h|^\beta}{8\rho^\gamma}\frac{|\delta_h\theta(x,t)|}{(\eta^2(t)+|h|^2)^{\frac\beta2}}
\leq\frac{2^{1-\frac{\beta}{2}}(\eta^2(t)+|h|^2)^{\frac\beta2}}{8\rho^\gamma}|\vartheta(x,t;h)|\\
%\frac{|\delta_h\theta(x,t)|}{(\eta^2(t)+|h|^2)^{\frac\beta2}}
&=\frac{|\delta_h\theta(x,t)|}{2^{2+\frac{\beta}{2}}\rho^\gamma}\leq\frac{|\delta_h\theta(x,t)|}{4\rho^\gamma}.
\end{split}
\end{eqnarray*}
Substituting this estimate of $\mathcal{K}$ into \eqref{5.5}
and utilizing the choice of $\rho$, we arrive at
\begin{eqnarray*}
D_\gamma(\delta_h\theta)(x)
\geq c_{n,\gamma}\frac{|\delta_h\theta(x)|^2}{2\rho^\gamma}=
\frac{c_{n,\gamma}}{2(8C_{n,\gamma})^{\frac{\gamma}{1-\beta}}|h|^\gamma}\Big[\frac{|\vartheta(x,t;h)|}{\|\vartheta(x,t;\cdot)\|_{L^\infty}}\Big]^{\frac{\gamma}{1-\beta}}|\delta_h\theta(x)|^2,
\end{eqnarray*}
which along with \eqref{5.2} concludes the proof of Lemma \ref{nonlinear-lower-bound-difference-of-solution}.

We proceed to derive the differential equation that $\eta$ should solve to control the first term in the right-hand side of \eqref{5.3}.
%with a fraction of the nonlinear lower bound \eqref{}.
\begin{lemma}\label{differential-equation}
Let $n\geq2$, $\gamma_0>0$, $\gamma\in[\gamma_0,1)$ and $\beta\in(1-\gamma,1)$.
If $\eta(t)$ solves
\begin{eqnarray}\label{5.7}
\frac{d}{dt}\eta(t)=-\frac{1}{16\mathfrak{c}_0\beta}\eta^{1-\gamma}(t),
\end{eqnarray}
then
\begin{eqnarray*}
-\frac{2\beta\eta'\eta}{\eta^2+|h|^2}\vartheta^2\leq\frac{\vartheta^2}{8\mathfrak{c}_0|h|^\gamma}~{\rm for~ all}~ x,h\in \mathbb{T}^n,
\end{eqnarray*}
where $\mathfrak{c}_0$ is the constant appearing in Lemma \ref{nonlinear-lower-bound-difference-of-solution}.
\end{lemma}
\textbf{Proof}.
If $\eta$ solves \eqref{5.7}, then
\begin{eqnarray*}
-\frac{2\beta\eta'\eta}{\eta^2+|h|^2}\vartheta^2
=\frac{\eta^{2-\gamma}(t)}{8\mathfrak{c}_0(\eta^2+|h|^2)}\vartheta^2
\leq\frac{\vartheta^2}{8\mathfrak{c}_0(\eta^2+|h|^2)^{\frac{\gamma}{2}}}\leq\frac{\vartheta^2}{8\mathfrak{c}_0|h|^\gamma},
\end{eqnarray*}
which completes the proof of Lemma \ref{differential-equation}.

Following \cite{[CotiZelati-Vicol]}, we now consider the second term in the right-hand side of \eqref{5.3} to derive a suitable upper bound in terms of the dissipation.
We begin with the estimate of $\delta_hu$.
\begin{lemma}\label{bound-difference-velocity}
Let $n\geq2$, $\alpha\in(0,\frac12]$, $0<\gamma_0\leq\gamma<2\alpha$ and $\beta\in(2\alpha-\gamma,2\alpha)$. If $\rho\geq 4|h|$, then we have, for all $x,h\in\mathbb{T}^n$,
%we have the estimate
\begin{eqnarray*}
|\delta_hu(x)|\leq C\Big[\rho^{\frac{\gamma}{2}+1-2\alpha}(D_\gamma(\delta_h\theta)(x))^{\frac12}+|h|\|\vartheta\|_{L^\infty_{x,h}}\Big(\frac{\eta^\beta(t)}{\rho^{2\alpha}}+\frac{1}{\rho^{2\alpha-\beta}}\Big)\Big].
\end{eqnarray*}
\end{lemma}
\textbf{Proof}.
By \eqref{velocity-difference}, we divide $\delta_hu(x)$ into two parts
\begin{eqnarray*}
\delta_hu(x)
=\underbrace{\int\limits_{\R^{n}}\Big[1-\chi(\frac{y}{\rho})\Big]\frac{y[\delta_h\theta(x-y)-\delta_h\theta(x)]}{|y|^{n+2\alpha}}dy}_{\mathcal{K}_1}
+\underbrace{\int\limits_{\R^{n}}\chi(\frac{y}{\rho})\frac{y[\delta_h\theta(x-y)-\delta_h\theta(x)]}{|y|^{n+2\alpha}}dy}_{\mathcal{K}_2}.
\end{eqnarray*}
Similar to \eqref{4.2}, we can estimate $\mathcal{K}_1$ as
\begin{eqnarray}\label{5.8-K1}
|\mathcal{K}_1|\leq C\rho^{\frac{\gamma}{2}+1-2\alpha}(D_\gamma(\delta_h\theta)(x))^{\frac12}.
\end{eqnarray}
By the mean value theorem, we can derive that, for $\rho\geq8|h|$,
%\begin{eqnarray*}
%|\mathcal{K}_2|
%&=&\Big|\int_{\R^{n}}\chi\Big(\frac{y}{\rho}\Big)\frac{y}{|y|^{n+2\alpha}}\delta_h(\theta(x-y)-\theta(x))dy\Big|\\
%&=&\Big|\int_{\R^n}\delta_{-h}\Big(\chi\Big(\frac{y}{\rho}\Big)\frac{y}{|y|^{n+2\alpha}}\Big)(\theta(x-y)-\theta(x))dy\Big|\\
%&\leq&\int_{|y|>\frac{\rho}{4}}\Big|\delta_{-h}\Big(\chi\Big(\frac{y}{\rho}\Big)\frac{y}{|y|^{n+2\alpha}}\Big)\Big||\delta_{-y}\theta(x)|dy\\
%&\leq&C_{n,\alpha}\int_{|y|>\frac{\rho}{4}}\frac{|h|}{|y|^{n+2\alpha}}|\delta_{-y}\theta(x)|dy\\
%&\leq&C_{n,\alpha}|h|\|\vartheta\|_{L^\infty_{x,h}}\int_{|y|>\frac{\rho}{4}}\frac{(\eta^2(t)+|y|^2)^{\frac\beta2}}{|y|^{n+2\alpha}}dy\\
%&\leq&C_{n,\alpha}|h|\|\vartheta\|_{L^\infty_{x,h}}\int_{|y|>\frac{3}{4}\rho}\frac{\eta^\beta(t)+|y|^\beta}{|y|^{n+2\alpha}}dy\\
%&=&C_{n,\alpha}|h|\|\vartheta\|_{L^\infty_{x,h}}\Big[\frac{\eta^\beta(t)}{\rho^{2\alpha}}+\frac{1}{\rho^{2\alpha-\beta}}\Big],
%\end{eqnarray*}
\begin{eqnarray*}
&&|\mathcal{K}_2|
%&=&\Big|\int_{\R^{n}}\chi\Big(\frac{y}{\rho}\Big)\frac{y}{|y|^{n+2\alpha}}\delta_h(\theta(x-y)-\theta(x))dy\Big|\\
%&=&\Big|\int_{\R^n}\delta_{-h}\Big(\chi\Big(\frac{y}{\rho}\Big)\frac{y}{|y|^{n+2\alpha}}\Big)(\theta(x-y)-\theta(x))dy\Big|\\
\leq\int_{|y|>\frac{\rho}{4}}\Big|\delta_{-h}\Big(\chi\Big(\frac{y}{\rho}\Big)\frac{y}{|y|^{n+2\alpha}}\Big)\Big||\delta_{-y}\theta(x)|dy\\
&&\qquad\leq C_{n,\alpha}\int_{|y|>\frac{\rho}{4}}\frac{|h|}{|y|^{n+2\alpha}}|\delta_{-y}\theta(x)|dy\\
&&\qquad\leq C_{n,\alpha}|h|\|\vartheta\|_{L^\infty_{x,h}}\int_{|y|>\frac{\rho}{4}}\frac{(\eta^2(t)+|y|^2)^{\frac\beta2}}{|y|^{n+2\alpha}}dy\\
%&\leq&C_{n,\alpha}|h|\|\vartheta\|_{L^\infty_{x,h}}\int_{|y|>\frac{3}{4}\rho}\frac{\eta^\beta(t)+|y|^\beta}{|y|^{n+2\alpha}}dy\\
&&\qquad=C_{n,\alpha}|h|\|\vartheta\|_{L^\infty_{x,h}}\Big[\frac{\eta^\beta(t)}{\rho^{2\alpha}}+\frac{1}{\rho^{2\alpha-\beta}}\Big],
\end{eqnarray*}
which along with \eqref{5.8-K1} concludes the proof of Lemma \ref{bound-difference-velocity}.

In view of Lemma \ref{bound-difference-velocity}, we are able to compare the nonlinear term in \eqref{5.3} with the lower bound on the dissipation term given in Lemma \ref{nonlinear-lower-bound-difference-of-solution}.
\begin{lemma}\label{condition-on-the-initial-data}
Let $n\geq2$, $\alpha\in(0,\frac12]$, $0<\gamma_0\leq\gamma<2\alpha$ and $\beta\in(2\alpha-\gamma,2\alpha)$. Suppose that
\begin{eqnarray*}
\|\vartheta\|_{L^\infty_{x,h}}\leq\frac{4\|\theta_0\|_{L^\infty}}{\eta^\beta_0}.
\end{eqnarray*}
There exists a constant $\mathfrak{c}=\mathfrak{c}(\gamma_0)\geq1$ such that if we choose $\eta_0=(\mathfrak{c}\beta\|\theta_0\|_{L^\infty})^{\frac{1}{2\alpha-\gamma}}$,
%\begin{eqnarray*}
%\eta_0=\Big(c_1\beta\|\theta_0\|_{L^\infty}\Big)^{\frac{1}{2\alpha-\gamma}},
%\end{eqnarray*}
then, for every $x,h\in\mathbb{T}^n$ with $|h|\leq\eta_0$, it holds that
\begin{eqnarray*}
\frac{2\beta|h||\delta_h u|}{\eta^2+|h|^2}\vartheta^2
\leq
\frac{D_\gamma(\delta_h\theta)(x)}{2(\eta^2+|h|^2)^\beta}+\frac{\vartheta^2}{8\mathfrak{c}_0|h|^\gamma}
\end{eqnarray*}
where $\mathfrak{c}_0$ is the constant appearing in Lemma \ref{nonlinear-lower-bound-difference-of-solution}.
\end{lemma}
\textbf{Proof}.
For $|h|\leq\eta_0$, letting $\rho=4(\eta^2+|h|^2)^{\frac12}\geq4|h|$
%\begin{eqnarray*}
%\rho=4(\eta^2+|h|^2)^{\frac12}\geq4|h|
%\end{eqnarray*}
in Lemma \ref{bound-difference-velocity} and utilizing Young's inequality lead us to
%Lemma \ref{bound-difference-velocity} and Young's inequality, it follows that
\begin{eqnarray*}
\frac{2\beta|h||\delta_h u|}{\eta^2+|h|^2}\vartheta^2
%&\leq&\frac{2C\beta|h|}{\eta^2+|h|^2}\Big[\rho^{\frac{\gamma}{2}+1-2\alpha}(D_\gamma(\delta_h\theta)(x))^{\frac12}+|h|\|\vartheta\|_{L^\infty_{x,h}}\Big(\frac{\eta^\beta}{\rho^{2\alpha}}+\frac{1}{\rho^{2\alpha-\beta}}\Big)\Big]\vartheta^2\\
%&=&\Big[\frac{(D_\gamma(\delta_h\theta)(x))^{\frac12}}{(\eta^2+|h|^2)^{\frac{\beta}{2}}} \frac{2C\beta|h|\rho^{\frac{\gamma}{2}+1-2\alpha}}{(\eta^2+|h|^2)^{1-\frac{\beta}{2}}}\vartheta^2\Big]\\
%&\leq&\frac{D_\gamma(\delta_h\theta)(x)}{2(\eta^2+|h|^2)^{\beta}}+\frac{C\beta^2|h|^2\rho^{\gamma+2-4\alpha}}{(\eta^2+|h|^2)^{2-\beta}}\vartheta^4
%+\frac{2C\beta|h|^2}{\eta^2+|h|^2}\Big[\|\vartheta\|_{L^\infty_{x,h}}\Big(\frac{\eta^\beta}{\rho^{2\alpha}}+\frac{1}{\rho^{2\alpha-\beta}}\Big)\Big]\vartheta^2\\
\leq\frac{D_\gamma(\delta_h\theta)(x)}{2(\eta^2+|h|^2)^{\beta}}
+\frac{C\beta|h|^2}{\eta^2+|h|^2}\Big[\frac{\beta\rho^{\gamma+2-4\alpha}\vartheta^2}{(\eta^2+|h|^2)^{1-\beta}}+\|\vartheta\|_{L^\infty_{x,h}}\Big(\frac{\eta^\beta}{\rho^{2\alpha}}+\frac{1}{\rho^{2\alpha-\beta}}\Big)\Big]\vartheta^2,
%&\leq&\frac{D_\gamma(\delta_h\theta)(x)}{2(\eta^2+|h|^2)^{\beta}}
%+C\beta\Big[\frac{\beta\rho^{\gamma+2-4\alpha}\vartheta^2}{(\eta^2+|h|^2)^{1-\beta}}+\|\vartheta\|_{L^\infty_{x,h}}\Big(\frac{\eta^\beta}{\rho^{2\alpha}}+\frac{1}{\rho^{2\alpha-\beta}}\Big)\Big]\vartheta^2,
\end{eqnarray*}
which along with the assumption of $\|\vartheta\|_{L^\infty_{x,h}}$, the monotonicity of $\eta$ and $2\alpha-\gamma<\beta$ implies
\begin{eqnarray*}
&&\frac{2\beta|h||\delta_h u|}{\eta^2+|h|^2}\vartheta^2-\frac{D_\gamma(\delta_h\theta)(x)}{2(\eta^2+|h|^2)^{\beta}}
%&\leq&\frac{2C\beta|h|}{\eta^2+|h|^2}\Big[\rho^{\frac{\gamma}{2}+1-2\alpha}(D_\gamma(\delta_h\theta)(x))^{\frac12}+|h|\|\vartheta\|_{L^\infty_{x,h}}\Big(\frac{\eta^\beta}{\rho^{2\alpha}}+\frac{1}{\rho^{2\alpha-\beta}}\Big)\Big]\vartheta^2\\
%&=&\Big[\frac{(D_\gamma(\delta_h\theta)(x))^{\frac12}}{(\eta^2+|h|^2)^{\frac{\beta}{2}}} \frac{2C\beta|h|\rho^{\frac{\gamma}{2}+1-2\alpha}}{(\eta^2+|h|^2)^{1-\frac{\beta}{2}}}\vartheta^2\Big]\\
%&\leq&\frac{D_\gamma(\delta_h\theta)(x)}{2(\eta^2+|h|^2)^{\beta}}+\frac{C\beta^2|h|^2\rho^{\gamma+2-4\alpha}}{(\eta^2+|h|^2)^{2-\beta}}\vartheta^4
%+\frac{2C\beta|h|^2}{\eta^2+|h|^2}\Big[\|\vartheta\|_{L^\infty_{x,h}}\Big(\frac{\eta^\beta}{\rho^{2\alpha}}+\frac{1}{\rho^{2\alpha-\beta}}\Big)\Big]\vartheta^2\\
%&\leq&
%\frac{C\beta|h|^2}{\eta^2+|h|^2}\Big[\frac{\beta\rho^{\gamma+2-4\alpha}\vartheta^2}{(\eta^2+|h|^2)^{1-\beta}}+\|\vartheta\|_{L^\infty_{x,h}}\Big(\frac{\eta^\beta}{\rho^{2\alpha}}+\frac{1}{\rho^{2\alpha-\beta}}\Big)\Big]\vartheta^2\\
%&\leq&C\beta\Big[\frac{\beta\rho^{\gamma+2-4\alpha}\vartheta^2}{(\eta^2+|h|^2)^{1-\beta}}+\|\vartheta\|_{L^\infty_{x,h}}\Big(\frac{\eta^\beta}{\rho^{2\alpha}}+\frac{1}{\rho^{2\alpha-\beta}}\Big)\Big]\vartheta^2\\
\leq C\beta\Big[\frac{\beta\rho^{\gamma+2-4\alpha}\|\vartheta\|^2_{L^\infty_{x,h}}}{(\eta^2+|h|^2)^{1-\beta}}
+\frac{\|\vartheta\|_{L^\infty_{x,h}}}{\rho^{2\alpha-\beta}}\Big]\vartheta^2\\
&&\quad\quad
\leq
C\beta\|\vartheta\|_{L^\infty}\Big[\beta\|\vartheta\|_{L^\infty}\frac{(\eta^2+|h|^2)^{\gamma+\beta-2\alpha}}{(\eta^2+|h|^2)^{\frac\gamma2}}+\frac{(\eta^2+|h|^2)^{\frac{\beta+\gamma}{2}-\alpha}}{(\eta^2+|h|^2)^{\frac\gamma2}}\Big]
\vartheta^2\\
&&\quad\quad
\leq
C\frac{4\beta\|\theta_0\|_{L^\infty}}{\eta^\beta_0}
\Big[\frac{4\beta\|\theta_0\|_{L^\infty}(2\eta^2_0)^{\gamma+\beta-2\alpha}}{\eta^\beta_0|h|^\gamma}
+\frac{(2\eta^2_0)^{\frac{\gamma+\beta}{2}-\alpha}}{|h|^\gamma}\Big]\vartheta^2\\
&&\quad\quad
\leq
C_{\alpha,\beta,\gamma}\Big[\frac{\beta^2\|\theta_0\|^2_{L^\infty}}{\eta^{2(2\alpha-\gamma)}_0}
+\frac{\beta\|\theta_0\|_{L^\infty}}{\eta^{2\alpha-\gamma}_0}\Big]\frac{\vartheta^2}{|h|^\gamma}
\leq\frac{\vartheta^2}{8\mathfrak{c}_0|h|^\gamma},
\end{eqnarray*}
provided we choose $\eta_0$ such that
\begin{eqnarray*}
\eta_0=(16\bar{C}\mathfrak{c}_0\beta\|\theta_0\|_{L^\infty})^{\frac{1}{2\alpha-\gamma}},
\end{eqnarray*}
for some constant $\bar{C}=\max\{C,\frac{1}{\mathfrak{c}_0}\}$. We then complete the proof of Lemma \ref{condition-on-the-initial-data}.

With the help of Lemmas \ref{nonlinear-lower-bound-difference-of-solution}, \ref{differential-equation} and \ref{condition-on-the-initial-data}, we are ready to conclude the proof of Proposition \ref{eventual-regularity}.

\textbf{Proof of Proposition \ref{eventual-regularity}}.
Given $\eta_0$ as in Lemma \ref{condition-on-the-initial-data} and define
\begin{eqnarray*}
t_\star=\sup\{t\in[0,T):\|\vartheta(\tau)\|_{L^\infty_{x,h}}<M,{\rm ~for~ all}~ \tau\in[0,t]\},
\end{eqnarray*}
where $M=:\frac{4\|\theta_0\|_{L^\infty}}{\eta^\beta_0}$.
By \eqref{5.2}, we can see that
\begin{eqnarray*}
\|\vartheta(0)\|_{L^\infty_{x,h}}\leq\frac{2\|\theta_0\|_{L^\infty}}{\eta^\beta_0}=\frac{M}{2},
\end{eqnarray*}
which follows from the continuity  of $\|\vartheta\|_{L^\infty_{x,h}}$ in time that $t_\star>0$ is well-defined.

We claim that $t_\star=T$. Assume, contrary to our claim, that $t_\star<T$.
By the continuity of $\vartheta$ in time and the definition of $t_\star$,
we have $\|\vartheta(t_\star)\|_{L^\infty_{x,h}}=M$. By Lemmas \ref{nonlinear-lower-bound-difference-of-solution}, \ref{differential-equation} and \ref{condition-on-the-initial-data}, \eqref{5.3} becomes
\begin{eqnarray}\label{5.8}
\mathcal{L_{\alpha,\gamma}}\vartheta^2+\frac{1}{4\mathfrak{c}_0|h|^\gamma}\Big[\Big[\frac{\vartheta(x,t;h)}{\|\vartheta\|_{L^\infty}}\Big]^{\frac{\gamma}{1-\beta}}-1\Big]\vartheta^2
+\frac{1}{4\mathfrak{c}_0|h|^\gamma}\Big[\frac{\vartheta(x,t;h)}{\|\vartheta\|_{L^\infty}}\Big]^{\frac{\gamma}{1-\beta}}\vartheta^2\leq0,
\end{eqnarray}
for all $x,h\in\mathbb{T}^n$ with $|h|\leq\eta_0$ and $t\in(0,t_\star]$.
%At this stage, we note that the maximum being attained at $(x_\ast,h_\ast)$ imposes an upper bound for $|h_\ast|$.
Since $\|\vartheta(t_\star)\|_{L^\infty_{x,h}}=M$, we can choose a sufficiently small $\varepsilon>0$ such that
\begin{eqnarray}\label{5.9}
\|\vartheta(t)\|_{L^\infty_{x,h}}\geq\frac{2M}{3}~{\rm~for~all}~t\in[t_\star-\varepsilon,t_\star).
\end{eqnarray}
%Let $t\in[t_\star-\varepsilon,t_\star)$ be arbitrary, where $\varepsilon>0$ is sufficiently small such that $\|\vartheta(t)\|_{L^\infty_{x,h}}\geq\frac{7M}{8}$ for all $t$ in this interval.
%which ensures that the maximum of $|\vartheta(x,t;h)|$ cannot be attained at an $h$ with $|h|>\eta_0$.
For $t\in[t_\star-\varepsilon,t_\star)$, the function $(x,h)\mapsto\vartheta^2(x,t;h)$ has a maximum value which is reached at some point $(x_t,h_t)\in\mathbb{T}^n\times\mathbb{T}^n$. Then
%So, it follows that
\begin{eqnarray*}
\nabla_x\vartheta^2(x_t,t,h_t)=\nabla_h\vartheta^2(x_t,t,h_t)=0,~\Lambda^\gamma\vartheta^2(x_t,t,h_t)\geq0,
%~{\rm and}~|h_t|\leq\eta_0.
\end{eqnarray*}
which follows from the definition of $\mathcal{L_{\alpha,\gamma}}$ that
\begin{eqnarray}\label{5.10}
\mathcal{L_{\alpha,\gamma}}\vartheta^2(x_t,t,h_t)\geq(\partial_t\vartheta^2)(x_t,t,h_t).
\end{eqnarray}
We also have $|h_t|\leq\eta_0$.
Indeed, for every $|h|>\eta_0$, it holds that, for any $(x,t)\in\mathbb{T}^n\times[0,T)$,
\begin{eqnarray*}
|\vartheta(x,t;h)|\leq\frac{2\|\theta(t)\|_{L^\infty}}{|h|^\beta}\leq\frac{2\|\theta_0\|_{L^\infty}}{\eta_0^\beta}=\frac{M}{2},
\end{eqnarray*}
which along with \eqref{5.9} implies $|h_t|\leq\eta_0$. Reading \eqref{5.8} at $(x_t,t,h_t)$ and using $\frac{|\vartheta(x_t,t,h_t)|}{\|\vartheta(t)\|_{L^\infty_{x,h}}}=1$ lead us to
\begin{eqnarray*}
%(\partial_t\vartheta^2)(x_t,t,h_t)+\frac{49M^2}{256c_0\eta^\gamma_0}
%\leq\mathcal{L}\vartheta^2(x_t,t,h_t)+\frac{\vartheta^2(x_t,t,h_t)}{4c_0|h_t|^\gamma}.
\mathcal{L}\vartheta^2(x_t,t,h_t)+\frac{\vartheta^2(x_t,t,h_t)}{4\mathfrak{c}_0|h_t|^\gamma}\leq0,
\end{eqnarray*}
which together with \eqref{5.10} gives
\begin{eqnarray*}
%(\partial_t\vartheta^2)(x_t,t,h_t)+\frac{49M^2}{256c_0\eta^\gamma_0}
%\leq\mathcal{L}\vartheta^2(x_t,t,h_t)+\frac{\vartheta^2(x_t,t,h_t)}{4c_0|h_t|^\gamma}.
(\partial_t\vartheta^2)(x_t,t,h_t)\leq\mathcal{L}\vartheta^2(x_t,t,h_t)+\frac{\vartheta^2(x_t,t,h_t)}{4\mathfrak{c}_0|h_t|^\gamma}
\leq0.
\end{eqnarray*}
%For $t\in[t_\star-\varepsilon,t_\star)$, the function $(x,h)\mapsto\vartheta^2(x,t;h)$ has a maximum value which is reached at some point $(x_t,h_t)\in\mathbb{T}^n\times\mathbb{T}^n$. Then
%So, it follows that
%\begin{eqnarray*}
%\nabla_x\vartheta^2(x_t,t,h_t)=\nabla_h\vartheta^2(x_t,t,h_t)=0,~\Lambda^\gamma\vartheta^2(x_t,t,h_t)\geq0.
%~{\rm and}~|h_t|\leq\eta_0.
%\end{eqnarray*}
%$\nabla_x\vartheta^2(x_t,t,h_t)=\nabla_h\vartheta^2(x_t,t,h_t)=0$, $\Lambda^\gamma\vartheta^2(x_t,t,h_t)\geq0$ and $|h_t|\leq\eta_0$.
%Therefore,
%\begin{eqnarray}\label{5.10}
%\mathcal{L_{\alpha,\gamma}}\vartheta^2(x_t,t,h_t)\geq(\partial_t\vartheta^2)(x_t,t,h_t).
%\end{eqnarray}
%By \eqref{5.9} and $|h_t|\leq\eta_0$, we arrive at
%\begin{eqnarray*}
%\frac{\vartheta^2(x_t,t,h_t)}{4c_0|h_t|^\gamma}\geq\frac{49M^2}{256c_0\eta^\gamma_0},
%\end{eqnarray*}
%which together with \eqref{5.10} yields that
%\begin{eqnarray*}
%(\partial_t\vartheta^2)(x_t,t,h_t)+\frac{49M^2}{256c_0\eta^\gamma_0}
%\leq\mathcal{L}\vartheta^2(x_t,t,h_t)+\frac{\vartheta^2(x_t,t,h_t)}{4c_0|h_t|^\gamma}.
%\mathcal{L}\vartheta^2(x_t,t,h_t)+\frac{\vartheta^2(x_t,t,h_t)}{4c_0|h_t|^\gamma}\geq(\partial_t\vartheta^2)(x_t,t,h_t)+\frac{49M^2}{256c_0\eta^\gamma_0}.
%\end{eqnarray*}
%Utilizing this inequality, reading \eqref{5.8} at $(x_t,t,h_t)$ and using $\frac{|\vartheta(x_t,t,h_t)|}{\|\vartheta(t)\|_{L^\infty_{x,h}}}=1$ lead us to
%\begin{eqnarray*}
%(\partial_t\vartheta^2)(x_t,t,h_t)
%\leq-\frac{49M^2}{256c_0\eta^\gamma_0},{\rm ~for~ all}~ t\in[t_\star-\varepsilon,t_\star).
%\end{eqnarray*}
Following an same argument as the Appendix B in \cite{[Constantin-Tarfulea-Vicol]}, we can show that
\begin{eqnarray*}
\frac{d}{dt}\|\vartheta(t)\|^2_{L^\infty_{x,h}}\leq(\partial_t\vartheta^2)(x_t,t,h_t)\leq0,{\rm ~for~ all}~ t\in[t_\star-\varepsilon,t_\star).
\end{eqnarray*}
%Moreover, by Lemma \ref{} with $f(t,\lambda)=\vartheta^2(x,t,h)$ and $\lambda=(x,h)\in_{}\mathcal{K}=\mathbb{T}^n\times\mathbb{T}^n$,
%it follows that
%\begin{eqnarray*}
%\frac{d}{dt}\|\vartheta(t)\|^2_{L^\infty_{x,h}}\leq(\partial_t\vartheta^2)(x_t,t,h_t)\leq-\frac{49M^2}{256c_0\eta^\gamma_0},{\rm ~for~ all}~ t\in[t_\star-\varepsilon,t_\star),
%\end{eqnarray*}
%By direct integration in time, we obtain that, for all $t\in[t_\star-\varepsilon,t_\star)$,
%\begin{eqnarray*}
%\|\vartheta(t)\|_{L^\infty_{x,h}}
%\leq
%\|\vartheta(t_\star-\varepsilon)\|_{L^\infty_{x,h}}-\frac{49M^2}{256c_0\eta^\gamma_0}(t-t_\star+\varepsilon)
%\end{eqnarray*}
Then, for all $t\in[t_\star-\varepsilon,t_\star)$, $\|\vartheta(t)\|_{L^\infty_{x,h}}\leq\|\vartheta(t_\star-\varepsilon)\|_{L^\infty_{x,h}}.$
%\begin{eqnarray*}
%\|\vartheta(t)\|_{L^\infty_{x,h}}\leq\|\vartheta(t_\star-\varepsilon)\|_{L^\infty_{x,h}}.
%\end{eqnarray*}
Letting $t\rightarrow t_\star$ and utilizing $\|\vartheta(t_\star-\varepsilon)\|_{L^\infty_{x,h}}<M$ thus implies that
$\|\vartheta(t_\star)\|_{L^\infty_{x,h}}<M$,
%\begin{eqnarray*}
%\|\vartheta(t_\star)\|_{L^\infty_{x,h}}<M,
%\end{eqnarray*}
which is a contradiction with the definition of $t_\star$.  Then $t_\star=T$ and $\vartheta\in L^\infty(\mathbb{T}^n\times(0,T))$.
Note that the solution of \eqref{5.7} is
%Given $\eta_0$ as in Lemma \ref{condition-on-the-initial-data}, we know that the solution of \eqref{5.7} is
\begin{eqnarray*}
\eta(t)
=
\begin{cases}
(\eta^\gamma_0-\frac{\gamma t}{16\mathfrak{c}_0\beta})^{\frac1\gamma},
&  \mbox{if $0\leq t\leq T^\ast_\alpha,$ }\\
0,
& \mbox{if $T^\ast_\alpha<t<T,$ } \\
\end{cases}
\end{eqnarray*}
where
\begin{eqnarray*}
T^\ast_\alpha=\frac{16\mathfrak{c}_0\beta\eta^\gamma_0}{\gamma}
=\frac{16\mathfrak{c}_0\mathfrak{c}^{\frac{\gamma}{2\alpha-\gamma}}}{\gamma }\beta^{\frac{2\alpha}{2\alpha-\gamma}}\|\theta_0\|^{\frac{\gamma}{2\alpha-\gamma}}_{L^\infty}.
\end{eqnarray*}
It follows that
%Noting that $\eta(t)\equiv0$ for $t\in(T^\ast_\alpha,T)$, we obtain that
\begin{eqnarray*}
[\theta(t)]_{C^\beta}=\|\vartheta(t)\|_{L^\infty_{x,h}}\leq M=\frac{4\|\theta_0\|_{L^\infty}}{\eta^\beta_0},{\rm ~for~ all}~ t\in(T^\ast_\alpha,T),
\end{eqnarray*}
which finishes the proof of Proposition \ref{eventual-regularity}.

%We now turn to prove Theorem \ref{the-2}.
\textbf{Proof of Theorem \ref{the-2}}.
Given nonnegative $\theta_0\in H^s(\mathbb{T}^n)$ with $s>\frac{n}{2}+1$.
In view of Proposition \ref{local}, there exists a unique local solution to \eqref{1.1}
in the class \eqref{solution-class}.
By Corollary \ref{lower-bound-time}, the local existence $T$ satisfies
\begin{eqnarray}\label{5.11}
T_1:=\frac{1}{C_1\|\theta_0\|^{1-\frac{n+4\alpha}{2s}}_{L^2(\mathbb{T}^n)}\|\Lambda^s\theta_0\|^{\frac{n+4\alpha}{2s}}_{L^2(\mathbb{T}^n)}}\leq T,
\end{eqnarray}
where $C_1$ is a positive constant independent of $\gamma$, $\beta$ and $\theta_0$.

Next we will choose $\gamma$ sufficiently close to $2\alpha$ such that $T^\ast_\alpha\leq T_1$, where $T^\ast_\alpha$ is as in \eqref{5.1}. By Gagliardo-Nirenberg inequality, $T^\ast_\alpha$ can be bounded as
\begin{eqnarray}\label{5.12}
\begin{split}
T^\ast_\alpha
&\leq C\beta^{\frac{2\alpha}{2\alpha-\gamma}}\Big[C_2\|\theta_0\|^{\frac{n}{2s}}_{\dot{H}^s(\mathbb{T}^n)}\|\theta_0\|^{1-\frac{n}{2s}}_{L^2(\mathbb{T}^n)}\Big]^{\frac{\gamma}{2\alpha-\gamma}}\\
&\leq C_3\beta^{\frac{2\alpha}{2\alpha-\gamma}}\|\theta_0\|^{\frac{n\gamma}{2s(2\alpha-\gamma)}}_{\dot{H}^s(\mathbb{T}^n)}
\|\theta_0\|^{\frac{\gamma(2s-n)}{2s(2\alpha-\gamma)}}_{L^2(\mathbb{T}^n)},
\end{split}
\end{eqnarray}
where $C_2>0$ is independent of $\gamma$, $\beta$ and $\theta_0$, $C_3=CC^{\frac{\gamma}{2\alpha-\gamma}}_2$ and $C$ is a constant that can be bounded from below and above for all $\gamma\in[\gamma_0,1)$.
%By Gagliardo-Nirenberg inequality, there exists a constant $C_2>0$ independent of $\gamma$, $\beta$ and $\theta_0$ such that
%\begin{eqnarray*}
%\|\theta_0\|_{L^\infty(\mathbb{T}^n)}\leq C_2\|\theta_0\|^{\frac{n}{2s}}_{\dot{H}^s(\mathbb{T}^n)}\|\theta_0\|^{1-\frac{n}{2s}}_{L^2(\mathbb{T}^n)}.
%\end{eqnarray*}
%Therefore,
%\begin{eqnarray*}
%T^\ast_\alpha
%&\leq&C\beta^{\frac{2\alpha}{2\alpha-\gamma}}\Big(C_2\|\theta_0\|^{\frac{n}{2s}}_{\dot{H}^s(\mathbb{T}^n)}\|\theta_0\|^{1-\frac{n}{2s}}_{L^2(\mathbb{T}^n)}\Big)^{\frac{\gamma}{2\alpha-\gamma}}\\
%&\leq&C_3\beta^{\frac{2\alpha}{2\alpha-\gamma}}\|\theta_0\|^{\frac{n\gamma}{2s(2\alpha-\gamma)}}_{\dot{H}^s(\mathbb{T}^n)}
%\|\theta_0\|^{\frac{\gamma(2s-n)}{2s(2\alpha-\gamma)}}_{L^2(\mathbb{T}^n)},
%\end{eqnarray*}
%where $C_3=CC^{\frac{\gamma}{2\alpha-\gamma}}_2$ and $C$ is a constant that can be bounded from below and above for all $\gamma\in[\gamma_0,1)$.
%In order to get $T^\ast_\alpha\leq T_1$, we observe by \eqref{} and \eqref{} that it is sufficient to show
%\begin{eqnarray*}
%C_1C_3\beta^{\frac{2\alpha}{2\alpha-\gamma}}
%\Big(\|\theta_0\|^{\frac{n+4\alpha-2\gamma}{2s}}_{\dot{H}^s(\R^n)}
%\|\theta_0\|^{1-\frac{n+4\alpha-2\gamma}{2s}}_{L^2(\R^n)}\Big)^{\frac{2\alpha}{2\alpha-\gamma}}
%\leq1.
%\end{eqnarray*}
%In order to get $T^\ast_\alpha\leq T_1$, we observe
By \eqref{5.11}, \eqref{5.12} and \eqref{assumption-initial-data}, we have
\begin{eqnarray*}
\frac{T^\ast_\alpha}{T_1}
\leq C_1C_3\beta^{\frac{2\alpha}{2\alpha-\gamma}}
\Big[\|\theta_0\|^{\frac{n+4\alpha-2\gamma}{2s}}_{\dot{H}^s(\mathbb{T}^n)}
\|\theta_0\|^{1-\frac{n+4\alpha-2\gamma}{2s}}_{L^2(\mathbb{T}^n)}\Big]^{\frac{2\alpha}{2\alpha-\gamma}}
\leq C_1C_3\beta^{\frac{2\alpha}{2\alpha-\gamma}}
R^{\frac{2\alpha}{2\alpha-\gamma}}.
\end{eqnarray*}
In order to obtain $T^\ast_\alpha\leq T_1$, we need to show $C_1C_3\beta^{\frac{2\alpha}{2\alpha-\gamma}}
R^{\frac{2\alpha}{2\alpha-\gamma}}
\leq1$,
%\begin{eqnarray*}
%C_1C_3\beta^{\frac{2\alpha}{2\alpha-\gamma}}
%R^{\frac{2\alpha}{2\alpha-\gamma}}
%\leq1,
%\end{eqnarray*}
which is equivalent to
\begin{eqnarray}\label{5.13}
\beta
\leq R^{-1}\Big(\frac{1}{C_1C_3}\Big)^{\frac{2\alpha-\gamma}{2\alpha}}.
\end{eqnarray}
To do this, we choose $\beta=\min\{1-2\gamma+2\alpha,\alpha\}$. By \eqref{5.13} it follows that there exists $\gamma_1=\gamma_1(R)\in[\gamma_0,2\alpha)$
such that $T^\ast_\alpha\leq T_1$ for all $\gamma\in[\gamma_1,2\alpha)$.
To finish, let $T_{max}$ be the maximal existence time for the solution to \eqref{1.1} in the class. For the moment, suppose that $T_{max}<\infty$, then
$\theta\in C([0,T_{max});H^s(\mathbb{T}^n))$ with $T^\ast_\alpha<T_1\leq T_{max}$.
By Proposition \ref{eventual-regularity}, $\theta\in C^\infty(\mathbb{T}^n\times(T^\ast_\alpha,T_{max}])$ and $\theta(T_{max})\in H^s(\mathbb{T}^n)$. Thus, by Proposition \ref{local}, we can extend $\theta$ in the class
\eqref{solution-class} beyond $T_{max}$. This leads to a contradiction with the maximality of $T_{max}$. Consequently, $T_{max}=\infty$ and $\theta$ is a global $H^s$-solution.
\section{Finite-time blow-up for the supercritical case $\gamma<\alpha$}
In this section, we will prove the finite-time blow-up of radial smooth solutions to \eqref{1.1} with $\gamma\in(0,\alpha)$.
We start by recalling a weighted inequality involving the nonlinear term in \eqref{1.1}, whose proof can be found in \cite{[Zhang]}.
%The one-dimensional version of analogous weighted inequality was first established by Li and Rodrigo in \cite{[Li-Rodrigo20]}.
\begin{lemma}\label{nonlinear-inequality-exponential-weight}
Let $n\geq2$ and $\alpha\in(0,1)$.
Let $f:\R^n\rightarrow\R$ be a radially symmetric Schwartz function. Then
\begin{eqnarray*}
\int_{\R^n}\frac{\Lambda^{-2+2\alpha}\nabla f(x)\cdot\nabla f(x)}{|x|^{n}}e^{-|x|}dx
\geq C'_{n,\alpha}\int_{\R^n}\frac{(f(0)-f(x))^2}{|x|^{n+2\alpha}}dx-C''_{n,\alpha}\|f\|^2_{L^\infty},
\end{eqnarray*}
where $C'_{n,\alpha}$ and $C''_{n,\alpha}$ are constants depending only on $n$ and $\gamma$.
\end{lemma}
%\textbf{Proof.} We refer to \cite{[Zhang]} for the proof of this lemma. The one-dimensional version of this weighted inequality was first established by Li and Rodrigo in \cite{[Li-Rodrigo20]}.
In order to handle the dissipative term, we establish a type of weighted inequality below, whose one-dimensional version was proved by Li and Rodrigo in \cite{[Li-Rodrigo20]}.
\begin{lemma}\label{dissipation}
Let $n\geq2$ and $\gamma\in(0,1)$.
Let $f:\R^n\rightarrow\R$ be a Schwartz function.
%such that
%\begin{eqnarray}\label{symmetry}
%f(-x)=f(x),~{\rm for ~all}~x\in\R^n.
%\end{eqnarray}
Then
\begin{eqnarray*}
\Big|\int_{\R^n}
\frac{\Lambda^\gamma f(0)-\Lambda^\gamma f(x)}{|x|^{n}}e^{-|x|}dx\Big|
\leq C_{n,\gamma}\int_{\R^n}\frac{|f(0)-f(x)|}{|x|^{n+\gamma}}\ln\Big(e+\frac{1}{|x|}\Big)dx,
\end{eqnarray*}
where $C_{n,\gamma}$ is a constant depending only on $n$ and $\gamma$.
\end{lemma}
\textbf{Proof.}
By \eqref{Fractional-Laplacian}, a change of variables and the use of parity, we rewrite $\Lambda^\gamma f(x)$ as
\begin{eqnarray*}
\Lambda^\gamma f(x)
=\frac12C_{n,\gamma}P.V.\int_{\mathbb{R}^n}\frac{2f(x)-f(x+y)-f(x-y)}{|y|^{n+\gamma}}dy,
\end{eqnarray*}
which along with
\begin{eqnarray*}
(\Lambda^\gamma f)(0)=C_{n,\gamma}P.V.\int_{\R^n}\frac{f(0)-f(y)}{|y|^{n+\gamma}}dy,
\end{eqnarray*}
implies that
\begin{eqnarray*}
\begin{split}
&\frac{2}{C_{n,\gamma}}\int_{\R^n}
\frac{\Lambda^\gamma f(0)-\Lambda^\gamma f(x)}{|x|^{n}}e^{-|x|}dx\\
&=\iint_{\R^n\times\R^n}\frac{2(f(0)-f(x))+f(x+y)+f(x-y)-2f(y)}{|x|^n|y|^{n+\gamma}}e^{-|x|}dxdy\\
&=\iint_{\R^n\times\R^n}\underbrace{\frac{-2g(x)+g(x+y)+g(x-y)-2g(y)}{|x|^n|y|^{n+\gamma}}e^{-|x|}}_{H(x,y)}dxdy\\
&=\Big(\underbrace{\iint_{\frac{|x|}{e}\leq|y|\leq e|x|}}_{X_1}
+\underbrace{\iint_{|y|\geq e|x|}}_{X_2}
+\underbrace{\iint_{|y|\leq\frac{|x|}{e}}}_{X_3}\Big)H(x,y)dxdy,
\end{split}
\end{eqnarray*}
where $g(x)=f(x)-f(0)$.

We estimate $X_1$, $X_2$ and $X_3$ separately. First, the term $X_1$ can be estimated as
\begin{eqnarray}\label{K1}
\begin{split}
|X_1|
&\leq\iint\limits_{\frac{|x|}{e}\leq|y|\leq e|x|}\frac{2|g(x)|}{|x|^{n}|y|^{n+\gamma}}dxdy
+e^{n+\gamma}\iint\limits_{\frac{|x|}{e}\leq|y|\leq e|x|}\frac{|g(x+y)|}{|x|^{2n+\gamma}}dxdy\\
&\
+e^{n+\gamma}\iint\limits_{\frac{|x|}{e}\leq|y|\leq e|x|}\frac{|g(x-y)|}{|x|^{2n+\gamma}}dxdy
+\iint\limits_{\frac{|x|}{e}\leq|y|\leq e|x|}\frac{2|g(y)|}{|x|^{n}|y|^{n+\gamma}}dxdy\\
&=\int_{\R^n}\frac{2|g(x)|}{|x|^{n}}\Big[\int_{\frac{|x|}{e}\leq|y|\leq e|x|}\frac{dy}{|y|^{n+\gamma}}\Big]dx
+2e^{n+\gamma}\iint\limits_{\frac{|x|}{e}\leq|x-z|\leq e|x|}\frac{|g(z)|}{|x|^{2n+\gamma}}dxdz\\
&\ \ \ \ \ \ \
+\int_{\R^n}\frac{2|g(y)|}{|y|^{n+\gamma}}\Big[\int_{\frac{|y|}{e}\leq|x|\leq e|y|}\frac{dx}{|x|^{n}}\Big]dy\\
&\leq C'_{n,\gamma}\int_{\R^n}\frac{|g(x)|}{|x|^{n+\gamma}}dx
+2e^{n+\gamma}\int_{\R^n}|g(z)|\Big[\int\limits_{|x|\geq\frac{|z|}{e+1}}
\frac{dx}{|x|^{2n+\gamma}}\Big]dz\\
&\leq C''_{n,\gamma}\int_{\R^n}\frac{|g(x)|}{|x|^{n+\gamma}}dx,
\end{split}
\end{eqnarray}
where in above identity we used a change of variables.
Secondly, by a change of variables, we rewrite $X_2$  as
\begin{eqnarray*}
X_{2}
=\underbrace{\iint\limits_{|y|\geq e|x|}\frac{-2g(x)e^{-|x|}}{|x|^n|y|^{n+\gamma}}dxdy}_{X_{21}}
+\underbrace{\iint\limits_{|y-x|\geq e|x|}\frac{2g(y)e^{-|x|}dxdy}{|x|^n|y-x|^{n+\gamma}}
-\iint\limits_{|y|\geq e|x|}\frac{2g(y)e^{-|x|}dxdy}{|x|^n|y|^{n+\gamma}}}_{X_{22}}.
\end{eqnarray*}
The estimate of $X_{21}$ is straightforward:
\begin{eqnarray*}
|X_{21}|
%&\leq&\iint\limits_{|y|\geq e|x|}\frac{2|g(x)|}{|x|^n|y|^{n+\gamma}}dxdy\\
\leq\int_{\R^n}\frac{2|g(x)|}{|x|^{n}}\Big[\int\limits_{|y|\geq e|x|}\frac{dy}{|y|^{n+\gamma}}\Big]dx
=C'_{n,\gamma}\int_{\R^n}\frac{|g(x)|}{|x|^{n+\gamma}}dx.
\end{eqnarray*}
By an elementary formula
\begin{eqnarray}\label{expression}
\int_{\mathbf{A}}h_1dz-\int_{\mathbf{B}}h_2dz
=\int_{\mathbf{A}\cap\mathbf{B}}(h_1-h_2)dz
+\int_{\mathbf{A}\cap\mathbf{B}^c}h_1dz
-\int_{\mathbf{B}\cap\mathbf{A}^c}h_2dz
\end{eqnarray}
we express the term $X_{22}$ as
\begin{eqnarray*}
\begin{split}
X_{22}
&=\underbrace{\iint_{\mathbf{A_1}\cap\mathbf{B_1}}\frac{2g(y)e^{-|x|}}{|x|^n}\Big[\frac{1}{|y-x|^{n+\gamma}}
-\frac{1}{|y|^{n+\gamma}}\Big]dxdy}_{X_{221}}\\
&\ \ \
+\underbrace{\iint_{\mathbf{A}_1\cap\mathbf{B}^c_1}\frac{2g(y)e^{-|x|}dxdy}{|x|^n|y-x|^{n+\gamma}}}_{X_{222}}
-\underbrace{\iint_{\mathbf{B}_1\cap\mathbf{A}^c_1}\frac{2g(y)e^{-|x|}dxdy}{|x|^n|y|^{n+\gamma}}}_{X_{223}}.
\end{split}
\end{eqnarray*}
where
\begin{eqnarray*}
\mathbf{A}_1=\{(x,y)\in\mathbb{R}^{2n}:|y-x|\geq e|x|\},~\mathbf{B}_1=\{(x,y)\in\mathbb{R}^{2n}:|y|\geq e|x|\}.
\end{eqnarray*}
By the mean-value theorem, we estimate the term $X_{221}$ as
\begin{eqnarray*}
\begin{split}
|X_{221}|
&\leq
\iint_{\mathbf{A}_1\cap\mathbf{B}_1}\frac{2|g(y)|}{|x|^n}\Big|\frac{1}{|y-x|^{n+\gamma}}
-\frac{1}{|y|^{n+\gamma}}\Big|dxdy\\
&\leq(n+\gamma)\Big(\frac{e}{e-1}\Big)^{n+\gamma+1}\iint\limits_{|y|\geq e|x|}\frac{|g(y)|dxdy}{|x|^{n-1}|y|^{n+\gamma+1}}\\
&=C_{n,\gamma}\int_{\R^n}\frac{|g(y)|}{|y|^{n+\gamma+1}}\Big[\int_{|x|\leq\frac{|y|}{e}}\frac{dx}{|x|^{n-1}}\Big]dy\\
&\leq C'_{n,\gamma}\int_{\R^n}\frac{|g(y)|}{|y|^{n+\gamma}}dy.
\end{split}
\end{eqnarray*}
The estimates of $X_{222}$ and $X_{223}$ are standard:
\begin{eqnarray*}
\begin{split}
|X_{222}|+|X_{223}|
&\leq\iint\limits_{(e-1)|x|\leq|y|\leq e|x|}\frac{2|g(y)|dxdy}{|x|^n|y|^{n+\gamma}}
+\iint\limits_{e|x|\leq|y|\leq(e+1)|x|}\frac{2|g(y)|dxdy}{|x|^n|y|^{n+\gamma}}\\
&=\int_{\R^n}\frac{2|g(y)|}{|y|^{n+\gamma}}\Big[\int_{\frac{|y|}{e+1}\leq|x|\leq\frac{|y|}{e-1}}\frac{dx}{|x|^{n}}\Big]dy
\leq C_n\int_{\R^n}\frac{|g(y)|}{|y|^{n+\gamma}}dy.
\end{split}
\end{eqnarray*}
Thus, we obtain that
\begin{eqnarray}\label{K2}
|X_{2}|
\leq C'_{n,\gamma}\int_{\R^n}\frac{|g(x)|}{|x|^{n+\gamma}}dx.
\end{eqnarray}
Similarly, by a change of variables, $X_{3}$ can be bounded by
\begin{eqnarray*}
\begin{split}
|X_{3}|
&\leq
\underbrace{\Big|\iint\limits_{|y|\leq\frac{|x-y|}{e}}\frac{g(x)e^{-|x-y|}}{|x-y|^n|y|^{n+\gamma}}dxdy
-\iint\limits_{|y|\leq\frac{|x|}{e}}\frac{g(x)e^{-|x|}}{|x|^n|y|^{n+\gamma}}dxdy\Big|}_{X_{31}}\\
&\
+\underbrace{\Big|\iint\limits_{|y|\leq\frac{|x+y|}{e}}\frac{g(x)e^{-|x+y|}dxdy}{|x+y|^n|y|^{n+\gamma}}-\iint\limits_{|y|\leq\frac{|x|}{e}}\frac{g(x)e^{-|x|}dxdy}{|x|^n|y|^{n+\gamma}}\Big|}_{X_{32}}
+\underbrace{\iint\limits_{|y|\leq\frac{|x|}{e}}\frac{2|g(y)|dxdy}{|x|^n|y|^{n+\gamma}e^{|x|}}}_{X_{33}}.
\end{split}
\end{eqnarray*}
By the identity \eqref{expression} and the mean value theorem, we can derive that
\begin{eqnarray*}
\begin{split}
X_{31}
&\leq
\iint_{\mathbf{A_2}\cap\mathbf{B_2}}\frac{|g(x)|}{|y|^{n+\gamma}}
\Big|\frac{e^{-|x-y|}}{|x-y|^n}-\frac{e^{-|x|}}{|x|^n}\Big|dxdy\\
&\ \ \
+\iint_{\mathbf{A_2}\cap\mathbf{B}^c_2}\frac{|g(x)|e^{-|x-y|}}{|x-y|^n|y|^{n+\gamma}}dxdy
+\iint_{\mathbf{A^c_2}\cap\mathbf{B}_2}\frac{|g(x)|e^{-|x|}}{|x|^n|y|^{n+\gamma}}dxdy\\
&\leq
C_n\iint_{\mathbf{A_2}\cap\mathbf{B_2}}
\frac{|g(x)|}{|y|^{n-1+\gamma}}
e^{-(1-\frac{1}{e})|x|}\Big[\frac{1}{|x|^n}+\frac{1}{|x|^{n+1}}\Big]dxdy\\
&\ \ \
+\iint\limits_{\mathbf{A_2}\cap\mathbf{B}^c_2}\frac{|g(x)|}{|x|^n|y|^{n+\gamma}}dxdy
+\iint\limits_{\mathbf{A^c_2}\cap\mathbf{B}_2}\frac{|g(x)|}{|x|^n|y|^{n+\gamma}}dxdy\\
&\leq
C_n
\iint_{\mathbf{A_2}\cap\mathbf{B_2}}
\frac{|g(x)|dxdy}{|x|^{n+1}|y|^{n-1+\gamma}}+2\iint_{|y|\geq\frac{e-1}{e^2}|x|}\frac{|g(x)|}{|x|^n|y|^{n+\gamma}}dxdy\\
&\leq C_n\int_{\R^n}\frac{|g(x)|}{|x|^{n+1}}\int_{|y|\leq\frac{|x|}{e}}\frac{dy}{|y|^{n-1+\gamma}}dx
+2\int_{\R^n}\frac{|g(x)|}{|x|^{n}}\int_{|y|\geq\frac{e-1}{e^2}|x|}\frac{dy}{|y|^{n+\gamma}}dx\\
&\leq C_{n,\gamma}\int_{\R^n}\frac{|g(x)|}{|x|^{n+\gamma}}dx,
\end{split}
\end{eqnarray*}
where
\begin{eqnarray*}
\mathbf{A}_2=\Big\{(x,y)\in\mathbb{R}^{2n}:|y|\leq\frac{|x-y|}{e}\Big\},~\mathbf{B}_2=\Big\{(x,y)\in\mathbb{R}^{2n}:|y|\leq\frac{|x|}{e}\Big\}.
\end{eqnarray*}
Completely similar to the estimate of $X_{31}$, $X_{32}$ can be bounded by
\begin{eqnarray*}
X_{32}
\leq C_{n,\gamma}\int_{\R^n}\frac{|g(x)|}{|x|^{n+\gamma}}dx.
\end{eqnarray*}
Finally, by direct computation and Lemma \ref{Inequality}, we can obtain that
\begin{eqnarray*}
X_{33}
\lesssim\int_{\R^n}\frac{|g(y)|}{|y|^{n+\gamma}}\Big[\int^\infty_{|y|}\frac{e^{-r}}{r}dr\Big]dy
\leq C_n\int_{\R^n}\frac{|g(y)|}{|y|^{n+\gamma}}\ln\Big(e+\frac{1}{|y|}\Big)dy.
%&\leq&C_{n,\gamma}\int_{\R^n}\frac{|g(y)|}{|y|^{n+\gamma}}\ln\Big(e+\frac{1}{|y|}\Big)dy.
\end{eqnarray*}
Thus, we have
\begin{eqnarray*}
|X_{3}|
\leq C_{n,\gamma}\int_{\R^n}\frac{|g(y)|}{|y|^{n+\gamma}}\ln\Big(e+\frac{1}{|y|}\Big)dy,
\end{eqnarray*}
which along with \eqref{K1} and \eqref{K2} then complete the proof of Lemma \ref{dissipation}.

With the help of Lemmas \ref{nonlinear-inequality-exponential-weight} and \ref{dissipation}, we are ready to prove Theorem \ref{the-3}.

\textbf{Proof of Theorem \ref{the-3}}.
Taking into account the quantity
\begin{eqnarray*}
J(t)=\int_{\R^n}\frac{\theta(0,t)-\theta(x,t)}{|x|^{n}}e^{-|x|}dx
\end{eqnarray*}
and computing $J^\prime(t)$ from $\eqref{1.1}_1$, we arrive at the differential identity
\begin{eqnarray*}
J^\prime(t)=\int_{\R^n} \frac{\Lambda^{-2+2\alpha}\nabla\theta\cdot\nabla\theta}{|x|^{n}}e^{-|x|}dx
-\kappa\int_{\R^n}\frac{\Lambda^\gamma \theta(0,t)-\Lambda^\gamma \theta(x,t)}{|x|^{n}}e^{-|x|}dx,
\end{eqnarray*}
where we used the fact that the velocity vanishes at the origin, that is, by Lemma \ref{radial-property-preserved},
\begin{eqnarray*}
u(0,t)=C_{n,\alpha}\int_{\R^n}\frac{\nabla\theta(y,t)}{|y|^{n-2+2\alpha}}dy=C_{n,\alpha}\int_{\R^n}\frac{y\theta^\prime(|y|,t)}{|y|^{n-1+2\alpha}}dy\equiv0.
\end{eqnarray*}
It follows from Lemmas \ref{nonlinear-inequality-exponential-weight}, \ref{maximum} and \ref{dissipation} that
\begin{eqnarray}\label{J-derivative-6.4}
\begin{split}
J^\prime(t)
&\geq C_1\int_{\R^n} \frac{(\theta(0,t)-\theta(x,t))^2}{|x|^{n+2\alpha}}dx
-C_2\|\theta_0\|^2_{L^\infty}\\
&\ \ \
-\kappa C_3\int_{\R^n}\frac{ \theta(0,t)- \theta(x,t)}{|x|^{n+\gamma}}\ln\Big(e+\frac{1}{|x|}\Big)dx.
\end{split}
\end{eqnarray}
By H\"{o}lder's inequality, we derive that
\begin{eqnarray*}
\begin{split}
|J(t)|
%&\leq&\int_{\R^n}\frac{|\theta(0,t)-\theta(x,t)|}{|x|^{n}}e^{-|x|}dx\\
&\leq\Big[\int_{\R^n}\frac{|\theta(0,t)-\theta(x,t)|^2}{|x|^{n+2\alpha}}dx\Big]^{\frac12}
\Big[\int_{\R^n}\frac{e^{-2|x|}}{|x|^{n-2\alpha}}dx\Big]^{\frac12}\\
&\leq C_4(n,\alpha)\Big[\int_{\R^n}\frac{|\theta(0,t)-\theta(x,t)|^2}{|x|^{n+2\alpha}}dx\Big]^{\frac12},
\end{split}
\end{eqnarray*}
which is equivalent to
\begin{eqnarray}\label{lower-bound-6.5}
\int_{\R^n}\frac{|\theta(0,t)-\theta(x,t)|^2}{|x|^{n+2\alpha}}dx
\geq
\frac{(J(t))^2}{C^2_4}.
\end{eqnarray}
Furthermore, by H\"{o}lder's inequality, Lemma \ref{maximum} and the restriction $\gamma\in(0,\alpha)$, we have
\begin{eqnarray}\label{dissipation-6.6}
&&\Big|\int_{\R^n}\frac{ \theta(0,t)- \theta(x,t)}{|x|^{n+\gamma}}\ln\Big(e+\frac{1}{|x|}\Big)dx\Big|\nonumber\\
&\leq&\Big[\int_{|x|\leq1}\frac{e^{|x|}\ln^2(e+\frac{1}{|x|})}{|x|^{n+2(\gamma-\alpha)}}dx\Big]^{\frac12}
(Q(t))^{\frac12}+2\|\theta_0\|_{L^\infty}\int_{|x|>1}\frac{\ln(e+\frac{1}{|x|})}{|x|^{n+\gamma}}dx\nonumber\\
&\leq&C_5(n,\alpha,\gamma)(Q(t))^{\frac12}
+C_6(n,\gamma)\|\theta_0\|_{L^\infty},
\end{eqnarray}
where
\begin{eqnarray*}
Q(t)=\int_{|x|\leq1}\frac{ |\theta(0,t)- \theta(x,t)|^2}{|x|^{n+2\alpha}}e^{-|x|}dx.
\end{eqnarray*}
Substituting \eqref{dissipation-6.6} and \eqref{lower-bound-6.5} into \eqref{J-derivative-6.4} and utilizing the inequality $a\tau-b\sqrt{\tau}\geq-\frac{b^2}{4a^2}$ for $a,b>0$, $\tau\geq0$, we conclude that
\begin{eqnarray*}
\begin{split}
J^\prime(t)
&\geq\frac{C_1}{2C^2_4}(J(t))^2
-C_2\|\theta_0\|^2_{L^\infty}-\kappa C_3C_6\|\theta_0\|_{L^\infty}
+\frac{C_1}{2}Q(t)-\kappa C_3C_5(Q(t))^{\frac12}\\
&\geq\frac{C_1}{2C^2_4}(J(t))^2
-C_2\|\theta_0\|^2_{L^\infty}-\kappa C_3C_6\|\theta_0\|_{L^\infty}
-\frac{\kappa^2C^2_3C^2_5}{C^2_1}.
\end{split}
\end{eqnarray*}
which implies that
\begin{eqnarray}\label{differential-J-derivative}
J^\prime(t)
\geq C_7(n,\alpha)(J(t))^2-C_8(\kappa,n,\alpha,\gamma)(1+\|\theta_0\|_{L^\infty})^2.
\end{eqnarray}
Now, choosing $A>\sqrt{\frac{C_8}{C_7}}$ and considering $J(0)\geq A(1+\|\theta_0\|_{L^\infty})$, \eqref{differential-J-derivative} implies $J^\prime(0)>0$ and that $J$ blows up in finite time. Since
\begin{eqnarray*}
\begin{split}
|J(t)|
%&\leq& \int_{|x|\leq1}\frac{|\theta(0,t)-\theta(x,t)|}{|x|^{n}}dx
%+2\|\theta_0\|_{L^\infty}\int_{|x|>1}\frac{e^{-|x|}}{|x|^{n}}dx\\
&\leq\|\nabla\theta(t)\|_{L^\infty}\int_{|x|\leq1}\frac{dx}{|x|^{n-1}}
+2\|\theta_0\|_{L^\infty}\int_{|x|>1}\frac{e^{-|x|}}{|x|^{n}}dx\\
&\leq C(n)(\|\nabla\theta\|_{L^\infty}+\|\theta_0\|_{L^\infty}),
\end{split}
\end{eqnarray*}
we conclude that $\|\nabla\theta(t)\|_{L^\infty}$ must blow up in finite time. The proof of Theorem \ref{the-3} is finished.
\appendix
\section{Proof of Lemma \ref{gradient-velocity}}
In this appendix, we give the proof of Lemma \ref{gradient-velocity}.

\textbf{Proof of Lemma \ref{gradient-velocity}}.
%Without loss of generality, we only consider the case of whole space.
By the integral representation of the Riesz potential \cite{[Stein]}, we have
\begin{eqnarray}\label{u-representation-A.1}
u(x,t)=\frac{\Gamma(\frac{n}{2}-1+\alpha)}{\pi^{\frac{n}{2}}2^{2-2\alpha}\Gamma(1-\alpha)}\int_{\R^{n}}\frac{\nabla\theta(y)}{|x-y|^{n-2+2\alpha}}dy.
\end{eqnarray}
Following \cite{[Beale-Kato-Majda]}, we introduce a cut-off function $\zeta_\rho(x)$ satisfying
\begin{eqnarray*}
\zeta_\rho(x)
=
\begin{cases}
1,
&  \mbox{{\rm for} $|x|<\rho,$ }\\
0,
& \mbox{{\rm for} $|x|>2\rho,$ } \\
\end{cases}
~{\rm and}~ \nabla\zeta_\rho(x)\leq C/\rho.
\end{eqnarray*}
Here $\rho\leq1$ is a radius to be chosen suitably small later.
By introducing a factor $\zeta_\rho(x-y)+[1-\zeta_\rho(x-y)]$ under the integral sign in \eqref{u-representation-A.1}, we can compute $\nabla u(x)$ as
\begin{eqnarray}\label{A.1}
\nabla u(x)=\underbrace{\int_{\R^{n}}\frac{\zeta_\rho(x-y)}{|x-y|^{n-2+2\alpha}}\nabla^2\theta(y)dy}_{\nabla u^{(1)}(x)}
+\underbrace{\int_{\R^{n}}\nabla_x\Big[\frac{1-\zeta_\rho(x-y)}{|x-y|^{n-2+2\alpha}}\Big]\nabla\theta(y)dy}_{\nabla u^{(2)}(x)}.
\end{eqnarray}
%When $s>\frac{n}{2}+2$, by H\"{o}lder's inequality, we can estimate $\nabla u^{(1)}(x)$ as
%\begin{eqnarray*}
%|\nabla u^{(1)}(x)|
%&\leq&\|\zeta_\rho\|_{L^\infty}\int_{|y-x|\leq2\rho}\frac{|\nabla^2\theta(y)|}{|x-y|^{n-2+2\alpha}}dy\\
%&\leq&\|\nabla^2\theta\|_{L^\infty}\int_{|z|<2\rho}\frac{dz}{|z|^{n-2+2\alpha}}\\
%&\leq&C_{n,\alpha}\rho^{2-2\alpha}\|\nabla^2\theta\|_{H^{s-2}}\\
%&\leq& C_{n,\alpha}\rho^{2-2\alpha}\|\theta\|_{H^s}.
%\end{eqnarray*}
%When $\frac{n}{2}+1<s<\frac{n}{2}+2$, by Sobolev's inequality, we also have that
%\begin{eqnarray*}
%|\nabla u^{(1)}(x)|
%&\leq&\|\zeta_\rho\|_{L^\infty}\int_{|y-x|\leq2\rho}\frac{|\nabla^2\theta(y)|}{|x-y|^{n-2+2\alpha}}dy\\
%&\leq&\Big(\int_{|z|<2\rho}\frac{dz}{|z|^{\frac{2n(n-2+2\alpha)}{n+2(s-2)}}}\Big)^{\frac{n+2(s-2)}{2n}}\|\nabla^2\theta\|_{L^{\frac{2n}{n-2(s-2)}}}\\
%&=&C_{n,\alpha,p}\rho^{s-\frac{n}{2}-2\alpha}\|\nabla^2\theta\|_{L^{\frac{2n}{n-2(s-2)}}}\\
%&\leq&C_{n,\alpha,p}\rho^{s-\frac{n}{2}-2\alpha}\|\nabla^2\theta\|_{H^{s-2}}\\
%&\leq& C_{n,\alpha,p}\rho^{s-\frac{n}{2}-2\alpha}\|\theta\|_{H^s}.
%\end{eqnarray*}
For $\alpha\in(0,\frac12]$ and $s\neq\frac{n}{2}+2$, by H\"{o}lder's inequality and the Sobolev embedding, we can estimate $\nabla u^{(1)}(x)$ as
\begin{eqnarray}\label{A.2}
\begin{split}
|\nabla u^{(1)}(x)|
%&\leq&\|\zeta_\rho\|_{L^\infty}\int_{|y-x|\leq2\rho}\frac{|\nabla^2\theta(y)|}{|x-y|^{n-2+2\alpha}}dy\nonumber\\
&\leq
\begin{cases}
\|\nabla^2\theta\|_{L^\infty}\int\limits_{|z|<2\rho}\frac{dz}{|z|^{n-2+2\alpha}},
&  \mbox{for $s>\frac{n}{2}+2,$ }\\
\Big[\int\limits_{|z|<2\rho}\frac{dz}{|z|^{\frac{2n(n-2+2\alpha)}{n+2(s-2)}}}\Big]^{\frac12+\frac{s-2}{n}}\|\nabla^2\theta\|_{L^{\frac{2n}{n-2(s-2)}}},
& \mbox{for $s\in(\frac{n}{2}+1,\frac{n}{2}+2),$ } \\
\end{cases}\\
%&=&
%\begin{cases}
%C_{n,\alpha}\rho^{2-2\alpha}\|\nabla^2\theta\|_{L^\infty},
%&  \mbox{for $s>\frac{n}{2}+2,,$ }\\
%C_{n,\alpha,s}\rho^{s-\frac{n}{2}-2\alpha}\|\nabla^2\theta\|_{L^{\frac{2n}{n-2(s-2)}}},
%& \mbox{for $\frac{n}{2}+1<s<\frac{n}{2}+2,$ } \\
%\end{cases}\\
%&\leq&
%\begin{cases}
%C_{n,\alpha}\rho^{2-2\alpha}\|\nabla^2\theta\|_{H^{s-2}},
%&  \mbox{for $s>\frac{n}{2}+2,$ }\\
%C_{n,\alpha,s}\rho^{s-\frac{n}{2}-2\alpha}\|\nabla^2\theta\|_{H^{s-2}},
%& \mbox{for $\frac{n}{2}+1<s<\frac{n}{2}+2,$ } \\
%\end{cases}\nonumber\\
&\leq
\begin{cases}
C_{n,\alpha}\rho^{2-2\alpha}\|\theta\|_{H^s},
&  \mbox{for $s>\frac{n}{2}+2,$ }\\
C_{n,\alpha,s}\rho^{s-\frac{n}{2}-2\alpha}\|\theta\|_{H^s},
& \mbox{for $s\in(\frac{n}{2}+1,\frac{n}{2}+2).$ } \\
\end{cases}
\end{split}
\end{eqnarray}
For the remaining case $s=\frac n2+2$ and $\alpha=\frac12$, by H\"{o}lder's inequality and the embedding $H^{\frac{n}{2}}(\R^n)\hookrightarrow L^q(\R^n)$ for any $q\in[2,\infty)$, we can derive that
\begin{eqnarray}\label{A.3}
|\nabla u^{(1)}(x)|
%&\leq&\|\zeta_\rho\|_{L^\infty}\int_{|y-x|\leq2\rho}\frac{|\nabla^2\theta(y)|}{|x-y|^{n-1}}dy\nonumber\\
\lesssim\Big[\int_{|z|<2\rho}\frac{dz}{|z|^{\frac{(n+1)(n-1)}{n}}}\Big]^{\frac{n}{n+1}}\|\nabla^2\theta\|_{L^{n+1}}
%&=&C_{n,\alpha,p}\rho^{\frac{n}{p}-(n-1)}\|\nabla^2\theta\|_{L^{\frac{p}{p-1}}}\\
\leq C_{n}\rho^{\frac{1}{n+1}}\|\theta\|_{H^{\frac{n}{2}+2}}.
\end{eqnarray}
Similarly, for $s=\frac n2+2$ and $\alpha\in(0,\frac12)$, the term $\nabla u^{(1)}(x)$ can be estimated as
\begin{eqnarray}\label{A.4}
|\nabla u^{(1)}(x)|
\leq\Big[\int_{|z|<2\rho}\frac{dz}{|z|^{\frac{n(n-2+2\alpha)}{n-1}}}\Big]^{\frac{n-1}{n}}\|\nabla^2\theta\|_{L^n}
\leq C_{n,\alpha}\rho^{1-2\alpha}\|\theta\|_{H^{\frac{n}{2}+2}}.
\end{eqnarray}
Letting $\varepsilon>0$ depending only on $s,n,\alpha$ given by, for $n\geq2$,
\begin{eqnarray*}
\varepsilon
=
\begin{cases}
2-2\alpha,
&  \mbox{for $s>\frac{n}{2}+2$ and $\alpha\in(0,\frac12]$, }\\
s-\frac{n}{2}-2\alpha,
& \mbox{for $s\in(\frac{n}{2}+1,\frac{n}{2}+2)$ and $\alpha\in(0,\frac12]$, }\\
\frac{1}{n+1},
& \mbox{for $s=\frac{n}{2}+2$ and $\alpha=\frac12$,} \\
1-2\alpha,
& \mbox{for $s=\frac{n}{2}+2$ and $\alpha\in(0,\frac12)$, } \\
\end{cases}
\end{eqnarray*}
and gathering up \eqref{A.2}-\eqref{A.4} yield that, for $\alpha\in(0,\frac12]$ and $s>\frac{n}{2}+1$ with $n\geq2$,
\begin{eqnarray}\label{A.5}
|\nabla u^{(1)}(x)|
\leq C\rho^{\varepsilon}\|\theta\|_{H^s}.
\end{eqnarray}
We proceed to split $\nabla u^{(2)}(x)$ into two parts as
\begin{eqnarray}\label{A.6}
\nabla u^{(2)}(x)
%&=&\int_{|x-y|\geq\rho}\nabla_x\Big\{\frac{1-\zeta_\rho(x-y)}{|x-y|^{n-2+2\alpha}}\Big\}\nabla^\bot\theta(y)dy\nonumber\\
=\Big[\underbrace{\int_{\rho\leq|x-y|\leq1}}_{\nabla u^{(21)}(x)}+\underbrace{\int_{|x-y|>1}}_{\nabla u^{(22)}(x)}\Big]\nabla_x\Big[\frac{1-\zeta_\rho(x-y)}{|x-y|^{n-2+2\alpha}}\Big]\nabla\theta(y)dy.
\end{eqnarray}
For the term $\nabla u^{(21)}(x)$, it can be bounded as, for $\alpha\in(0,\frac12]$ and $\rho\in(0,1]$,
\begin{eqnarray}\label{A.7}
\begin{split}
|\nabla u^{(21)}(x)|
&\leq C\Big[\int^1_\rho \frac{r^{n-1}dr}{r^{n-1+2\alpha}}+\rho^{-1}\int^{2\rho}_\rho \frac{r^{n-1}dr}{r^{n-2+2\alpha}}\Big]\|\nabla\theta\|_{L^\infty}\\
&\leq
\begin{cases}
C(-\log\rho+1)\|\nabla\theta\|_{L^\infty},
&  \mbox{for $\alpha=\frac12,$ }\\
C\Big[\frac{1}{1-2\alpha}+\frac{2^{2-2\alpha}-1}{2-2\alpha}\Big]\|\nabla\theta\|_{L^\infty},
& \mbox{for $\alpha\in(0,\frac12).$ } \\
\end{cases}
\end{split}
\end{eqnarray}
When $\rho\leq\frac12$, by the choice of $\zeta$ and H\"{o}lder inequality, we have, for $\alpha\in(0,\frac12]$,
\begin{eqnarray}\label{A.8}
|\nabla u^{(22)}(x)|
%&\leq&\int_{|x-y|>1}\frac{|\nabla\theta(y)|}{|x-y|^{n-1+2\alpha}}dy\nonumber\\
%&\leq&\int_{|x-y|>1}\frac{|\nabla^\bot\theta(y)|}{|x-y|^{n-1+2\alpha}}dy\nonumber\\
\lesssim\Big[\int_{|z|>1}\frac{dz}{|z|^{2(n-1+2\alpha)}}\Big]^{\frac12}\|\nabla\theta\|_{L^2}
\leq C_{n,\alpha}\|\nabla\theta\|_{L^2}.
\end{eqnarray}
For $\rho\in(\frac12,1]$, the term $\nabla u^{(22)}(x)$ can be estimated as
\begin{eqnarray*}
&&|\nabla u^{(22)}(x)|
=\Big|\Big[\int_{1<|x-y|\leq2\rho}+\int_{|x-y|>2\rho}\Big]\nabla_x\Big[\frac{1-\zeta_\rho(x-y)}{|x-y|^{n-2+2\alpha}}\Big]\nabla^\bot\theta(y)dy\Big|\nonumber\\
%&\leq&C\|\nabla\theta\|_{L^\infty}\Big(\int^{2\rho}_1\frac{r^{n-1}dr}{r^{n-1+2\alpha}}+\rho^{-1}\int^{2\rho}_1\frac{r^{n-1}dr}{r^{n-2+2\alpha}}\Big)
%+\int\limits_{|x-y|>2\rho}\frac{|\nabla^\bot\theta(y)|dy}{|x-y|^{n-1+2\alpha}}\nonumber\\
&&
\lesssim\|\nabla\theta\|_{L^\infty}\Big[\int^{2\rho}_1\frac{r^{n-1}dr}{r^{n-1+2\alpha}}+\frac{1}{\rho}\int^{2\rho}_1\frac{r^{n-1}dr}{r^{n-2+2\alpha}}\Big]
+\|\nabla\theta\|_{L^2}\Big[\int^\infty_{2\rho}\frac{r^{n-1}dr}{r^{2(n-1+2\alpha)}}\Big]^{\frac12}\nonumber\\
%&\leq&
%\begin{cases}
%C_n\Big(\ln(2\rho)+2-\frac{1}{\rho}\Big)\|\nabla\theta\|_{L^\infty}
%+C_n\Big(\int\limits_{|z|>2\rho}\frac{dz}{|z|^{2n}}\Big)^{\frac12}\|\nabla\theta\|_{L^2},
%&  \mbox{for $\alpha=\frac12,$ }\nonumber\\
%C\cdot\frac{3(2^{1-2\alpha}-1)}{1-2\alpha}\cdot\|\nabla\theta\|_{L^\infty}
%+\|\nabla\theta\|_{L^2}\Big(\int\limits_{|z|>2\rho}\frac{dz}{|z|^{2(n-1+2\alpha)}}\Big)^{\frac12},
%& \mbox{for $\alpha\in(0,\frac12),$ }\\
%\end{cases}\\
%&\leq&C\|\nabla\theta\|_{L^\infty}\Big(\int^{2}_1\frac{dr}{r^{2\alpha}}+2\int^{2}_1r^{1-2\alpha}dr\Big)
%+\|\nabla\theta\|_{L^2}\Big(\int^\infty_{1}\frac{dr}{r^{n-1+4\alpha}}\Big)^{\frac12}\nonumber\\
%&\leq&
%\begin{cases}
%C_n\Big(\ln2+2\Big)\|\nabla\theta\|_{L^\infty}
%+\frac{C_n\|\nabla\theta\|_{L^2}}{(2\rho)^{\frac n2}},
%&  \mbox{for $\alpha=\frac12,$ }\\
%\frac{3 C(2^{1-2\alpha}-1)}{1-2\alpha}\cdot\|\nabla\theta\|_{L^\infty}
%+\frac{C\|\nabla\theta\|_{L^2}}{(2\rho)^{\frac{n}{2}-1+2\alpha}},
%& \mbox{for $\alpha\in(0,\frac12),$ }\\
%\end{cases}\nonumber\\
&&
\leq C_{n,\alpha}(\|\nabla\theta\|_{L^\infty}+\|\nabla\theta\|_{L^2}),
\end{eqnarray*}
which along with \eqref{A.8} yields that for any $\rho\in(0,1]$,
%Gathering up \eqref{} and \eqref{} yields that, for any $0<\rho\leq1$,
\begin{eqnarray*}
|\nabla u^{(22)}(x)|\leq C_{n,\alpha}(\|\nabla\theta\|_{L^\infty}+\|\nabla\theta\|_{L^2}).
\end{eqnarray*}
It follows from \eqref{A.7} and \eqref{A.6} that, for any $\rho\in(0,1]$,
\begin{eqnarray*}
|\nabla u^{(2)}(x)|
%&\leq&|\nabla u^{(21)}(x)|+|\nabla u^{(22)}(x)|\\
\leq
\begin{cases}
C(-\log\rho+1)\|\nabla\theta\|_{L^\infty}+C\|\nabla\theta\|_{L^2},
&  \mbox{for $\alpha=\frac12,$ }\\
C(\|\nabla\theta\|_{L^\infty}+\|\nabla\theta\|_{L^2}),
& \mbox{for $\alpha\in(0,\frac12).$ } \\
\end{cases}
\\
\end{eqnarray*}
Combining this estimate and \eqref{A.5} implies that, for $s>\frac{n}{2}+1$ and any $\rho\in(0,1]$,
\begin{eqnarray}\label{A.9}
\|\nabla u\|_{L^\infty}
\leq
\begin{cases}
C\rho^{\varepsilon}\|\theta\|_{H^s}+C(-\log\rho+1)\|\nabla\theta\|_{L^\infty}+C\|\nabla\theta\|_{L^2},
&  \mbox{for $\alpha=\frac12,$ }\\
C\rho^{\varepsilon}\|\theta\|_{H^s}+C(\|\nabla\theta\|_{L^\infty}+\|\nabla\theta\|_{L^2}),
& \mbox{for $\alpha\in(0,\frac12).$ } \\
\end{cases}
\end{eqnarray}
If $\|\theta\|_{H^s}\leq1$, we take $\rho=1$; otherwise we choose $\rho$
such that $\rho^{\varepsilon}\|\theta\|_{H^s}=1$, namely, $\rho=\|\theta\|^{-\frac{1}{\epsilon}}_{H^s}$,
%\begin{eqnarray*}
%\rho^{\varepsilon}\|\theta\|_{H^s}=1,{\rm namely},~\rho=\|\theta\|^{-\frac{1}{\epsilon(s)}}_{H^s}
%\end{eqnarray*}
%$\rho^{\varepsilon}\|\theta\|_{H^s}=1$, namely,
%$\rho=\|\theta\|^{-\frac{1}{\epsilon(s)}}_{H^s}$,
and \eqref{A.9} becomes \eqref{gradient-velocity-3.10}. In either case \eqref{gradient-velocity-3.10} holds.

With only minor changes, the same proof applies to periodic fluid
flow. We then complete the proof of Lemma \ref{gradient-velocity}.
%\begin{eqnarray*}
%\|\nabla u\|_{L^\infty}
%\leq
%\begin{cases}
%C\Big[1+(1+\ln^+\|\theta\|_{H^s})\|\nabla\theta\|_{L^\infty}+\|\nabla\theta\|_{L^2}\Big],
%&  \mbox{for $\alpha=\frac12,$ }\\
%C(1+\|\nabla\theta\|_{L^\infty}+\|\nabla\theta\|_{L^2}),
%& \mbox{for $\alpha\in(0,\frac12),$ } \\
%\end{cases}
%\end{eqnarray*}

{\bf Acknowledgements.} W. Zhang was partially supported
by the Science and Technology Research Project of Department of Education of Jiangxi
Province (No. GJJ2200363).
%%reference

\end{document}